\documentclass[12pt,a4paper]{article}
\usepackage[a4paper,top=2cm,bottom=2cm,left=2cm,right=2cm,marginparwidth=1.75cm]{geometry}
\usepackage{indentfirst}
\usepackage{natbib, authblk}
\usepackage{hyperref}
\hypersetup{
colorlinks =true,
citecolor  = blue,
linkcolor = blue,
urlcolor = black
}
\usepackage{xcolor}

\usepackage{hyperref}
\hypersetup{
colorlinks =true,
citecolor  = blue,
linkcolor = blue,
urlcolor = black
}
\usepackage{amsmath, bbm, dsfont, amsfonts, xcolor, graphicx, bm}
\usepackage{enumitem}
\usepackage{amsthm, natbib, comment}
\usepackage{etoc}

\newcommand{\quotes}[1]{``#1''}
\newcommand{\indep}{\perp \!\!\! \perp}

\newenvironment{funding}{\section*{Funding}}{}

\usepackage{fancyhdr}
\title{Asymptotic Normality of the Pearson Correlation between Two Weakly Dependent Random Fields}
\author[1*]{Jiyue Qin}
\author[1,2]{Armin Schwartzman}
\author[1]{Samuel Davenport}

\affil[1]{\small Division of Biostatistics and Bioinformatics, Herbert Wertheim School of Public Health and Human Longevity Science, University of California, San Diego}
\affil[2]{\small Halıcıoğlu Data Science Institute, University of California, San Diego}

\affil[*]{Corresponding author: \texttt{j5qin@ucsd.edu}}

\date{}

\begin{document}
\maketitle

\begin{abstract}
We study the asymptotic behavior of the sample Pearson correlation coefficient between two weakly dependent random fields on a lattice of any fixed dimension. Under strict stationarity, $(BL,\theta)$-weak dependence and mild moment conditions, we establish the consistency and asymptotic normality of the estimator on regularly growing domains, which need not be rectangular, with an asymptotic variance expressed through joint fourth-order moments. When the correlation is zero, asymptotic normality holds under weaker conditions. When the fields are independent, the asymptotic variance simplifies to a sum over spatial lags of products of their autocorrelations. The variance approximations behind the classic modified $t$-test, although derived for independent Gaussian fields, target the right variance even without Gaussianity. We also give sufficient conditions for asymptotic normality that are easy to check and cover bounded, Gaussian, $m$-dependent and independent fields. For Gaussian fields, summable auto- and cross-covariances suffice. These conditions rely on new closure results. Squares and products, although not Lipschitz, preserve $(BL,\theta)$-dependence for Gaussian fields, and products of independent fields do so when one is positively or negatively associated. Finally, we propose kernel-based estimators of the asymptotic variance, enabling inference without Gaussianity or parametric spatial models. Simulation studies illustrate their finite-sample performance.
\end{abstract}


\theoremstyle{plain}
\newtheorem{theorem}{Theorem}
\newtheorem*{theorem*}{Theorem}
\newtheorem{lemma}{Lemma}
\newtheorem{assumption}{Assumption}
\newtheorem{corollary}{Corollary}
\newtheorem{proposition}{Proposition}
\newtheorem{conditionB}{Condition}
\renewcommand{\theconditionB}{B\arabic{conditionB}}

\newtheorem{conditionBs}{Condition}
\renewcommand{\theconditionBs}{B\arabic{conditionBs}\textsuperscript{*}}

\newtheorem{conditionM}{Condition}
\renewcommand{\theconditionM}{M\arabic{conditionM}}

\newtheorem{conditionG}{Condition}
\renewcommand{\theconditionG}{G\arabic{conditionG}}
\newtheorem{assumptionStar}{Assumption}
\renewcommand{\theassumptionStar}{\arabic{assumption}\textsuperscript{*}}

\theoremstyle{definition}
\newtheorem{definition}{Definition}
\newtheorem{remark}{Remark}
\newtheorem{example}{Example}
\theoremstyle{plain}

\newcommand{\cov}{\operatorname{Cov}}

\newcommand{\rfgeneral}{\texttt{RFCorr}}
\newcommand{\rfgaussian}{\texttt{RFCorr-\allowbreak Gauss}}
\newcommand{\rfzero}{\texttt{RFCorr-\allowbreak 0}}
\newcommand{\rfzerogaussian}{\texttt{RFCorr-\allowbreak Gauss-\allowbreak 0}}
\newcommand{\rfindep}{\texttt{RFCorr-\allowbreak Indep}}

\section{Introduction}

The correlation between two variables observed over time or space serves as a fundamental
measure of global association and is of broad interest across diverse scientific domains.
We model such variables as real-valued random fields on the lattice $\mathbb{Z}^d$ for some integer
$d\geq1$.
For $d=1$, typical examples include correlation between two financial time series such as stock
log-returns \citep{campbell1998econometrics}, between two genotype profiles along chromosomes used to
assess genetic relatedness between individuals \citep{weir1996genetic}, and between two functional
Magnetic Resonance Imaging (fMRI) time series
used to construct functional brain networks  \citep{fox2005human}.
For $d=2$, examples include correlation between environmental variables across geographic locations
such as biodiversity and remoteness \citep{viladomat2014assessing}, concentrations of the heavy
metals arsenic and lead \citep{vallejos2020spatial}, and correlation between two digital images
\citep{martens1998image}.
For $d=3$, an example is correlation between brain activation maps in neuroimaging such as stress
and fear \citep{alexander2018testing}.

The most commonly used estimator for correlation is the Pearson correlation coefficient $\hat\rho$,
due to its simplicity, interpretability, and invariance to affine transformations (see e.g.,
\citealt{afyouni2019effective, viladomat2014assessing}).
For $N$ independent and identically distributed (i.i.d.) bivariate Gaussian observations, the
classical test of zero correlation (equivalently, of independence) is based on $\hat\rho\sqrt{N-2}/\sqrt{1-\hat\rho^2}$, which follows
a $t_{N-2}$ distribution under the null.
For random fields, however, observations are typically autocorrelated.
A common strategy for approximate inference in this setting modifies results under no
autocorrelation by replacing sample size with a typically deflated effective sample size (ESS).
Various formulas of the ESS have been proposed, including \cite{bartlett1935some},
\cite{quenouille1947notes},
\cite{medrano2024assessing} for $d=1$ and
\citet{clifford1989assessing}, \citet{dutilleul1993modifying} for $d\geq1$.
However, these methods were derived under the null
hypothesis that the two series or fields are independent, typically with further strong assumptions
such as Gaussianity or parametric autoregressive structures.
Nonparametric alternatives such as the permutation procedures of \citet{viladomat2014assessing} and
\citet{alexander2018testing} also test the null of independence.
They preserve spatial structure empirically but lack theoretical guarantees, can have inflated
false positive rates \citep{markello2021comparing} and can be computationally intensive.

While much of the existing literature focuses on testing the null hypothesis of independence,
hypothesis testing alone is often insufficient for scientific inference.
In many applications, the primary goal is not merely to detect whether correlation is nonzero, but
to estimate its magnitude and quantify uncertainty, for example through confidence intervals.
Consequently, it is essential to understand the large-sample behavior of the Pearson correlation, in
particular, its consistency and asymptotic normality (with an expression for asymptotic variance).
For i.i.d.\ observations, the large-sample behavior of the Pearson correlation is classical: it is
consistent and asymptotically normal (CAN) when the fourth moments are finite, with an
asymptotic variance that depends on moments up to fourth order
\citep[Example~3.6]{van2000asymptotic}.
For bivariate Gaussian observations, this variance reduces to $(1-\rho^2)^2$
\citep[Theorem~4.2.4]{anderson2003introduction}.

In contrast, when observations exhibit temporal or spatial dependence, the validity of these
properties is far less straightforward.
In the time series literature, building on \cite{hannan1976asymptotic}, \citet{roy1989asymptotic}
established the consistency and asymptotic normality of Pearson correlation for an ergodic, strictly
stationary, linearly purely nondeterministic, near-Gaussian time series.
While mathematically rigorous, these assumptions are difficult to interpret or verify in
applications, and they do not readily generalize to spatially indexed data.
For lattices of dimension $d\geq1$, \citet{bhattacharya2025nonsense} showed that the Pearson
correlation between two independent Ising models is asymptotically normal, without giving the
asymptotic variance explicitly.
Like the effective sample size methods above, their results are restricted to independent fields.
For random fields, \citet{li2008asymptotic} established the
asymptotic normality of sample autocovariances of a univariate field under strong mixing and a moment condition imposed on the
estimators themselves, and \citet{li2008testing} obtained the analogous result for cross-covariances
of multivariate fields, but only for fields observed at fixed locations over a growing time domain.
Neither result covers the variances and cross-covariance of two fields jointly over a growing
spatial domain, which the sample correlation requires.
Moreover, their strong mixing conditions are hard to verify, and for $d>1$ few examples of mixing random fields are known beyond $m$-dependent
ones \citep[Chapter 1, Section 5]{bulinski2007limit}.
To the best of our knowledge, no prior work has established the consistency and asymptotic
normality of $\hat\rho$ for general $\rho$ in random fields with $d\geq1$ under weak-dependence
conditions that are easy to verify.

The main contributions of this paper are as follows.
First, we fill the literature gap by establishing the asymptotic theory for the Pearson correlation
between two random fields under $(BL,\theta)$-weak dependence.
The $(BL,\theta)$ framework, introduced by \citet{doukhan1999new} for time series and extended to
random fields by \citet{bulinski2001normal} (see also \citealp{bulinski2007limit}), is formulated using bounded Lipschitz functionals and
provides a unified treatment of weak dependence.
Under mild moment and dependence conditions, we show that Pearson correlation is CAN, and derive its asymptotic variance, expressed in terms of joint fourth-order moments of the
standardized fields, which reduce to auto- and cross-correlations in the Gaussian case.
When $\rho=0$, we obtain asymptotic normality under weaker conditions, which are the ones relevant
for testing.

Second, when $X$ and $Y$ are independent, the asymptotic variance reduces to a sum of products of
the autocorrelations of $X$ and $Y$.
We relate this variance to the approximations used in the modified $t$-test of
\citet{clifford1989assessing} and its refinement by \citet{dutilleul1993modifying}, which were
both derived for independent Gaussian fields.
These approximations are expressed through the true covariance matrices of $X$ and $Y$ over the
observed grid points, which are estimated by sample counterparts in practice.
We show that both approximations, evaluated at the true covariance matrices, converge to our
asymptotic variance in any fixed dimension $d$, even
without Gaussianity.
Thus, although these approximations were derived under Gaussianity, they target the correct
variance for non-Gaussian independent fields as well.
This also proves the large-sample agreement between the two approximations that
\citet{dutilleul1993modifying} observed numerically.

Third, our CLT for general $\rho$ requires a weak dependence condition on the
squares and cross-products of the standardized fields.
We establish easy-to-verify conditions on the original fields $X$ and $Y$ under
which this condition holds, for four broad classes.
For Gaussian fields, it suffices that the auto- and cross-covariances of $X$ and $Y$ are
summable.
For $m$-dependent fields, no further condition is needed.
For bounded fields, it suffices that $(X,Y)$ is $(BL,\theta)$, for example when it is
positively or negatively associated (Definitions~\ref{def:pa} and \ref{def:na}) with
summable auto- and cross-covariances.
For independent $X$ and $Y$ that need not be bounded or Gaussian, it suffices that $X$ and $Y$ are each positively or negatively associated with
summable autocovariances.
Summability holds, for example, whenever the covariances decay exponentially.
To obtain these conditions, we establish closure results for $(BL,\theta)$-dependence. We show
that it is preserved under
Lipschitz maps, which covers squares and products of bounded fields. Without boundedness these
maps are not Lipschitz, but they still preserve $(BL,\theta)$-dependence for Gaussian fields, and
the product
of two independent $(BL,\theta)$ fields is again $(BL,\theta)$ when one of them is positively or
negatively associated.

{Finally, we propose kernel-based estimators of the asymptotic variance, built on a spatial
Newey--West estimator with a product Bartlett kernel.
A single general construction yields five estimators.
Two estimate the variance at an arbitrary $\rho$ for confidence intervals, one without
distributional assumptions and one for Gaussian fields.
The other three impose a null hypothesis and are used for testing. Two impose $\rho=0$, again one
without distributional assumptions and one for Gaussian fields, and the third imposes independence
of $X$ and $Y$, under which the asymptotic variance takes the same form whether or not the fields
are Gaussian.
All of these variance estimators apply in any fixed dimension $d$ and on domains that need not be
rectangular, and a simple modification
guarantees nonnegative variance estimates.
The truncation lag grows as the $(d+2)$-th root of the number of observed grid points.
This rate minimizes the asymptotic truncated mean squared error of the Bartlett estimator
\citep{kim2011spatial} and reduces to the classical cube-root rule when $d=1$ \citep{andrews1991heteroskedasticity}.}

Together, these results yield a theoretically justified and computationally feasible alternative to
heuristic effective-sample-size adjustments, enabling asymptotically valid inference for
correlations between weakly dependent random fields without requiring Gaussianity or a specific
parametric model.

The rest of the paper is organized as follows.
Section \ref{sec:pre} introduces the notation, definitions and assumptions in the
$(BL,\theta)$-dependence framework.
Section \ref{sec:theorem} presents the main asymptotic results.
Section \ref{sec:condition} presents sufficient conditions for asymptotic normality.
Section \ref{sec:var} introduces the kernel-based estimator of the asymptotic variance.
Section \ref{sec:simu} reports numerical results illustrating finite-sample performance.
{For both Gaussian and non-Gaussian fields, the confidence intervals and hypothesis tests
based on the proposed variance estimators attain nominal levels.
In contrast, ignoring autocorrelation, or using Gaussian variance formulas for non-Gaussian fields,
causes undercoverage or inflated Type~I error that does not vanish as the domain grows.
The exception is testing for independent fields, where, as our theory shows, the Gaussian formula
remains valid without Gaussianity.}
Section \ref{sec:discussion} concludes with a discussion and possible extensions.

\section{Preliminaries}\label{sec:pre}

\subsection{Notation}
Let $X=\{X_s\in \mathbb{R}:s\in \mathbb{Z}^d\}$, $Y = \{Y_s\in \mathbb{R}:s\in \mathbb{Z}^d\}$, $d
\in \mathbb{N}$ be two real-valued random fields.
Let $Z = \{\bm{Z}_s=[X_s, Y_s]'\in \mathbb{R}^2:s\in \mathbb{Z}^d\}$ be the bivariate random field.
Let $U_n$ be a sequence of sets in $\mathbb{Z}^d$ and $\{\bm{Z}_s:s\in U_n\subset \mathbb{Z}^d\}$ be
a realization of size $|U_n|$, where $|U_n|$ denotes the cardinality of a finite set $U_n$.
We assume that $Z$ is strictly stationary and has finite fourth moments.
Then the mean and variance of $X_s$ are constant for all $s$, denoted as $\mu_x < \infty$,
$\sigma^2_x < \infty$, respectively.
Similarly, the mean and variance of $Y_s$ are denoted as $\mu_y< \infty$, $\sigma^2_y<\infty$.
{We also assume $\sigma^2_x>0$ and $\sigma^2_y>0$.}
The correlation between $X_s$ and $X_{s+k}$ only depends on the lag $k\in \mathbb{Z}^d$ and is denoted as
$\rho_{xx, k}$.
Similarly, the correlation between $Y_s$ and $Y_{s+k}$ is denoted as $\rho_{yy,k}$, the
correlation between $X_s$ and $Y_{s+k}$ as $\rho_{xy, k}$, and the correlation between $Y_s$ and
$X_{s+k}$ as $\rho_{yx, k}$.
Note that $\rho_{xx,k} = \rho_{xx,-k}$, $\rho_{yy,k} = \rho_{yy,-k}$, $\rho_{xy,k} = \rho_{yx,-k}$.
Our parameter of interest is $\rho \equiv \rho_{xy,0}$.
Let $\bar{X}_n =|U_n|^{-1}\sum_{s \in U_n} X_s $ and $\bar{Y}_n =|U_n|^{-1}\sum_{s \in U_n} Y_s$
denote the sample means of $X$ and $Y$, respectively.
Our estimator of interest, the Pearson correlation coefficient, is defined as:
$$\hat \rho =\frac{\sum_{s \in U_n} (X_s-\bar{X}_n) (Y_s-\bar{Y}_n)}{\sqrt{\sum_{s \in U_n}
(X_s-\bar{X}_n)^2 \sum_{s \in U_n}(Y_s-\bar{Y}_n)^2}}.$$
Let $\tilde{X}_s$ and $\tilde{Y}_s$ be centered and scaled versions of $X_s$ and $Y_s$ to have mean
0 and variance 1, i.e., $\tilde{X}_s = (X_s-\mu_x)/\sigma_x$ and $\tilde{Y}_s =
(Y_s-\mu_y)/\sigma_y$.
Let $\tilde X=\{\tilde X_s\in \mathbb{R}:s\in \mathbb{Z}^d\}$, $\tilde Y = \{\tilde Y_s\in \mathbb{R}:s\in
\mathbb{Z}^d\}$.

Throughout the paper, vectors are denoted in boldface.
For $\bm{v}=[v_1, \dots, v_d]'\in \mathbb{R}^d$, we use $\|\bm{v}\|_p$ to denote its $\ell^p$-norm,
defined as $\|\bm{v}\|_p = (\sum_{i=1}^d |v_i|^p)^{1/p}, 1\leq p \in \mathbb{R}$.
We use $\|\bm{v}\|_\infty$ to denote its sup-norm, defined as $\|\bm{v}\|_\infty= \max_{1\le i\le d}|v_i|$.
For a function $f:\mathbb{R}^d\to \mathbb{R}$, we use $\|f\|_\infty$ to denote its sup-norm, defined
as $\|f\|_\infty = \sup_x|f(x)|$.
If $f$ is Lipschitz continuous, we define the Lipschitz constant as $\operatorname{Lip}(f)=\sup
\{|f(\bm{x})- f(\bm{y})|/{\|\bm{x}- \bm{y}\|_1}: \bm{x}, \bm{y} \in \mathbb{R}^d , \bm{x} \neq
\bm{y}\}.$
We use $\cov(\bm{A},\bm{B})$ to denote the covariance between random vectors $\bm{A}$ and $\bm{B}$,
and $\operatorname{Var}(\bm{A})$ to denote the variance-covariance matrix of $\bm{A}$.
Convergence in distribution to a Gaussian random variable with mean 0 and variance $\sigma^2$ is
denoted by $\xrightarrow{d} \mathcal{N}(0, \sigma^2)$.

\subsection{Weak law of large numbers (WLLN) for the sample mean}

We state and prove two convergence results for sample means of random fields on $\mathbb{Z}^d$.
The first gives the WLLN under vanishing covariance, without a rate.
The second gives the rate $|U_n|^{-1/2}$ under the stronger summable covariance condition and is
used throughout.
Both follow from a direct variance bound and Chebyshev's inequality.

\begin{definition}[Vanishing covariance]\label{def:dying_cov}
A real-valued weakly stationary random field $\{V_s, s\in \mathbb{Z}^d\}$ has vanishing covariance if
$\cov(V_0, V_s)\to0$ as $\|s\|_\infty\to\infty$.
\end{definition}

The following lemma is an analog of the WLLN of \citet{Berstein1918law} \citep[see also][Exercise
254]{cacoullos1989exercises} for random fields on $\mathbb{Z}^d$. It does not require
stationarity, but the covariances must vanish uniformly in the locations.
See also \citet[Proposition 1]{hormann2013consistency} for a closely related result for spatially
distributed functional data.
\begin{lemma}[WLLN under vanishing covariance]\label{lemma:bernstein_wlln}
Let $\{V_s, s\in \mathbb{Z}^d\}$ be a real-valued random field with $\sup_{s\in\mathbb{Z}^d}
\operatorname{Var}(V_s)\leq c<\infty$ and
$$
\sup\left\{|\cov(V_s, V_t)| : s,t\in\mathbb{Z}^d,\ \|s-t\|_\infty\geq r\right\}\to0 \text{ as }
r\to\infty.
$$
Then for any sequence of nonempty finite sets $U_n\subset \mathbb{Z}^d$ with $|U_n|\to\infty$, we have
$$
\frac{1}{|U_n|}\sum_{s\in U_n} \left(V_s-E(V_s)\right) \xrightarrow{p} 0 \text{ as } n\to \infty.
$$
In particular, if $E(V_s)=\mu$ for all $s\in\mathbb{Z}^d$, then $|U_n|^{-1}\sum_{s\in U_n}V_s
\xrightarrow{p}\mu$ as $n\to\infty$.
\end{lemma}

\begin{corollary}[WLLN for weakly stationary fields]\label{cor:bernstein_wlln}
If $\{V_s, s\in \mathbb{Z}^d\}$ is a real-valued weakly stationary random field with vanishing
covariance, then for any sequence of nonempty finite sets $U_n\subset \mathbb{Z}^d$ with
$|U_n|\to\infty$, we have
$$
\frac{1}{|U_n|}\sum_{s\in U_n} V_s \xrightarrow{p} E(V_0) \text{ as } n\to \infty.
$$
\end{corollary}

\begin{definition}[Summable covariance]\label{def:summable_cov}
A real-valued weakly stationary random field $\{V_s, s\in\mathbb{Z}^d\}$ has summable covariance if
$\sum_{s\in\mathbb{Z}^d}|\cov(V_0,V_s)|<\infty$. This is also called the finite susceptibility
condition \citep[see e.g.][Chapter 10]{schmidt2014stochastic}. More generally, a weakly stationary
random field $\{\bm{V}_s\in\mathbb{R}^p, s\in\mathbb{Z}^d\}$ has summable covariance if
$\sum_{s\in\mathbb{Z}^d}|\cov(V_0^{(l)},V_s^{(m)})|<\infty$ for all $l,m\in\{1,\dots,p\}$, where
$V_s^{(l)}$ denotes the $l$th component of $\bm{V}_s$, that is, if all its auto- and
cross-covariances are absolutely summable.
\end{definition}

Summable covariance implies vanishing covariance, since the terms of an
absolutely convergent series tend to $0$.
The converse does not hold. For example, fields with long-range dependence can have covariances
that tend to $0$ too slowly to be summable.
Under summable covariance, the following lemma strengthens
Corollary~\ref{cor:bernstein_wlln} with a rate.

\begin{lemma}[Convergence of sample mean under summable covariance]\label{lemma:finite_sus_rate}
If $\{V_s, s\in \mathbb{Z}^d\}$ is a real-valued weakly stationary random field that has summable covariance, then for any {sequence of nonempty finite sets $U_n
\subset \mathbb{Z}^d$}, we have:
$$
\frac{1}{|U_n|}\sum_{s\in U_n} V_s = E(V_0) + O_p\left(\frac{1}{\sqrt{|U_n|}}\right) \text{ as } n\to \infty.
$$
In particular, {if $|U_n|\to\infty$,} this implies
the WLLN: $$\frac{1}{|U_n|}\sum_{s\in U_n}
V_s \xrightarrow{p} E(V_0) \text{ as } n\to \infty.$$
\end{lemma}
The proofs of Lemmas~\ref{lemma:bernstein_wlln} and \ref{lemma:finite_sus_rate} and
Corollary~\ref{cor:bernstein_wlln} are given in Appendix~\ref{append:proofs_pre}.

\subsection{Central limit theorem (CLT) for the sample mean}
\citet{bulinski2007limit} established the
CLT for $(BL, \theta)$-dependent random fields, which is a key result underpinning our theory.
{We first recall the distance between sets, regular growth, and the dependence notion,
then restate the theorem for vector-valued fields.}

\begin{definition}[Distance between two sets] For two finite sets $I, J\subset \mathbb{Z}^d$, the
distance between $I$ and $J$ is defined as $\operatorname{dist}(I, J) = \inf\{\|\bm{i}-\bm{j}\|_\infty: \bm{i}\in
I, \bm{j} \in J\}$.
\end{definition}

{The CLT also requires the domains to grow in a balanced way, so that the boundary of $U_n$ is
negligible compared with its volume.}
\begin{definition}[Regularly growing sets] \citep[see e.g.][Chapter 3, Definition 1.4]{bulinski2007limit} Let
$\delta U$ be the boundary of a set $U \subset \mathbb{Z}^d$, defined as $\delta U=\{\bm{j} \in
\mathbb{Z}^d, \bm{j} \notin U: \operatorname{dist}(\bm{j}, U) =1\}$.
A sequence of finite sets $U_n \subset \mathbb{Z}^d$, $n\in \mathbb{N}$, is called regularly growing
(to infinity) if
$$
|U_n|\to \infty \text{ and } |\delta U_n|/|U_n| \to 0 \text{ as } n \to  \infty.
$$
\end{definition}

\begin{example}
When $d=2$,
$U_n = \{1, 2, \dots, n\}\times \{1, 2, \dots, n\}$ is regularly growing, whereas $U_n = \{1\}
\times \{1, 2, \dots, n\}$ is not regularly growing because its boundary-to-volume ratio does not
vanish.
\end{example}

{The dependence notion bounds the covariance between Lipschitz functions of the field on two
separated sets by a coefficient $\theta_r$ that decays with their distance $r$.}
\begin{definition}[$(BL,\theta)$-dependent random field] \citep[see e.g.][Chapter 1, Definitions
5.12 and 5.14]{bulinski2007limit} Let $BL(n)$ be the class of bounded Lipschitz functions
$f:\mathbb{R}^n\rightarrow \mathbb{R}, n \in \mathbb{N}$.
A vector-valued random field $\{\bm{V}_s \in \mathbb{R}^p, s \in \mathbb{Z}^d\}$ is called $(BL,
\theta)$-dependent if there exists a nonincreasing sequence $\theta=\left(\theta_r\right)_{r \in
\mathbb{N}}, \theta_r \rightarrow 0$ as $r \rightarrow \infty$, such that for all finite disjoint
sets $I, J \subset \mathbb{Z}^d$ with $\operatorname{dist}(I, J)=r \in \mathbb{N}$, and any functions $f \in
BL(p|I|), g \in B L(p|J|)$, the following holds:
\begin{equation}|\cov(f(\bm{V}_I), g(\bm{V}_J))| \leq \operatorname{Lip}(f)
\operatorname{Lip}(g)(|I| \wedge|J|) \theta_r.
\label{eq:bl}
\end{equation}

Here, for a finite set $I = \{s_1, \dots, s_k\} \subset \mathbb{Z}^d$ listed in a fixed order,
$\bm{V}_I = [\bm{V}_{s_1}', \dots, \bm{V}_{s_k}']' \in \mathbb{R}^{p|I|}$.
\end{definition}

\begin{remark}\label{BL-general}
{\citet[Chapter 1, Definitions 5.12 and 5.14]{bulinski2007limit} appear to require $\theta_r<\infty$ for
all $r$, since their proof of the CLT uses a bound involving $\theta_1$.
We allow $\theta_r=\infty$ for finitely many $r$, so that $m$-dependent random fields
(Definition~\ref{def:mdep}) are $(BL,\theta)$-dependent, with $\theta_r=\infty$ for $r\leq m$ and
$\theta_r=0$ for $r>m$.
Lemma~\ref{lemma:clt_general} in Appendix~\ref{append:clt} shows that the CLT below still holds
under this generalized definition.}
\end{remark}

\begin{remark}\label{remark:finite_sus}
{A real-valued weakly stationary $(BL,\theta)$ random field has summable covariance \citep[Chapter 3, Remark 1.10]{bulinski2007limit}.
Lemma~\ref{lemma:abs_sum} in Appendix~\ref{append:clt} shows that this remains true under the
generalized definition
(Remark~\ref{BL-general}).}
\end{remark}

We show in the following lemma that the $(BL,\theta)$ property is preserved under Lipschitz
transformations.
\begin{lemma}
\label{lemma:bl}
Suppose a vector-valued random field $\{\bm{V}_s \in \mathbb{R}^p:s\in \mathbb{Z}^d\}$ is
$(BL,\theta)$, and $h: \mathbb{R}^p \rightarrow \mathbb{R}^q$ is a function such that $h = [h_1,
\dots, h_q]$ and $h_i$ is $l_i$-Lipschitz for $i = 1, \dots, q$.
Then $\{h(\bm{V}_s)\in \mathbb{R}^q:s\in \mathbb{Z}^d\}$ is $(BL,\theta)$.
\end{lemma}

We now state the CLT of \citet{bulinski2007limit} that underlies our results. We
state it for vector-valued fields, since Section~\ref{sec:theorem} applies it to the field of
squares and cross-products of the standardized fields.
\begin{theorem*}[\cite{bulinski2007limit}, Chapter 3, Theorem 1.12 {and Corollary 1.13}]
{Let $V=\{\bm{V}_s\in\mathbb{R}^p, s\in \mathbb{Z}^d\}$ be a $(BL,\theta)$-dependent
strictly stationary centered square-integrable random field.
Let $\bm{S}(U_n) = \sum_{s\in U_n}\bm{V}_s$. Then, for any sequence
of regularly growing sets $U_n\subset \mathbb{Z}^d$, one has
\[
\frac{\bm{S}(U_n)}{\sqrt{|U_n|}}
\xrightarrow{d}
\mathcal{N}\Big(\bm{0},\sum_{s \in \mathbb{Z}^d}\cov(\bm{V}_0, \bm{V}_s)\Big)
\quad \text{as } n \to \infty ,
\]
where $\cov(\bm{V}_0, \bm{V}_s)$ is the $p\times p$ matrix with $(l,m)$ entry $\cov(V_0^{(l)},
V_s^{(m)})$, and $V_s^{(l)}$ denotes the $l$th component of $\bm{V}_s$.
Every entry of the matrix series converges absolutely (Lemma~\ref{lemma:abs_sum} in
Appendix~\ref{append:clt}), so the limit covariance does not depend on the order of summation.}
\end{theorem*}

\subsection{Classes of random fields and dependence coefficients}
The following classes of random fields appear in the special cases of Section~\ref{sec:theorem}
and in the sufficient conditions of Section~\ref{sec:condition}.
{The $m$-dependent fields are the simplest fields with finite-range dependence.
Gaussian fields give the explicit asymptotic variances of Corollaries~\ref{cor:zd_gaussian} and
\ref{cor:zd_rho0_gaussian}.
Positively and negatively associated fields give a route to $(BL,\theta)$-dependence through
covariances alone.}

We generalize the definition of $m$-dependent real-valued random fields \citep[see e.g.][Page
91]{bulinski2007limit} to vector-valued random fields.
\begin{definition}[Vector-valued $m$-dependent random field]\label{def:mdep}
A vector-valued random field $\{\bm{V}_s \in \mathbb{R}^p,\allowbreak s \in \mathbb{Z}^d\}$ is
called $m$-dependent if for all finite sets $I,J\subset \mathbb{Z}^d$ with $\operatorname{dist}(I, J)>m$,
$\{\bm{V}_s, s\in I\}$ and $\{\bm{V}_s, s\in J\}$ are independent.
\end{definition}
{ Note that \citet[Page 91]{bulinski2007limit} instead require independence when $\operatorname{dist}(I, J)\geq
m$. We use $\operatorname{dist}(I, J)>m$ for consistency with the time series literature
\citep[see e.g.][]{brockwell1991time}.}

\begin{definition}[Vector-valued Gaussian random field]\label{def:gaussian}

A vector-valued random field $\{\bm{V}_s \in \mathbb{R}^p, s \in
\mathbb{Z}^d\}$ is Gaussian if for any finite set $I \subset \mathbb{Z}^d$, the concatenated
vector $\bm{V}_I \in \mathbb{R}^{p|I|}$
follows a (possibly degenerate) multivariate normal distribution.
{This is stronger than the definition of \citet[Page 11]{adler2007random}, which only
requires $\{\bm{\alpha}'\bm{V}_s, s\in\mathbb{Z}^d\}$ to be a real-valued Gaussian random field for every $\bm{\alpha}\in\mathbb{R}^p$.}
\end{definition}

{Positive and negative association require bounded, coordinate-wise nondecreasing functions
of the field on disjoint sets to be nonnegatively or nonpositively correlated, respectively.}
\begin{definition}[PA random field]\label{def:pa} \citep[see e.g.][Chapter 1, Definition 1.2 or 5.13]{bulinski2007limit} Let
$M(n)$ be the class of bounded coordinate-wise nondecreasing Borel functions
$f:\mathbb{R}^n\rightarrow \mathbb{R}, n \in \mathbb{N}$.
A vector-valued random field $\{\bm{V}_s \in \mathbb{R}^p, s \in \mathbb{Z}^d\}$ is called weakly
associated or positively associated (PA) if for all finite disjoint sets $I, J \subset \mathbb{Z}^d$
and all functions $f \in M(p|I|), g \in M(p|J|)$, the following holds:
$$\cov(f(\bm{V}_I), g(\bm{V}_J )) \geq 0.$$
\end{definition}

\begin{definition}[NA random field]\label{def:na} \citep[see e.g.][Chapter 1, Definition 1.3 or 5.13]{bulinski2007limit}
A vector-valued random field $\left\{\bm{V}_s \in \mathbb{R}^p, s \in \mathbb{Z}^d\right\}$ is
called negatively associated (NA) if for all finite disjoint sets $I, J \subset \mathbb{Z}^d$ and
all functions $f \in M(p|I|), g \in M(p|J|)$, the following holds:
$$\cov(f(\bm{V}_I), g(\bm{V}_J )) \leq 0.$$
\end{definition}

Nondecreasing transformations preserve the PA (NA) property, as the following lemma shows.
\begin{lemma}
\label{lemma:pa}
If a vector-valued random field $\left\{\bm{V}_s \in \mathbb{R}^p:s\in \mathbb{Z}^d\right\}$ is PA
(NA), and $h: \mathbb{R}^p \rightarrow \mathbb{R}^q$ is a function such that $h = [h_1, \dots, h_q]$
and $h_i$ is coordinate-wise nondecreasing Borel for $i = 1, \dots, q$, then $\left\{h(\bm{V}_s)\in
\mathbb{R}^q:s\in \mathbb{Z}^d\right\}$ is PA (NA).
This includes the special case $q=1$ and $h(\bm{\eta})=\bm{a}'\bm{\eta}$ for a fixed
$\bm{a}\in\mathbb{R}^p$ with nonnegative components, noted by \citet{burton1986invariance}.
\end{lemma}

{For PA, NA and Gaussian fields, $(BL,\theta)$-dependence can be checked through covariances
alone, as the following remark shows.}

\begin{remark}\label{remark:cox_bl}
By \citet[Chapter 1, Theorem 5.17]{bulinski2007limit}, a weakly stationary random
field $V=\{\bm{V}_s,s\in\mathbb{Z}^d\}$ is $(BL,\theta)$ if it has summable covariance and either $V$
is PA or NA with $E\|\bm{V}_s\|_2^2<\infty$ for all $s$, or $V$ is Gaussian.
\citet{bulinski2007limit} state the covariance condition through the Cox--Grimmett coefficient
\[
u_r=\sum_{\|j\|_\infty\ge r}\ \sum_{l,m=1}^p|\cov(V_0^{(l)},V_j^{(m)})|,
\]
which is finite for every
$r$ and tends to $0$ if and only if $V$ has summable covariance.
Lemma~\ref{lemma:cox} gives a simple sufficient condition, exponential decay of the covariances.
\end{remark}

\begin{lemma}\label{lemma:cox}

Let $\{\bm{V}_s = [V^{(1)}_s, \dots, V^{(p)}_s]'\in \mathbb{R}^p, s\in\mathbb{Z}^d\}$ be a weakly
stationary random field whose auto- and cross-covariances decay exponentially, that is, there exist
constants $C>0$ and $\lambda>0$ such that
\begin{equation}\label{eq:exp_decay}
|\cov(V^{(l)}_t, V^{(m)}_s)|\le C\exp(-\lambda\|t-s\|_\infty) \quad \text{for all } l,m
=1,\dots,p \text{ and } t,s\in \mathbb{Z}^d.
\end{equation}
Then $\bm{V}$ has summable covariance.

Condition \eqref{eq:exp_decay} holds, for example, for squared-exponential covariances: if
{for some positive definite matrix $A$,
$|\cov(V^{(l)}_t, V^{(m)}_s)|\le C\exp(-(t-s)'A(t-s))$}, then
\eqref{eq:exp_decay} holds with $\lambda$ equal to the smallest eigenvalue of $A$.
\end{lemma}

The proofs of Lemmas~\ref{lemma:bl}, \ref{lemma:pa} and \ref{lemma:cox} are given in
Appendix~\ref{append:proofs_pre}.

\subsection{Assumptions}
We consider the following assumptions.
The first three are standard, and the last two control the dependence, with a weaker version for
the case $\rho=0$.
\begin{assumption}[Strict stationarity]\label{ass:zd_stationary}
The random field $Z$ is strictly stationary.
\end{assumption}
\begin{assumption}[Finite fourth moments {and nondegeneracy}]\label{ass:zd_squareint}
$E\|\bm{Z}_0\|_2^4<\infty${, $\operatorname{Var}(X_0)>0$ and
$\operatorname{Var}(Y_0)>0$}.
\end{assumption}

\begin{assumption}[Regularly growing domains]\label{ass:regular}
The sequence of domains $U_n \subset \mathbb{Z}^d$ is regularly growing.
\end{assumption}

\begin{assumption}[Weak dependence, general $\rho$]\label{ass:bl}
Let $\tilde{\bm{T}}_s = [\tilde X_s^2, \tilde Y_s^2, \tilde X_s \tilde Y_s]' \in \mathbb{R}^3$ for
$s\in\mathbb{Z}^d$ and $\tilde T = \{\tilde{\bm{T}}_s, s\in\mathbb{Z}^d\}$.
$\tilde T$ is $(BL,\theta)$, and $\tilde X$ and $\tilde Y$ each have summable covariance.
\end{assumption}

\begin{assumptionStar}[Weak dependence, $\rho=0$ case]\label{ass:finite}
Let $\tilde X^2 = \{\tilde{X}_s^2:s\in \mathbb{Z}^d\}$,
$\tilde Y^2 = \{\tilde{Y}_s^2:s\in \mathbb{Z}^d\}$ and $\tilde X\tilde Y =
\{\tilde{X}_s\tilde{Y}_s:s\in \mathbb{Z}^d\}$.
$\tilde X\tilde Y$ is $(BL,\theta)$, $\tilde X$ and $\tilde Y$ each have summable
covariance, and $\tilde X^2$ and $\tilde Y^2$ each have vanishing covariance.
\end{assumptionStar}

\section{Asymptotic normality of the sample correlation}\label{sec:theorem}
This section establishes that $\hat \rho$ is a CAN estimator of $\rho$ and connects the
result to existing methods.
Section~\ref{sec:theorem_general} treats general $\rho$.
Section~\ref{sec:theorem_rho0} treats $\rho=0$ under weaker
conditions, which is of interest in hypothesis testing.
Section~\ref{sec:theorem_indep} specializes to independent fields and shows that the variance
approximations behind the modified $t$-test of \citet{clifford1989assessing} and
\citet{dutilleul1993modifying} converge to our asymptotic variance.
The proofs of the results in this section are given in Appendix~\ref{append:proofs_theorem}.

\subsection{General correlation}\label{sec:theorem_general}
Theorem~\ref{theorem:zd} requires the field $\tilde{T}$ of squares and cross-products to be
$(BL,\theta)$ (Assumption~\ref{ass:bl}).
Section~\ref{sec:condition} shows when this holds in terms of $X$ and $Y$.
\begin{theorem}
\label{theorem:zd}
Under Assumptions \ref{ass:zd_stationary}, \ref{ass:zd_squareint}, \ref{ass:regular} and
\ref{ass:bl},
we have
$\sqrt{|U_n|}(\hat \rho-\rho) \xrightarrow{d} \mathcal{N}(0, \sigma^2)$ as $n \to \infty$.
The asymptotic variance is $\sigma^2=  \sum_{k\in\mathbb{Z}^d}h_k$, where

\begin{equation}\label{eq:variance_v}
h_k = \cov(L_0, L_k),
\qquad
L_s=\tilde X_s\tilde Y_s-\frac{\rho}{2}\big(\tilde X_s^2+\tilde Y_s^2\big),
\end{equation}
so that $\sigma^2$ is the long-run variance of the scalar field $L=\{L_s, s\in\mathbb{Z}^d\}$.
Since $E[L_s]=0$, $h_k=E[L_0L_k]$, which expands to
\begin{align}\label{eq:variance}
h_k &=
\frac{1}{4} \rho^2 E\left[\tilde{X}_0^2 \tilde{X}_{k}^2+\tilde{Y}_0^2 \tilde{Y}_{k}^2+\tilde{X}_0^2
\tilde{Y}_{k}^2+\tilde{X}_{k}^2 \tilde{Y}_0^2\right] \nonumber \\
&\quad - \frac{1}{2} \rho E\left[\left(\tilde{X}_0 \tilde{X}_{k}+\tilde{Y}_0
\tilde{Y}_{k}\right)\left(\tilde{X}_0 \tilde{Y}_{k}+\tilde{X}_{k} \tilde{Y}_0\right)\right]
\nonumber \\
&\quad + E\left[\tilde{X}_0 \tilde{X}_{k} \tilde{Y}_0 \tilde{Y}_{k}\right].
\end{align}
\end{theorem}

{
\begin{remark}\label{remark:sigma_zero}
Theorem~\ref{theorem:zd} allows $\sigma^2=0$, in which case the limit is the point mass at $0$.
This happens, for example, when $|\rho|=1$, because then $\tilde X_s=\pm\tilde Y_s$ almost surely
and $\hat\rho=\rho$ for every $n$.
\end{remark}
}

{For Gaussian fields, the fourth moments in \eqref{eq:variance} reduce to products of auto- and
cross-correlations, which gives a simpler formula.}
\begin{corollary}\label{cor:zd_gaussian}
{Under the assumptions of} Theorem \ref{theorem:zd}, if $Z$ is a Gaussian random field, then
$\sqrt{|U_n|}(\hat \rho-\rho) \xrightarrow{d} \mathcal{N}(0, \sigma^2)$ as $n \to \infty$, where
{$\sigma^2=\sum_{k\in\mathbb{Z}^d}h_k$ with}
\begin{align}\label{eq:variance_gaussian}
    {h_k}&=\frac{1}{2}\rho^2 \left(\rho_{x x, k}^2+\rho_{yy, k}^2+\rho_{xy,
k}^2+\rho_{xy,-k}^2\right) \nonumber\\
&\quad - \rho \left(\rho_{xx, k}+\rho_{yy, k}\right)\left(\rho_{xy,
k}+\rho_{xy,-k}\right) \nonumber\\
&\quad + \rho_{x x, k} \rho_{yy, k}+\rho_{x y, k} \rho_{x
y,-k}.
\end{align}
\end{corollary}

\begin{remark}
\label{remark:zd_iid}
One special case of Theorem~\ref{theorem:zd} is when the pairs $(X_s, Y_s)$, $s \in
\mathbb{Z}^d$, are i.i.d. Then the asymptotic variance in \eqref{eq:variance} reduces to \[\sigma^2 = h_0=\operatorname{Var}(L_0) = \frac{1}{4}\rho^2
E\left[\tilde{X}_0^4+\tilde{Y}_0^4\right]+\left(\frac{1}{2} \rho^2+1\right) E\left[\tilde{X}_0^2
\tilde{Y}_0^2\right]-\rho E\left[\tilde{X}_0^3 \tilde{Y}_0+\tilde{X}_0\tilde{Y}_0^3\right].\]
If the pairs are also Gaussian, Corollary~\ref{cor:zd_gaussian} with only the $k=0$ term gives
$\sigma^2=(1-\rho^2)^2$.
{In both cases $\sigma^2$ does not depend on $d$ or on the shape of $U_n$, since the observations
in $U_n$ are then $N=|U_n|$ i.i.d.\ pairs.
So $\sigma^2$ is the classical asymptotic variance of $\sqrt{N}(\hat\rho-\rho)$ for an i.i.d.\
sample of size $N$ (see, e.g., \citealp[Example~3.6]{van2000asymptotic}, and
\citealp[Theorem~4.2.4]{anderson2003introduction} for the Gaussian case).}
\end{remark}

\begin{remark}
\label{remark:1d}
For $d=1$, the asymptotic variance \eqref{eq:variance_gaussian} is the same as that of
\cite{roy1989asymptotic}.
Of note, \cite{roy1989asymptotic} used the assumption \quotes{all cumulants of fourth order are 0},
which is satisfied by a Gaussian distribution but not other well-known distributions.
{Expression \eqref{eq:variance_gaussian} is also asymptotically equivalent to the
finite-sample variance of \citet{afyouni2019effective} (Appendix~\ref{append:afyouni}), which they
derived from a first-order Taylor expansion of $\hat\rho$ rather than from asymptotic arguments.}
\end{remark}

To see what Corollary~\ref{cor:zd_gaussian} implies in practice, we compute $\sigma^2$ for a
simple pair of correlated Gaussian random fields, in which $X$ and $Y$ have the same autocorrelation and
their cross-correlation is proportional to it.
\begin{example}[Smoothed Gaussian fields]\label{example:smoothed}
Let $W_1$ and $W_2$ be independent, centered, unit-variance, stationary Gaussian random fields with
the same absolutely summable autocorrelation $r_k$, and let $X=W_1$ and
$Y=\rho W_1+\sqrt{1-\rho^2}\,W_2$ with $|\rho|<1$. Then $\rho_{xx,k}=\rho_{yy,k}=r_k$ and
$\rho_{xy,k}=\rho_{xy,-k}=\rho\, r_k$. The field $Z=\{[X_s,Y_s]',s\in\mathbb{Z}^d\}$ is stationary and Gaussian with unit
variances, so Assumptions~\ref{ass:zd_stationary} and \ref{ass:zd_squareint} hold, and
Assumption~\ref{ass:bl} is verified in Section~\ref{sec:gaussian}. Hence, for regularly growing
domains $U_n$ (Assumption~\ref{ass:regular}),
Corollary~\ref{cor:zd_gaussian} applies, and substituting into \eqref{eq:variance_gaussian} gives $\sigma^2=(1-\rho^2)^2\sum_{k\in\mathbb{Z}^d}r_k^2$.
For instance, suppose $W_1$ and $W_2$ are obtained by smoothing i.i.d.\ $\mathcal{N}(0,1)$ noise on
$\mathbb{Z}^d$ with the Gaussian kernel $g(u)=\exp\{-\|u\|_2^2/(2\ell^2)\}$, $u\in\mathbb{Z}^d$, and
rescaling to unit variance. Replacing lattice sums by integrals gives
$r_k\approx\exp\{-\|k\|_2^2/(4\ell^2)\}=\prod_{i=1}^d\exp\{-k_i^2/(4\ell^2)\}$, so
$\sum_{k\in\mathbb{Z}^d}r_k^2\approx(\sqrt{2\pi}\,\ell)^d$ and
$\sigma^2\approx(1-\rho^2)^2(\sqrt{2\pi}\,\ell)^d$. With $N=|U_n|$ observed grid points, the variance of
$\hat\rho$ is therefore approximately $\sigma^2/N\approx(1-\rho^2)^2(\sqrt{2\pi}\,\ell)^d/N$.
\end{example}

The approximation in Example~\ref{example:smoothed} is the i.i.d.\ variance $(1-\rho^2)^2/N$
inflated by the factor $(\sqrt{2\pi}\,\ell)^d$. Thus $N/(\sqrt{2\pi}\,\ell)^d$ acts as an effective number of independent
observations. The precision of $\hat\rho$ depends on $\rho$ through $(1-\rho^2)^2$, as in the
i.i.d.\ case, and on the smoothness $\ell$, the number of grid points $N$ and the dimension $d$ through
$(\sqrt{2\pi}\,\ell)^d/N$. The formula also gives the number of grid points needed
for a target precision. A standard error of at most $\varepsilon$ requires approximately
$N\ge(1-\rho^2)^2(\sqrt{2\pi}\,\ell)^d/\varepsilon^2$. For example, with $\rho=0.5$, $\ell=2$
and $\varepsilon=0.05$, this gives about $1{,}130$, $5{,}700$ and $28{,}000$ grid points for $d=1$, $2$ and $3$. Since
$\sqrt{2\pi}\,\ell>1$ for $\ell\ge1$, the number of grid points needed grows exponentially in $d$, a
form of the curse of dimensionality.

\subsection{Zero correlation under weaker conditions}\label{sec:theorem_rho0}
When $\rho=0$, Theorem~\ref{theorem:zd} can be proved under weaker conditions.
Assumption~\ref{ass:finite} only requires $\tilde X\tilde Y$ to be $(BL,\theta)$, the
fields $\tilde X$ and $\tilde Y$ each to have summable covariance, and the squared fields $\tilde X^2$
and $\tilde Y^2$ each to have vanishing covariance, whereas
Assumption~\ref{ass:bl} requires the three-dimensional field $\tilde{\bm{T}}$ to be $(BL,\theta)$.
This matters in practice: for independent unbounded non-Gaussian fields,
Assumption~\ref{ass:bl} is hard to verify, whereas Assumption~\ref{ass:finite} can be verified with
Lemmas~\ref{lemma:product} and \ref{lemma:square_dying} (see
Proposition~\ref{prop:indep}).
\begin{theorem}
\label{theorem:zd_rho0}

Under Assumptions \ref{ass:zd_stationary}, \ref{ass:zd_squareint}, \ref{ass:regular} and
\ref{ass:finite}, if $\rho=0$, then
we have $\sqrt{|U_n|}\hat \rho \xrightarrow{d} \mathcal{N}(0, \sigma^2)$ as $n\to\infty$.
The asymptotic variance is $\sigma^2=\sum_{k\in\mathbb{Z}^d}h_k$, where

\begin{equation}\label{eq:variance_v0}
h_k = \cov(L_0, L_k),
\qquad
L_s=\tilde X_s\tilde Y_s,
\end{equation}
which is \eqref{eq:variance_v} with $\rho=0$.
Since $E[L_s]=\rho=0$, $h_k=E[L_0L_k]$, that is,
\begin{equation}\label{eq:variance_0}
h_k = E\left[\tilde{X}_0 \tilde{X}_{k} \tilde{Y}_0 \tilde{Y}_{k}\right].
\end{equation}
\end{theorem}

\begin{remark}
The same results can be obtained by taking the special case of $\rho=0$ in Theorem \ref{theorem:zd}.
However, Theorem \ref{theorem:zd_rho0} has a simpler proof without using the delta method and uses
weaker assumptions than Theorem \ref{theorem:zd}.
To see that Assumption \ref{ass:finite} is implied by Assumption \ref{ass:bl}, note that if $\tilde
T$ is $(BL,\theta)$, then so is each of its coordinates by Lemma \ref{lemma:bl} applied with the
coordinate projections, and hence $\tilde X^2$ and $\tilde Y^2$ each have summable covariance by Remark \ref{remark:finite_sus}, which implies vanishing covariance.
\end{remark}

{As for general $\rho$, the asymptotic variance simplifies for Gaussian fields.}
\begin{corollary}\label{cor:zd_rho0_gaussian}
{Under the assumptions of} Theorem \ref{theorem:zd_rho0}, if $Z$ is a Gaussian random field, then
$\sqrt{|U_n|}\hat \rho \xrightarrow{d} \mathcal{N}(0, \sigma^2)$ as $n\to\infty$, where
{$\sigma^2=\sum_{k\in\mathbb{Z}^d}h_k$ with}
\begin{align}\label{eq:variance_gaussian_0}
    {h_k}&=\rho_{x x, k} \rho_{yy, k}+\rho_{x y, k} \rho_{x y,-k},
\end{align}
and $\rho_{xy,0}=\rho=0$, so that $h_0=1$.
\end{corollary}

\subsection{Independent fields and effective sample size methods}\label{sec:theorem_indep}
Testing $H_0^{\text{indep}}: X\indep Y$ is the setting of the effective sample size methods
reviewed in the Introduction.
Under independence, the asymptotic variance takes a particularly simple form.
\begin{corollary}\label{cor:zd_indep}
Under the assumptions of Theorem \ref{theorem:zd_rho0}, if $X$ and $Y$ are independent, then
$\sqrt{|U_n|}\hat \rho \xrightarrow{d} \mathcal{N}(0, \sigma^2)$ as $n\to\infty$, where
{$\sigma^2=\sum_{k\in\mathbb{Z}^d}h_k$ with}
\begin{align}\label{eq:variance_indep}
    {h_k}&=\rho_{x x, k} \rho_{yy, k}.
\end{align}
\end{corollary}

\begin{remark}\label{remark:ESS_1d}
A special case of Corollary~\ref{cor:zd_indep} is when $X$ and $Y$ are independent time
series ($d=1$), for which {the asymptotic variance} becomes $\sigma^2 =
1+2\sum_{k=1}^{\infty}\rho_{xx,k}\rho_{yy,k}$.
Thus, with $n$ observed time points (that is, $U_n=\{1,\dots,n\}$, so $|U_n|=n$), the asymptotic variance of $\sqrt n\,\hat\rho$ is $\sigma^2$, which agrees with the
approximation $\operatorname{Var}(\hat\rho)\approx\sigma^2/n$ of \citet{quenouille1947notes}.
\end{remark}

The equivalence in Remark~\ref{remark:ESS_1d} extends to random fields of any fixed dimension $d\geq1$.
For independent Gaussian fields $X$ and $Y$, \citet{clifford1989assessing} approximated $\operatorname{Var}(\hat\rho)$
by $V_n^{\mathrm{C}} =
\operatorname{tr}(\Sigma_X\Sigma_Y)/\{\operatorname{tr}(\Sigma_X)\operatorname{tr}(\Sigma_Y)\}$
\citep[eq.~3.4]{dutilleul1993modifying}, where $\Sigma_X$ and $\Sigma_Y$ are the $|U_n|\times|U_n|$
covariance matrices of $\{X_s: s\in U_n\}$ and $\{Y_s: s\in U_n\}$ under a common ordering of $U_n$.
\citet{dutilleul1993modifying} refined it to $V_n^{\mathrm{D}} = \operatorname{tr}(\Pi\Sigma_X
\Pi\Sigma_Y)/\{\operatorname{tr}(\Pi\Sigma_X) \operatorname{tr}(\Pi\Sigma_Y)\}$, where $\Pi = I_{|U_n|} -
|U_n|^{-1}\bm{1}_{|U_n|}\bm{1}_{|U_n|}'$ is the centering matrix, and found numerically that the two
differ for small $|U_n|$ but agree as $|U_n|$ grows.

\begin{theorem}[Asymptotic equivalence with Clifford and Dutilleul]\label{theorem:clifford}

Let $X$ and $Y$ be weakly stationary with positive variances, let each have summable covariance, and let $U_n$ be regularly growing.
Then $|U_n|\, V_n^{\mathrm{C}} \to \sum_{k\in\mathbb{Z}^d}\rho_{xx,k}\rho_{yy,k}$ and
$|U_n|\, V_n^{\mathrm{D}} \to \sum_{k\in\mathbb{Z}^d}\rho_{xx,k}\rho_{yy,k}$ as $n\to\infty$.
Under the assumptions of Corollary~\ref{cor:zd_indep}, this limit is the asymptotic variance
$\sigma^2$ in \eqref{eq:variance_indep}.
\end{theorem}

\begin{remark}\label{remark:ESS_zd}

By Theorem~\ref{theorem:clifford}, both approximations are asymptotically
equivalent to our result and to each other, which establishes the large-sample agreement observed
by \citet{dutilleul1993modifying}.
Theorem~\ref{theorem:clifford} does not use Gaussianity, so the population targets of
both approximations equal the asymptotic variance of $\hat\rho$ for non-Gaussian independent fields
as well. Proposition~\ref{prop:indep} gives conditions under which such fields
satisfy the assumptions of Corollary~\ref{cor:zd_indep}.
\end{remark}

\section{Conditions for sample correlation asymptotic normality}\label{sec:condition}

Theorems~\ref{theorem:zd} and \ref{theorem:zd_rho0} rely on Assumption~\ref{ass:bl} or
\ref{ass:finite}, which concern the transformed fields $\tilde X^2$, $\tilde Y^2$ and
$\tilde X\tilde Y$ and are hard to check directly.
This section gives sufficient conditions on $X$ and $Y$ for bounded, Gaussian,
$m$-dependent and independent fields.
Together with Assumptions~\ref{ass:zd_stationary}, \ref{ass:zd_squareint} and \ref{ass:regular},
each of these conditions guarantees that $\hat \rho$ is CAN.
{When $d=1$, limit theorems for time series, such as the CLT for the sample mean,
have historically been established under mixing conditions.
Such conditions involve $\sigma$-fields and are hard to check in practice
\citep{doukhan1999new}, and for $d>1$ few examples of mixing random fields are known beyond
$m$-dependent ones \citep[Chapter 1, Section 5]{bulinski2007limit}.
For completeness, Appendix~\ref{append:mixing} gives versions of our results for
$\hat\rho$ under strong mixing when $d=1$.}
The proofs of the results in this section are given in Appendix~\ref{append:proofs_condition}.

\subsection{\texorpdfstring{Bounded fields}{Bounded fields}}\label{sec:bounded}

We begin with bounded fields, the most tractable class. Squares and products are not
Lipschitz on $\mathbb{R}^2$, so Lemma~\ref{lemma:bl} does not apply to them directly. On a bounded
set, however, they are Lipschitz, so the $(BL,\theta)$ structure of $Z$ is inherited by the squares
and cross-products that enter Assumption \ref{ass:bl}, as the following lemma shows.
\begin{lemma}[Squares and products of bounded fields]\label{lemma:bounded_bl}
If $\{X_s \in \mathbb{R}:s\in \mathbb{Z}^d\}$ and $\{Y_s \in \mathbb{R}:s\in
\mathbb{Z}^d\}$ are both bounded and $\{[X_s, Y_s]'\in \mathbb{R}^2:s\in
\mathbb{Z}^d\}$ is $(BL,\theta)$, then $\{[X_s^2, Y_s^2, X_sY_s ]'\in
\mathbb{R}^3:s\in \mathbb{Z}^d\}$ is $(BL,\theta)$.
\end{lemma}
Boundedness also ensures that $\tilde{\bm{T}}$ is square integrable.
The three conditions below are nested, from the most general to the easiest to verify.
\begin{proposition}\label{prop:bounded}
Under Assumptions \ref{ass:zd_stationary} and \ref{ass:zd_squareint}, if $X$ and $Y$ are bounded,
then Assumption \ref{ass:bl} holds when any of the following holds:
\begin{enumerate}[label=(B\arabic*), ref=B\arabic*]
    \item \label{ass:B1} $Z$ is $(BL,\theta)$.
    \item (stronger than B1) \label{ass:B2} $Z$ is PA or NA and has summable
covariance.
    \item (stronger than B2)  \label{ass:B3} $Z$ is PA or NA, and its auto- and cross-covariances decay
exponentially as in \eqref{eq:exp_decay}.
\end{enumerate}
\end{proposition}

\subsection{Gaussian fields}\label{sec:gaussian}
We now consider Gaussian fields, for which Assumption~\ref{ass:bl} reduces to a condition on the
covariances of $Z$ alone.
The key step is the following Gaussian analogue of Lemma~\ref{lemma:bounded_bl}. Although squares
and products are not Lipschitz, they preserve the $(BL,\theta)$-dependence of Gaussian fields.
\begin{lemma}[Squares and products of Gaussian fields]\label{lemma:gauss_square}
If $Z=\{[X_s,Y_s]'\in\mathbb{R}^2:s\in\mathbb{Z}^d\}$ is a weakly stationary Gaussian random field
that is $(BL,\theta)$, then $\{[X_s^2,Y_s^2,X_sY_s]'\in\mathbb{R}^3:s\in\mathbb{Z}^d\}$ is
$(BL,\theta)$.
\end{lemma}
Applying Lemma~\ref{lemma:gauss_square} to the standardized fields gives the following result.
\begin{proposition}\label{prop:gaussian}
Under Assumptions~\ref{ass:zd_stationary} and \ref{ass:zd_squareint}, if $Z$ is a Gaussian random
field, then Assumption~\ref{ass:bl} holds when any of the following holds:
\begin{enumerate}[label=(N\arabic*), ref=N\arabic*]
    \item \label{ass:N1} $Z$ is $(BL,\theta)$.
    \item (equivalent to N1) \label{ass:N2} $Z$ has summable covariance.
    \item (stronger than N2) \label{ass:N3} The auto- and cross-covariances of $Z$ decay
exponentially as in \eqref{eq:exp_decay}.
\end{enumerate}
\end{proposition}

\begin{remark}\label{remark:gaussian}
Proposition~\ref{prop:gaussian} requires neither independence of $X$ and $Y$ nor positive or
negative association. The proof of Lemma~\ref{lemma:gauss_square} bounds the covariance of smooth
functions of a Gaussian vector by
the covariances of the vector itself, through a Gaussian covariance identity
\citep[Corollary~4.12]{bandeira2023matrix}, and passes to Lipschitz functions by convolution
approximation \citep{shashkin2002quasi}. Gaussian fields have finite fourth moments, and
Assumption~\ref{ass:bl} implies Assumption~\ref{ass:finite}. Hence, by
Proposition~\ref{prop:gaussian}, Corollaries~\ref{cor:zd_gaussian} and \ref{cor:zd_rho0_gaussian}
hold whenever $Z$ is a strictly stationary Gaussian field with summable covariance and positive
variances, observed on regularly growing domains, with $\rho=0$ for
Corollary~\ref{cor:zd_rho0_gaussian}. For the fields of Example~\ref{example:smoothed}, the auto- and
cross-covariances $r_k$ and $\rho\,r_k$ are absolutely summable, so Condition~\ref{ass:N2} holds
and Assumption~\ref{ass:bl} follows.
\end{remark}

\subsection{\texorpdfstring{$m$-dependent fields}{m-dependent fields}}\label{sec:mdep}
We next consider $m$-dependent fields, where the dependence range is finite.
Each $\tilde{\bm{T}}_s$ is a function of $\bm{Z}_s$ alone, so $\tilde T$ is
$m$-dependent whenever $Z$ is, and $m$-dependent fields are $(BL,\theta)$ under the generalized
definition
(Remark~\ref{BL-general}). Hence no condition beyond Assumptions~\ref{ass:zd_stationary} and
\ref{ass:zd_squareint} is needed to verify Assumption \ref{ass:bl}.

\begin{proposition}\label{prop:mdep}
Under Assumption{s} \ref{ass:zd_stationary} {and \ref{ass:zd_squareint}}, if $Z$ is $m$-dependent, then Assumption \ref{ass:bl} holds.
\end{proposition}

For example, a Gaussian field is $m$-dependent if and only if all its auto- and cross-covariances
vanish at lags $k$ with $\|k\|_\infty>m$ \citep[Theorem~2.4.2]{anderson2003introduction}.

\subsection{\texorpdfstring{Independent fields}{Independent fields}}\label{sec:indep}
Finally, we consider independent fields, the setting of Theorem~\ref{theorem:clifford} and of
methods such as \cite{clifford1989assessing}. When $X$ and $Y$ are independent, conditions on $Z$
can be checked on $X$ and $Y$ separately. For $(BL,\theta)$-dependence and association, this is
the following lemma.
\begin{lemma}
\label{lemma:indep}
If $\{X_s \in \mathbb{R}:s\in \mathbb{Z}^d\}$ and $\{Y_s \in \mathbb{R}:s\in \mathbb{Z}^d\}$ are
independent, then $\{X_s\}$ and $\{Y_s\}$ are both PA (NA, $(BL,\theta)$) if and only if the
bivariate random field $\{[X_s,Y_s]' \in \mathbb{R}^2, s \in \mathbb{Z}^d\}$ is PA (NA,
$(BL,\theta)$).
\end{lemma}
For summable covariance and exponential decay of covariances, the same holds
because the cross-covariances of $Z$ vanish. If $X$ and $Y$ are each $m$-dependent, then so is $Z$,
because for $\operatorname{dist}(I,J)>m$ the blocks $X_I$, $X_J$, $Y_I$ and $Y_J$ are mutually independent. If $X$
and $Y$ are each Gaussian, then so is $Z$. Hence Propositions~\ref{prop:bounded},
\ref{prop:gaussian} and \ref{prop:mdep} give Assumption~\ref{ass:bl} for independent bounded,
Gaussian and $m$-dependent fields, with their conditions checked on $X$ and $Y$ separately.

For independent fields that are unbounded and non-Gaussian, we verify the weaker
Assumption~\ref{ass:finite}. By Theorem~\ref{theorem:zd_rho0} and Corollary~\ref{cor:zd_indep},
Assumption~\ref{ass:finite} suffices for testing independence of $X$ and $Y$. It requires the product $\tilde X\tilde Y$ to be $(BL,\theta)$.
The following lemma shows that the product of two independent $(BL,\theta)$ fields is again
$(BL,\theta)$ when at least one of them is PA or NA.
\begin{lemma}[Products of independent weakly dependent fields]\label{lemma:product}

Let $X=\{X_s, s\in \mathbb{Z}^d\}$ and $Y=\{Y_s, s\in \mathbb{Z}^d\}$ be two independent random
fields.
Assume both are $(BL,\theta)$, at least one of them is PA or NA, and
\[
\sup_{s\in\mathbb{Z}^d} E[X_s^2] \le M_X < \infty,
\qquad
\sup_{s\in\mathbb{Z}^d} E[Y_s^2] \le M_Y < \infty.
\]
Define $Q_s = X_s Y_s$. Then $\{Q_s\}$ is $(BL,\theta)$.
\end{lemma}

\begin{remark}\label{remark:product}

Lemma~\ref{lemma:product} does not follow from Lemma~\ref{lemma:bl}, because the product
$(x,y)\mapsto xy$ is not Lipschitz on $\mathbb{R}^2$.
{For bounded fields, the result follows from Lemmas~\ref{lemma:indep} and
\ref{lemma:bounded_bl}, applied to the product component.}
Lemma~\ref{lemma:product} removes the boundedness condition at the cost of assuming that one
of the fields is PA or NA.
It is what makes Assumption~\ref{ass:finite} verifiable for independent unbounded fields, without Gaussianity (Proposition~\ref{prop:indep}).
\end{remark}

{Assumption~\ref{ass:finite} also requires the squared fields $\tilde X^2$ and $\tilde Y^2$ each to
have vanishing covariance.
Squares are not Lipschitz, so this does not follow from Lemma~\ref{lemma:bl}. The
following lemma shows that it follows from $(BL,\theta)$-dependence and finite fourth moments.}

\begin{lemma}[Squares of weakly dependent fields]\label{lemma:square_dying}
Let $V=\{V_s, s\in\mathbb{Z}^d\}$ be a strictly stationary real-valued $(BL,\theta)$ random field
with $E[V_0^4]<\infty$. Then the squared field $V^2=\{V_s^2, s\in\mathbb{Z}^d\}$ has vanishing
covariance.
\end{lemma}

Combining these results gives the following sufficient conditions, which require none of
boundedness, Gaussianity or $m$-dependence.
\begin{proposition}\label{prop:indep}
Under Assumptions~\ref{ass:zd_stationary} and \ref{ass:zd_squareint}, if $X$ and $Y$ are
independent, then Assumption~\ref{ass:finite} holds when any of the following holds:
\begin{enumerate}[label=(I\arabic*), ref=I\arabic*]
    \item \label{ass:I1} $X$ and $Y$ are both $(BL,\theta)$, and at least one of them is PA or NA.
    \item (stronger than I1) \label{ass:I2} $X$ and $Y$ are each PA or NA with summable
covariance.
    \item (stronger than I2) \label{ass:I3} $X$ and $Y$ are each PA or NA, and their
autocovariances decay exponentially as in \eqref{eq:exp_decay}.
\end{enumerate}
\end{proposition}

\section{Variance estimation}\label{sec:var}
Confidence intervals and tests based on Theorems~\ref{theorem:zd} and
\ref{theorem:zd_rho0} require an estimate of the asymptotic variance $\sigma^2$.
We assume $\sigma^2>0$.
Since $\sigma^2$ is the long-run variance of the field $L$ in \eqref{eq:variance_v}, we
estimate it with kernel estimators.
Such estimators for spatial data, known in econometrics
as spatial heteroskedasticity and autocorrelation consistent (HAC) estimators, go back to
\citet{conley1999gmm} for stationary mixing processes, and \citet{kelejian2007hac} extended
them to linear processes that may be heteroskedastic and locally nonstationary.
Using a common kernel form (Section~\ref{sec:var_general}), we propose five estimators that
differ only in how $h_k$ is estimated (Section~\ref{sec:var_special}).
One targets the general case, and the other four exploit the simpler forms of $\sigma^2$
derived in Section~\ref{sec:theorem} when the fields are Gaussian, when $\rho=0$, or when $X$ and
$Y$ are independent.
We also give a modification that guarantees nonnegative estimates (Remark~\ref{remark:nonneg}),
discuss the choice of bandwidth (Section~\ref{sec:var_bandwidth}), and describe an implementation
that applies in any fixed dimension and on domains that need not be rectangular
(Section~\ref{sec:var_impl}).

\subsection{{Common form}}\label{sec:var_general}
{Recall from \eqref{eq:variance} and \eqref{eq:variance_v} that
$\sigma^2=\sum_{k\in\mathbb{Z}^d}h_k$ with $h_k=\cov(L_0,L_k)$, the long-run variance of the scalar
field $L$.
For a field observed on a rectangular domain, \citet{prause2018estimation} estimate a long-run
covariance by a spatial Newey--West estimator with a product Bartlett kernel. We consider the same
form:
$$
\hat\sigma^2 = \sum_{|k_i|\leq m_i,\, i=1,\dots, d} w_{\bm{m}}(k)\,\hat h_k,
$$
where $\hat h_k$ is an estimate of $h_k$.
Here} $\bm{m}=[m_1, \dots, m_d]'\in \mathbb{Z}^{d}_{+}$
is the vector of truncation lags, and the weights are $w_{\bm{m}}(k) =
\prod_{i=1}^{d}\{1-|k_i|/(m_i+1)\}$, so that every lag with $|k_i|\le m_i$ for all $i$ receives a
positive weight.
The variance of $\hat\rho$ is then estimated by $\hat\sigma^2/|U_n|$.

\subsection{{Estimators}}\label{sec:var_special}
{Let $\hat X_s=(X_s-\bar X_n)/\hat\sigma_x$ and $\hat Y_s=(Y_s-\bar Y_n)/\hat\sigma_y$ be
standardized by the sample means and the sample standard deviations, where
$\hat\sigma_x^2=|U_n|^{-1}\sum_{s\in U_n}(X_s-\bar X_n)^2$ and
$\hat\sigma_y^2=|U_n|^{-1}\sum_{s\in U_n}(Y_s-\bar Y_n)^2$.}
{Without assuming Gaussianity or independence, we estimate $h_k=\cov(L_0,L_k)$ by the
sample autocovariance of $\hat L$, normalized by the number of pairs at lag
$k$~\citep{prause2018estimation},
\begin{gather*}
\hat h_k=\frac{1}{|A(k)|}\sum_{s\in A(k)}\hat L_s\hat L_{s+k},
\qquad
\hat L_s=\hat X_s\hat Y_s-\frac{\hat\rho}{2}\big(\hat X_s^2+\hat Y_s^2\big),\\
A(k)= \{ s \in U_n : s + k \in U_n\}.
\end{gather*}
No centering is needed, since $E[L_s]=0$.}
{This gives the estimator \rfgeneral{}, short for random field correlation.}

For Gaussian fields, Corollary~\ref{cor:zd_gaussian} gives {$h_k$ in
\eqref{eq:variance_gaussian}}, which depends only on the auto- and cross-correlations
at lag $k$.
The estimator \rfgaussian{} replaces these by the sample correlations
\begin{gather*}
\hat\rho_{xx,k} = \frac{1}{|A(k)|}\sum_{s\in A(k)}\hat X_s\hat X_{s+k},\qquad
\hat\rho_{yy,k} = \frac{1}{|A(k)|}\sum_{s\in A(k)}\hat Y_s\hat Y_{s+k},\\
\hat\rho_{xy,k} = \frac{1}{|A(k)|}\sum_{s\in A(k)}\hat X_s\hat Y_{s+k},
\end{gather*}
{and $\rho$ by $\hat\rho$.}
{For $d=1$, \rfgaussian{} is close to the estimator of \citet{afyouni2019effective},
which also plugs regularized sample auto- and cross-correlations into {\eqref{eq:variance_gaussian}}.}

{For testing $H_0:\rho=0$ or $H_0^{\text{indep}}: X\indep Y$, we estimate $h_k$ in
Theorem~\ref{theorem:zd_rho0} and Corollaries~\ref{cor:zd_rho0_gaussian} and \ref{cor:zd_indep}
in a similar way.}
This gives the three estimators for testing:
\begin{align*}
\text{\rfzero{}:}&\quad \hat\sigma^2 = \sum_{|k_i|\le m_i} w_{\bm{m}}(k)\,\frac{1}{|A(k)|}
\sum_{s\in A(k)}\hat X_s\hat Y_s\hat X_{s+k}\hat Y_{s+k},\\
\text{\rfzerogaussian{}:}&\quad \hat\sigma^2 = \sum_{|k_i|\le m_i}
w_{\bm{m}}(k)\big(\hat\rho_{xx,k}\hat\rho_{yy,k}+\hat\rho_{xy,k}\hat\rho_{xy,-k}\big)
{\ \text{ with } \hat\rho_{xy,0}=0},\\
\text{\rfindep{}:}&\quad \hat\sigma^2 = \sum_{|k_i|\le m_i}
w_{\bm{m}}(k)\,\hat\rho_{xx,k}\hat\rho_{yy,k}.
\end{align*}
Table~\ref{tab:methods} summarizes the five estimators, the asymptotic-variance formula each one
estimates, and the inference task for which it is intended.

\begin{remark}\label{remark:nonneg}
The normalization $1/|A(k)|$ reduces the bias of each {$\hat h_k$} but does not guarantee
$\hat\sigma^2\ge0$.
Replacing $1/|A(k)|$ by $1/|U_n|$ {throughout} guarantees $\hat\sigma^2\ge0$ for \rfgeneral{}, \rfzero{} and \rfindep{}
(Lemma~\ref{lemma:nonneg} in Appendix~\ref{append:nonneg}), at the cost of a downward bias at
large lags. In the rare case that $\hat\sigma^2<0$ for \rfgeneral{}, \rfzero{} or \rfindep{}, we recompute it
with the normalization $|U_n|$.
Lemma~\ref{lemma:nonneg} does not cover \rfgaussian{} and \rfzerogaussian{}, so if either gives
$\hat\sigma^2<0$, we use the $|U_n|$ version of \rfgeneral{} or \rfzero{}, respectively.
In our simulations, the fallback was never needed, since no estimate was negative.
\end{remark}

\subsection{Choice of bandwidth}\label{sec:var_bandwidth}
The choice of the bandwidth $\bm{m}$ is crucial to the estimator's performance.
For the Bartlett kernel, the bandwidth minimizing the asymptotic truncated MSE of
{$\hat\sigma^2$} grows as $|U_n|^{1/(d+2)}$ \citep{kim2011spatial}.
For $d=1$, this recovers the classical cube-root rate
\citep{andrews1991heteroskedasticity}.
Note that subsampling variance estimators, which for the sample mean are asymptotically
equivalent to Bartlett-kernel estimators, have the same MSE-optimal rate for the block side length
on rectangular domains \citep{politis1994large}.
We therefore use the same truncation lag in every direction, $m_i = \lceil \kappa\,|U_n|^{1/(d+2)}\rceil$
for $i=1,\dots,d$, with $\kappa=1$ throughout (to be discussed in Section \ref{sec:discussion}).
Note that when $U_n$ is a rectangle with side lengths $n_1, \dots, n_d$, $|U_n|^{1/(d+2)}$ is the geometric
mean of the direction-specific lags $n_1^{d/(d+2)}, \dots, n_d^{d/(d+2)}$.

\subsection{Implementation}\label{sec:var_impl}
The estimators apply in any dimension $d$ and to any finite domain $U_n\subset\mathbb{Z}^d$,
not only to rectangles.
To compute them, we embed $U_n$ in a rectangular grid $G_n\supseteq U_n$, such as its bounding box,
and set the standardized fields to zero at grid points outside $U_n$.
{Then $\hat L$ is also zero outside $U_n$, and every sum over $A(k)$ in
Section~\ref{sec:var_special} has the form $\sum_{s\in G_n}a_sb_{s+k}$ with
$a,b\in\{\hat L,\hat X\hat Y,\hat X,\hat Y\}$.
Similarly,} { $|A(k)|=\sum_{s\in G_n}\mathbbm{1}_{U_n}(s)\,\mathbbm{1}_{U_n}(s+k)$, where
$\mathbbm{1}_{U_n}$ is the indicator of $U_n$.
Such lagged sums can be computed for all lags $k$ at once, for example with the Python function
\texttt{scipy.signal.correlate}.}
Lags with $A(k)=\emptyset$, which can occur for non-rectangular domains, are omitted from
{$\hat\sigma^2$}.
The methods are implemented in the Python package \texttt{RFCorr} \citep{RFCorr}.

\section{Simulations}\label{sec:simu}
We conduct simulations to evaluate the performance of the various variance estimation methods on two
tasks: interval estimation and hypothesis testing.
For the interval estimation task, we evaluate the coverage rate of the confidence interval.
For the testing task, we evaluate the Type~I error rate of the test for the null
$H_0:\rho=0$ and for the null $H_0^{\text{indep}}: X\indep Y$.
{Since independent fields satisfy both nulls, we also include two designs in which $X$ and
$Y$ are uncorrelated but dependent, so that $H_0$ holds and $H_0^{\text{indep}}$ does not.}

\subsection{Simulation design}\label{sec:simu_design}
We simulate fields in $d=2$ observed on an elliptical domain, mimicking images that come with
non-rectangular masks.
For a side length $n\in\{25,35,50,70,100,140,200,300,420,600\}$, the domain $U_n$ is the set
of points of the $n\times n$ grid that lie in the ellipse centered in the grid with semi-axes $n/2$
and $0.8\,n/2$ along the two grid axes.
Consecutive side lengths differ by a factor of about $\sqrt2$, so that $N$ roughly
doubles from one to the next.
The number of points $N=|U_n|$ ranges from $389$ to $226{,}204$.
The fields are generated on the full grid and only the points in $U_n$ are used.

We first generate two independent stationary Gaussian fields $W_1$ and $W_2$ by smoothing i.i.d.\
$\mathcal{N}(0,1)$ noise with a Gaussian kernel, and rescale them to
unit variance.
To ensure stationarity and avoid edge effects, the noise is drawn on an enlarged grid
and cropped to the original grid after smoothing.
The kernel for $W_1$ has full width at half maximum (FWHM) of $1$, $2$, $3$ or $4$ pixels, and the kernel
for $W_2$ has half that FWHM, so a larger FWHM gives smoother fields with stronger autocorrelation.
{In the baseline design, $X$ and $Y$ depend on each other only through $\rho$.}
In the Gaussian setting we set $X=W_1$ and $Y=\rho W_1+\sqrt{1-\rho^2}\,W_2$.
In the $\chi^2_1$ setting we set $X=W_1^2$ and $Y=\rho W_1^2+\sqrt{1-\rho^2}\,W_2^2$, which are
non-Gaussian but still have correlation $\rho$, since $W_1^2$ and $W_2^2$ are independent with equal
variance.

{To separate the two nulls, we also use two designs with $\rho=0$ in which $Y$ depends on
$X$.
In the sum-difference design, the kernel for $W_2$ has a quarter of the FWHM of the kernel for
$W_1$, and we set $X=(W_1+W_2)/\sqrt2$ and $Y=(W_1-W_2)/\sqrt2$.
In the $\chi^2_1$ setting, $W_1$ and $W_2$ are replaced by $W_1^2$ and $W_2^2$.
Since the two summands are independent with equal variance, $X$ and $Y$ have equal variance and
$\cov(X_s,Y_t)=\{\cov(W_{1,s},W_{1,t})-\cov(W_{2,s},W_{2,t})\}/2$
for all $s,t$.
This is $0$ for $s=t$, so $\rho=0$, but it is positive for $s\neq t$ because $W_1$ is smoother
than $W_2$, so $X$ and $Y$ are cross-correlated at nonzero lags.
In the shift design, we instead set $X=W_1$ and define $Y$ by $Y_s=(W_{1,s+e_1}+W_{1,s-e_1})/2-a\,W_{1,s}$, where $e_1=(1,0)$ and $a$ is the
population lag-one autocorrelation $\operatorname{Corr}(W_{1,s},W_{1,s+e_1})$.
In the $\chi^2_1$ setting, $W_1$ is replaced by $W_1^2$, whose lag-one autocorrelation is $a^2$.
In both settings, $Y_s=(X_{s+e_1}+X_{s-e_1})/2-c\,X_s$, where $r(k)=\operatorname{Corr}(X_s,X_{s+k})$ and $c=r(e_1)$, which is $a$ in the Gaussian setting
and $a^2$ in the $\chi^2_1$ setting. Then
$$
\cov(X_s,Y_t)=\operatorname{Var}(X_0)\Big\{\frac{r(k+e_1)+r(k-e_1)}{2}-c\,r(k)\Big\},
\qquad k=t-s.
$$
This is $0$ for $t=s$, since $r(\pm e_1)=c$, so $\rho=0$, but it is generally nonzero for $t\neq s$,
so $X$ and $Y$ are dependent.}

{In the baseline design, we use
$\rho\in\{0,\pm0.1,\pm0.2,\pm0.4,\pm0.6,\pm0.8,\pm0.9\}$.
For each setting we run $1000$ replicates.
{The gray dashed lines in the figures mark the band} $p\pm1.96\sqrt{p(1-p)/1000}$, where $p$ is the nominal
level ($0.95$ for coverage and $0.05$ for Type~I error).}
All kernel-based methods use the Bartlett kernel with the same truncation lag
$m=\lceil |U_n|^{1/4}\rceil$ in both directions, the rule of Section~\ref{sec:var_bandwidth} for
$d=2$.

To improve the normal approximation for the correlation estimate, we apply Fisher's transformation, $F_{XY} =
\operatorname{arctanh}(\hat{\rho})$.
By the delta method, $F_{XY}$ has approximate variance
\begin{equation}\label{eq:fisher}
    \operatorname{Var}(F_{XY}) = \operatorname{Var}(\hat{\rho})\,(1 - \rho^2)^{-2}.
\end{equation}
In practice we replace $\rho$ by $\hat\rho$ for confidence intervals, and by $0$ for tests.

For comparison, in both tasks we consider a \texttt{naive} method that ignores
autocorrelation in the random fields, setting $\widehat{\operatorname{Var}}(F_{XY}) = 1/|U_n|$.
This is motivated by the result that, in the absence of autocorrelation and under Gaussianity,
$\sqrt{|U_n|}\,\hat{\rho}$ has asymptotic variance $(1-\rho^2)^2$, so that
$\operatorname{Var}(F_{XY}) = 1/|U_n|$ after the Fisher transformation.

\subsection{Interval estimation of the correlation}

We construct $95\%$ confidence intervals for $F_{XY}$ using a normal approximation, and then map the
endpoints back to the correlation scale via the inverse transformation $\tanh(\cdot)$.
To estimate $\operatorname{Var}(F_{XY})$ via \eqref{eq:fisher}, we first estimate
$\operatorname{Var}(\hat{\rho})$ using the two methods \rfgeneral{} and \rfgaussian{} described in
Section~\ref{sec:var}.

\begin{figure}[htbp]
    \centering
    \includegraphics[width=\textwidth]{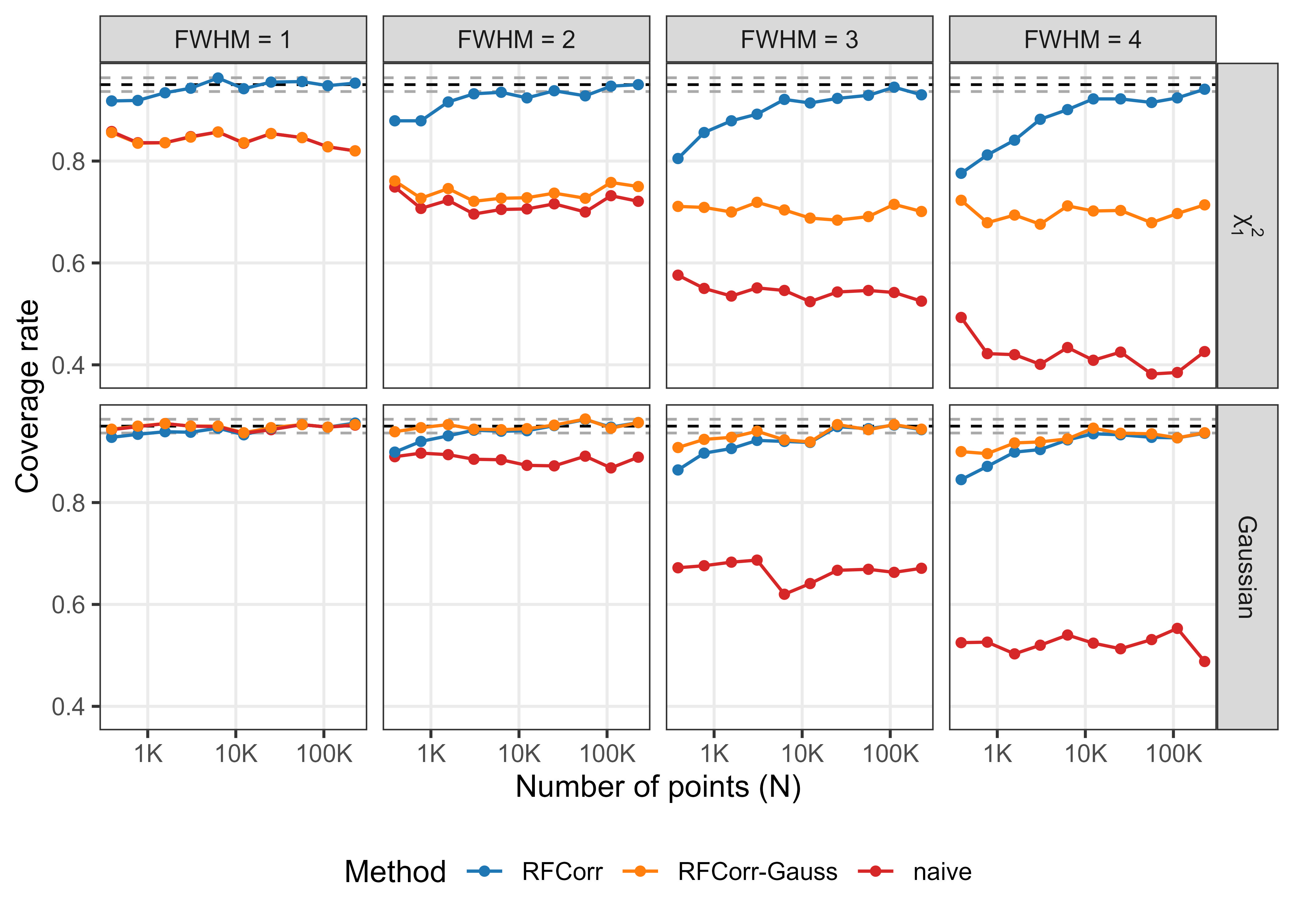}
    \caption{Coverage rate of $95\%$ confidence intervals for $\rho=0.4$ against the
    number of points $N$ (log scale) under the baseline design, for $\chi^2_1$ fields (top) and
    Gaussian fields (bottom) in $d=2$.
    Columns give the FWHM of the smoothing kernel, so fields are smoother from left to right.
    The black dashed line marks the nominal level $0.95$ and the gray dashed lines mark the range
    expected from Monte Carlo error alone, a pointwise $95\%$ binomial band for a method at the
    nominal level with $1000$ replicates (Section~\ref{sec:simu_design}).}
    \label{fig:coverage}
\end{figure}

Figure~\ref{fig:coverage} shows the coverage rates for $\rho=0.4$.
{Results for the other values of $\rho$ (Appendix~\ref{append:simu}) were similar, except
that \rfgeneral{} undercovered more on small domains when $|\rho|$ was large, and the undercoverage
of \rfgaussian{} for $\chi^2_1$ fields grew with $|\rho|$.}
For Gaussian fields, \rfgeneral{} and \rfgaussian{} both approached the nominal level as $N$ grew.
On small domains with smooth fields, \rfgeneral{} undercovered, whereas \rfgaussian{}
stayed closer to the nominal level, since it does not estimate fourth-order moments.
For $\chi^2_1$ fields, only \rfgeneral{} approached the nominal level as $N$ grew.
\rfgaussian{} undercovered at every sample size, because its variance formula assumes Gaussian
fourth-order moments.
The naive method undercovered whenever the fields were autocorrelated or non-Gaussian, increasingly
so for smoother fields, and did not improve as $N$ grew.
For non-Gaussian fields, its variance $(1-\rho^2)^2$ is incorrect even without autocorrelation
(Remark~\ref{remark:zd_iid}).

\subsection{Testing zero correlation and independence}
For hypothesis testing, we compute the $z$-score $Z =
F_{XY}/\sqrt{\widehat{\operatorname{Var}}(F_{XY})}$.
As for interval estimation, the key is to estimate $\operatorname{Var}(F_{XY})$.
Since both $H_0$ and $H_0^{\text{indep}}$ imply $\rho=0$, we plug in $\rho=0$.
For the null $H_0:\rho=0$, we estimate $\operatorname{Var}(\hat{\rho})$ using the two methods
\rfzero{} and \rfzerogaussian{} described in Section~\ref{sec:var}.
For the null $H_0^{\text{indep}}: X\indep Y$, we use \rfindep{}, also described in
Section~\ref{sec:var}. \rfindep{} estimates the same asymptotic variance as the
approximations of \citet{clifford1989assessing} and \citet{dutilleul1993modifying} underlying the
modified $t$-test (Theorem~\ref{theorem:clifford}), even though those were derived for Gaussian
fields.

\begin{figure}[htbp]
    \centering
    \includegraphics[width=\textwidth]{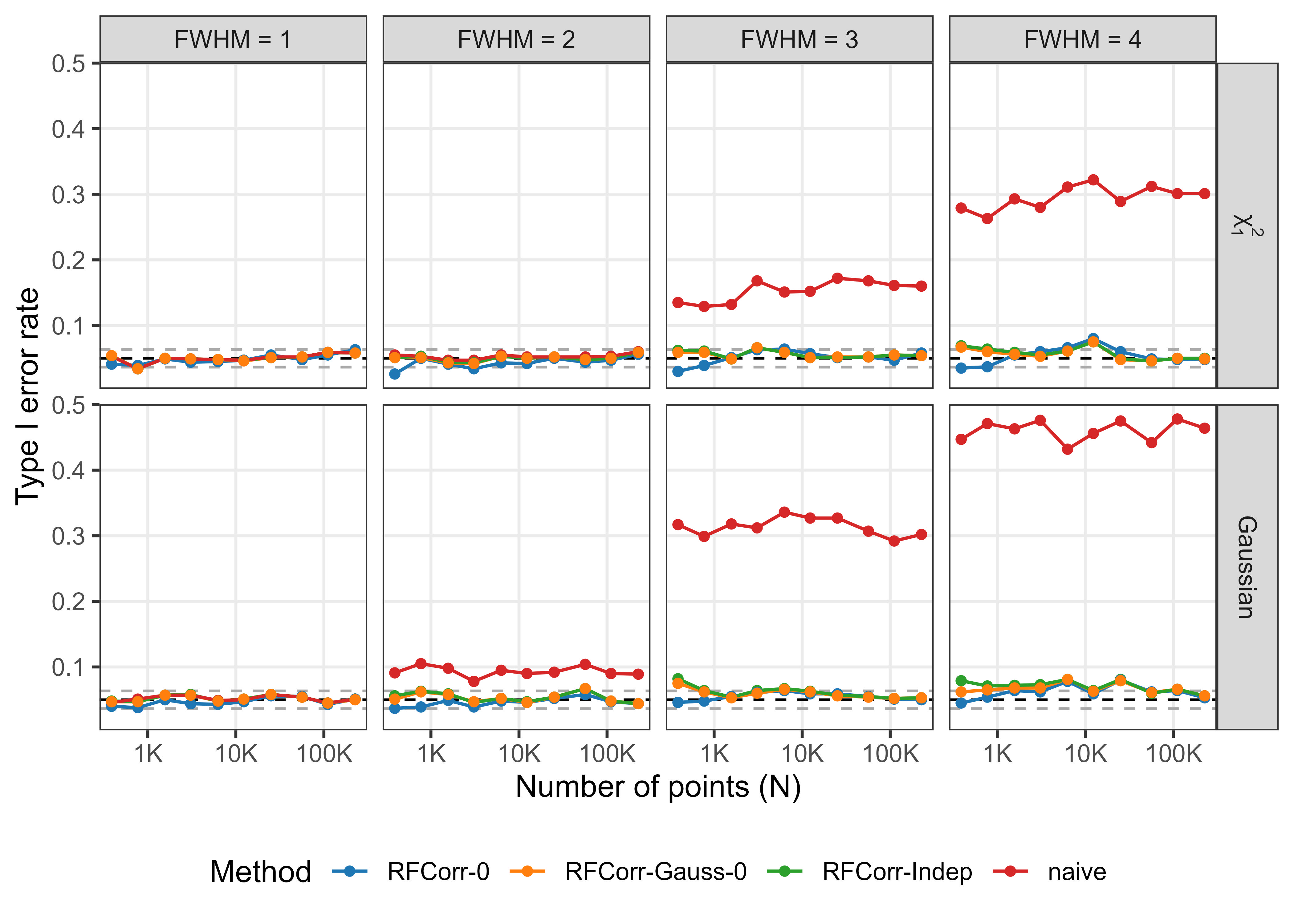}
    \caption{Type~I error rate of level-$0.05$ tests against the number of points $N$
    (log scale) under the baseline design with $\rho=0$, for independent $\chi^2_1$ fields (top)
    and Gaussian fields (bottom) in $d=2$.
    \rfzero{}, \rfzerogaussian{} and the naive method test $H_0:\rho=0$, and \rfindep{} tests
    $H_0^{\text{indep}}: X\indep Y$.
    Columns and dashed lines are as in Figure~\ref{fig:coverage}, with the nominal level $0.05$.}
    \label{fig:type1}
\end{figure}

Figure~\ref{fig:type1} shows the Type~I error rates {under the baseline design}.
{In this design}, $X$ and $Y$ are independent when $\rho=0$, so both $H_0$
and $H_0^{\text{indep}}$ hold, and \rfzero{}, \rfzerogaussian{} and \rfindep{} estimate the same
asymptotic variance.
All three tests kept the Type~I error rate close to $0.05$, with mild inflation for smooth
Gaussian fields.
The naive test had inflated Type~I error rates that grew with the smoothness of the fields and did
not decrease as $N$ grew.

\begin{figure}[htbp]
    \centering
    \includegraphics[width=\textwidth]{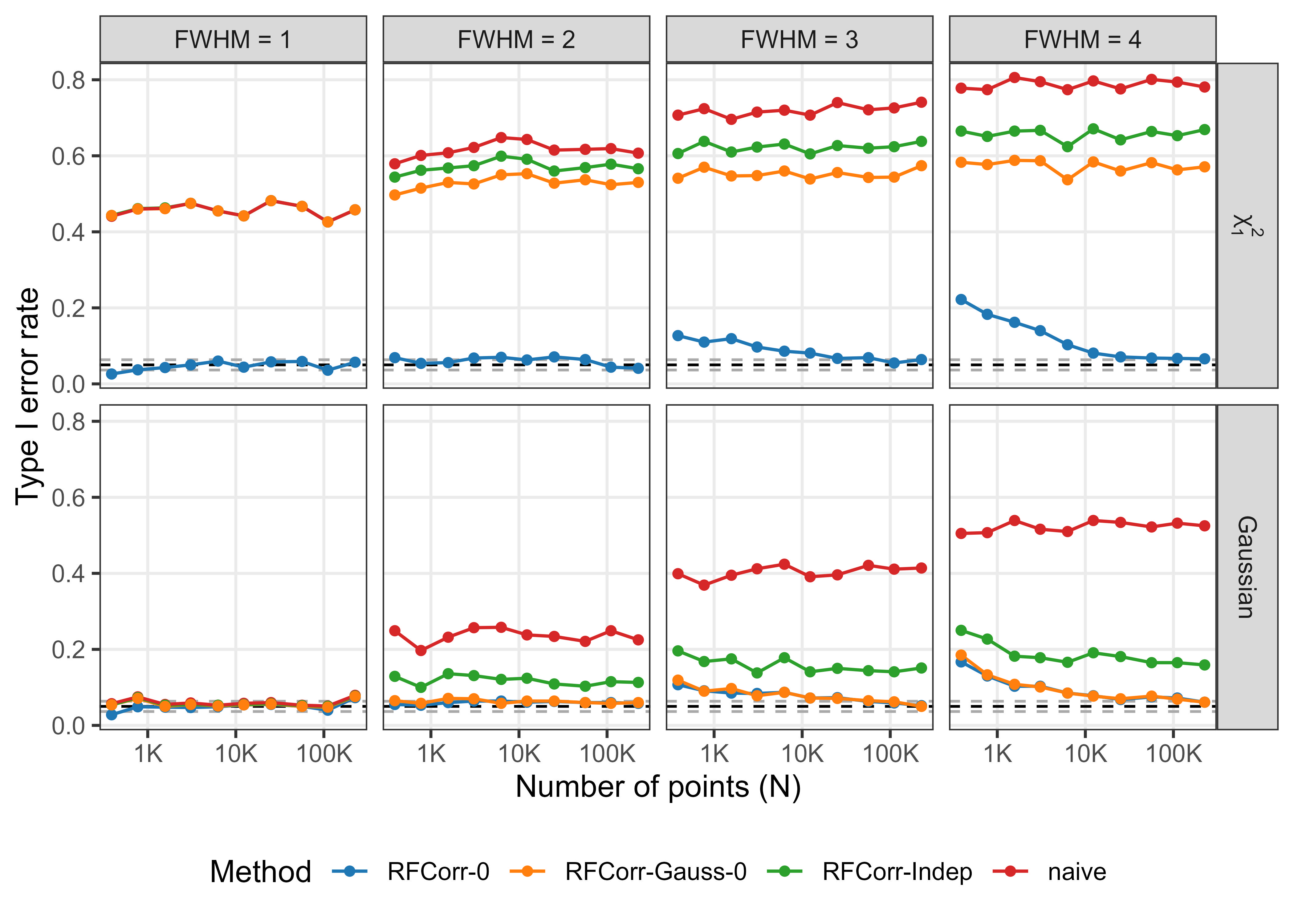}
    \caption{Rate of rejecting $H_0:\rho=0$ at level $0.05$ against the number of
    points $N$ (log scale) under the sum-difference design of Section~\ref{sec:simu_design}, where
    $X$ and $Y$ are uncorrelated but dependent, for $\chi^2_1$ fields (top) and Gaussian fields
    (bottom) in $d=2$.
    $H_0$ holds, so the rates of \rfzero{}, \rfzerogaussian{} and the naive method are Type~I error
    rates, whereas $H_0^{\text{indep}}$, tested by \rfindep{}, is false.
    {For $\chi^2_1$ fields with FWHM $1$, the line of \rfindep{} is hidden behind those of
    \rfzerogaussian{} and the naive method, since the three rates nearly coincide.}
    Columns and dashed lines are as in Figure~\ref{fig:type1}.}
    \label{fig:type1_sumdiff}
\end{figure}

{Figure~\ref{fig:type1_sumdiff} shows the rates of rejecting $H_0:\rho=0$ under the
sum-difference design, where $H_0$ holds but $H_0^{\text{indep}}$ does not.
The Type~I error rate of \rfzero{} approached $0.05$ as $N$ grew, although slowly for smooth
fields.
\rfzerogaussian{} did so only for Gaussian fields.
\rfindep{} rejected $H_0$ more often than $5\%$, at rates that did not decrease with $N$.
This is expected since it tests $H_0^{\text{indep}}$, which is false here.
The exception was Gaussian fields with FWHM $1$, where the cross-correlation of $X$ and $Y$ is
weak.
Like \rfindep{}, the modified $t$-test of \citet{clifford1989assessing} tests
independence, so it should not be used to test zero correlation when the fields may be dependent.
The naive test was the most inflated, except for Gaussian fields with FWHM $1$, which are nearly
free of autocorrelation.
The shift design gave similar results (Figure~\ref{fig:type1_symmetric} in
Appendix~\ref{append:simu}).}

\section{Discussion}\label{sec:discussion}

We studied the asymptotic properties of the Pearson correlation coefficient between
two random fields on
$\mathbb{Z}^d$ under $(BL,\theta)$-dependence, in any fixed dimension.
Since $\hat\rho$ is CAN with an explicit asymptotic variance, which our kernel-based estimators
estimate, one can construct confidence intervals for $\rho$, rather than only test whether
$\rho=0$, without Gaussianity or a parametric model for the spatial dependence.
For independent fields, the variance approximations of \citet{clifford1989assessing} and
\citet{dutilleul1993modifying}, which underlie the modified $t$-test, converge to our asymptotic
variance when computed from the true covariance matrices of $X$ and $Y$
(Theorem~\ref{theorem:clifford}), even though they were derived for Gaussian fields.

A strength of our theory is that, for broad classes of fields, the sufficient conditions for the
asymptotic normality of $\hat\rho$ are simple to state and to check.
For Gaussian fields, absolutely summable auto- and cross-covariances suffice, with no restriction
on the dependence between $X$ and $Y$.
For bounded fields that are positively or negatively associated, the same covariance condition
suffices.
These covariance conditions hold, for example, whenever the covariances decay exponentially
(Lemma~\ref{lemma:cox}).
For $m$-dependent fields, finite fourth moments suffice.
For independent fields, it suffices that both are $(BL,\theta)$ and at least one is positively or
negatively associated, without boundedness or Gaussianity.

The explicit form of the asymptotic variance $\sigma^2$ sheds light on what determines the
precision of $\hat\rho$.
In general, $\sigma^2$ depends on the joint fourth moments of the standardized fields, and for
Gaussian fields only on their auto- and cross-correlations.
For fields obtained by smoothing Gaussian white noise with a Gaussian kernel of standard deviation
$\ell$ (Example~\ref{example:smoothed}), $\operatorname{Var}(\hat\rho)$ is approximately
$(1-\rho^2)^2(\sqrt{2\pi}\,\ell)^d/N$, so the precision depends only on the number of grid points, the
correlation, the smoothness and the dimension.
The correlation enters through $(1-\rho^2)^2$, as for i.i.d.\ data, and the smoothness and
dimension through the inflation factor $(\sqrt{2\pi}\,\ell)^d$.
Because this factor grows exponentially in $d$, smooth fields in three dimensions, such as brain
activation maps, need many more grid points than one- or two-dimensional data for the same
precision.
The approximation can also be used at the planning stage to choose the number of grid points $N$
needed for a target standard error.

Our theory does not require rectangular domains.
Assumption~\ref{ass:regular} only requires $U_n$ to be regularly growing, which holds for many
non-rectangular shapes, such as brain masks \citep{qin2025simultaneous} or land areas \citep{sommerfeld2018confidence}, as long as their boundary is small relative
to their volume.
The variance estimators can likewise be computed on non-rectangular domains
(Section~\ref{sec:var_impl}).
In our simulation study, the fields were observed on elliptical domains, with the bandwidth
computed from $|U_n|$ as in Section~\ref{sec:var_bandwidth}, and both interval estimation of
$\rho$ and the tests of $\rho=0$ and of independence performed well.

The simulations also guide the choice of variance estimator, which should follow the question of
interest.
For a confidence interval for $\rho$, we recommend \rfgeneral{}, and for testing $H_0:\rho=0$,
\rfzero{}. When the fields are plausibly Gaussian, \rfgaussian{} and \rfzerogaussian{} can be used
instead. For testing
$H_0^{\text{indep}}: X\indep Y$, we recommend \rfindep{}, which estimates the same asymptotic variance as the
approximations of \citet{clifford1989assessing} and \citet{dutilleul1993modifying}.
The naive variance $1/|U_n|$ should be avoided. In our simulations, its confidence intervals
undercovered and its tests had inflated Type~I error, increasingly so for smoother fields, whereas
\rfgeneral{} and the three tests stayed close to the nominal levels for moderate to
large $|U_n|$, each under its own null hypothesis.
\rfgaussian{} undercovered substantially for $\chi^2_1$ fields at every sample size, so it should not be
used for clearly non-Gaussian fields.
For smooth fields on small domains, \rfgeneral{} undercovered somewhat, presumably because the
fourth-order moments in its variance are hard to estimate from few effectively independent
observations, whereas \rfgaussian{} stayed closer to the nominal level for Gaussian
fields.
More generally, smoother fields needed larger domains to reach the nominal coverage
and Type~I error rate.

Several questions remain open.
First, our bandwidth rule $m_i = \lceil \kappa\,|U_n|^{1/(d+2)}\rceil$ fixes the rate but not the
constant $\kappa$.
In our simulations, $\kappa=1$ performed well.
This agrees with the observation of \citet{andrews1991heteroskedasticity} and
\citet{kim2011spatial} that the variance estimator performs well as long as the bandwidth is close
to its optimal value, without needing to equal it exactly.
Data-driven rules for choosing the bandwidth are an alternative, and they have been developed for time
series \citep{andrews1991heteroskedasticity, newey1994automatic} and for spatial data
\citep{kim2011spatial}.
These rules require estimating weighted sums of autocovariances that are hard to estimate
nonparametrically \citep{kim2011spatial}.
In practice, we recommend checking how sensitive the conclusions are to the choice of $\kappa$. For example, one can recompute $\hat\sigma^2$ over a range such as $\kappa\in\{0.5,1,2,3\}$
and examine how it changes with $m$.
A value of $m$ at which $\hat\sigma^2$ has leveled off is preferable to one at which it is still
changing. The choice can also be calibrated by simulation. For example, one generates fields with the same domain and
similar smoothness as the data and checks the coverage of the resulting confidence intervals for several values of
$\kappa$.

Second, proving that the variance estimators in Section~\ref{sec:var} are consistent under the
$(BL,\theta)$-dependence assumed in Section~\ref{sec:theorem} is substantial work beyond the scope
of this paper, so we assessed their finite-sample performance in Section~\ref{sec:simu} instead.
Theorem~\ref{theorem:clifford} concerns the Clifford and Dutilleul approximations computed
from the true covariance matrices of $X$ and $Y$, whereas in practice these are estimated from the data, so the
consistency of the modified $t$-test is also open.
Existing consistency results for spatial kernel estimators \citep{kim2011spatial,
prause2018estimation} rely on mixing conditions and higher-order moment conditions that differ
from ours.
We expect that a proof under $(BL,\theta)$-dependence would need additional assumptions, such as
summable fourth-order {cumulants of $L$}.

Third, our theory could be extended in two directions.
{One is to fields indexed by a continuous domain, observed on increasingly fine lattices
over growing domains, for which the CLT of \citet{Bulinski2010central} could be a starting point.}
The other direction is to nonstationary fields, whose means or variances vary over space.
Central limit theorems for nonstationary arrays of random fields on irregular lattices
\citep{jenish2009central} could serve as a starting point, although $\rho$ would then have to be
interpreted as a spatially averaged correlation.

Finally, our results rely on increasing-domain asymptotics, in which the domain grows while the
spacing between observations stays fixed.
In practice, the domain of an image is fixed, and a higher resolution adds observations within
it, which is closer to infill asymptotics.
Under infill asymptotics, estimators of dependence parameters can behave very differently and
may even be inconsistent \citep{zhang2004inconsistent}.
Our results therefore give a good approximation when the range of spatial dependence is small
relative to the size of the domain, which is typical of large, high-resolution images.
Studying $\hat\rho$ under infill asymptotics is an interesting direction for future work.

\begin{funding}
J.Q., S.D. and A.S. were partially supported by NIH grant R01MH128923.
\end{funding}

\section*{Data and code availability}
The proposed methods are implemented in the Python package \texttt{RFCorr} \citep{RFCorr},
available at \url{https://github.com/JiyueQin/RFcorr}, together with the code to reproduce the
simulations and figures. This paper uses only simulated data, generated as described in
Section~\ref{sec:simu_design}.

\clearpage
\begin{table}[ht]
\centering
\setlength{\belowcaptionskip}{8pt}
\caption{Variance estimation methods by case, with the corresponding asymptotic-variance formula
and the inference task for which each is intended.}
\label{tab:methods}
\begin{tabular}{llll}
\hline
Case & Method name & Formula & Used for \\
\hline
General & \rfgeneral & \eqref{eq:variance} & CI for $\rho$ \\
Gaussian & \rfgaussian & \eqref{eq:variance_gaussian} & CI for $\rho$ \\
$\rho=0$ & \rfzero & \eqref{eq:variance_0} & test of $\rho=0$ \\
$\rho=0$, Gaussian & \rfzerogaussian & \eqref{eq:variance_gaussian_0} & test of $\rho=0$ \\
$X\indep Y$ & \rfindep & \eqref{eq:variance_indep} & test of independence \\
\hline
\end{tabular}
\end{table}

\clearpage
\bibliographystyle{apalike}
\bibliography{references}

\begin{thebibliography}{}

\bibitem[Adler and Taylor, 2007]{adler2007random}
Adler, R.~J. and Taylor, J.~E. (2007).
\newblock {\em Random Fields and Geometry}.
\newblock Springer.

\bibitem[Afyouni et~al., 2019]{afyouni2019effective}
Afyouni, S., Smith, S.~M., and Nichols, T.~E. (2019).
\newblock Effective degrees of freedom of the {P}earson's correlation
  coefficient under autocorrelation.
\newblock {\em NeuroImage}, 199:609--625.

\bibitem[Alexander-Bloch et~al., 2018]{alexander2018testing}
Alexander-Bloch, A.~F., Shou, H., Liu, S., Satterthwaite, T.~D., Glahn, D.~C.,
  Shinohara, R.~T., Vandekar, S.~N., and Raznahan, A. (2018).
\newblock On testing for spatial correspondence between maps of human brain
  structure and function.
\newblock {\em NeuroImage}, 178:540--551.

\bibitem[Anderson, 2003]{anderson2003introduction}
Anderson, T.~W. (2003).
\newblock {\em An Introduction to Multivariate Statistical Analysis}.
\newblock Wiley-Interscience, Hoboken, NJ, 3rd edition.

\bibitem[Andrews, 1991]{andrews1991heteroskedasticity}
Andrews, D. W.~K. (1991).
\newblock Heteroskedasticity and autocorrelation consistent covariance matrix
  estimation.
\newblock {\em Econometrica}, 59(3):817--858.

\bibitem[Bandeira et~al., 2023]{bandeira2023matrix}
Bandeira, A.~S., Boedihardjo, M.~T., and van Handel, R. (2023).
\newblock Matrix concentration inequalities and free probability.
\newblock {\em Inventiones Mathematicae}, 234:419--487.

\bibitem[Bartlett, 1935]{bartlett1935some}
Bartlett, M.~S. (1935).
\newblock Some aspects of the time-correlation problem in regard to tests of
  significance.
\newblock {\em Journal of the Royal Statistical Society}, 98(3):536--543.

\bibitem[Bernstein, 1918]{Berstein1918law}
Bernstein, S.~N. (1918).
\newblock On the law of large numbers.
\newblock {\em Communications of Kharkov. Mathematical Society. Second Series},
  16(1--2):82--87.
\newblock In Russian.

\bibitem[Bhattacharya et~al., 2025]{bhattacharya2025nonsense}
Bhattacharya, S., Mukherjee, R., and Ogburn, E.~L. (2025).
\newblock Nonsense associations in {M}arkov random fields with pairwise
  dependence.
\newblock {\em Biometrika}, 112(4).
\newblock Article asaf041.

\bibitem[Bradley, 2005]{bradley2005basic}
Bradley, R.~C. (2005).
\newblock Basic properties of strong mixing conditions: {A} survey and some
  open questions.
\newblock {\em Probability Surveys}, 2:107--144.

\bibitem[Brockwell and Davis, 1991]{brockwell1991time}
Brockwell, P.~J. and Davis, R.~A. (1991).
\newblock {\em Time Series: Theory and Methods}.
\newblock Springer.

\bibitem[Bulinski, 2010]{Bulinski2010central}
Bulinski, A. (2010).
\newblock Central limit theorem for random fields and applications.
\newblock In Skiadas, C.~H., editor, {\em Advances in Data Analysis: Theory and
  Applications to Reliability and Inference, Data Mining, Bioinformatics,
  Lifetime Data, and Neural Networks}, pages 141--150. Birkh{\"a}user, Boston.

\bibitem[Bulinski and Shashkin, 2007]{bulinski2007limit}
Bulinski, A. and Shashkin, A. (2007).
\newblock {\em Limit Theorems for Associated Random Fields and Related
  Systems}, volume~10 of {\em Advanced Series on Statistical Science and
  Applied Probability}.
\newblock World Scientific.

\bibitem[Bulinski and Suquet, 2001]{bulinski2001normal}
Bulinski, A. and Suquet, C. (2001).
\newblock Normal approximation for quasi-associated random fields.
\newblock {\em Statistics \& Probability Letters}, 54(2):215--226.

\bibitem[Burton et~al., 1986]{burton1986invariance}
Burton, R.~M., Dabrowski, A., and Dehling, H. (1986).
\newblock An invariance principle for weakly associated random vectors.
\newblock {\em Stochastic Processes and their Applications}, 23(2):301--306.

\bibitem[Cacoullos, 1989]{cacoullos1989exercises}
Cacoullos, T. (1989).
\newblock {\em Exercises in Probability}.
\newblock Springer.

\bibitem[Campbell et~al., 1997]{campbell1998econometrics}
Campbell, J.~Y., Lo, A.~W., and MacKinlay, A.~C. (1997).
\newblock {\em The Econometrics of Financial Markets}.
\newblock Princeton University Press, Princeton, NJ.

\bibitem[Clifford et~al., 1989]{clifford1989assessing}
Clifford, P., Richardson, S., and H{\'e}mon, D. (1989).
\newblock Assessing the significance of the correlation between two spatial
  processes.
\newblock {\em Biometrics}, 45(1):123--134.

\bibitem[Conley, 1999]{conley1999gmm}
Conley, T.~G. (1999).
\newblock {GMM} estimation with cross sectional dependence.
\newblock {\em Journal of Econometrics}, 92(1):1--45.

\bibitem[Doukhan and Louhichi, 1999]{doukhan1999new}
Doukhan, P. and Louhichi, S. (1999).
\newblock A new weak dependence condition and applications to moment
  inequalities.
\newblock {\em Stochastic Processes and their Applications}, 84(2):313--342.

\bibitem[Durrett, 2019]{durrett2019probability}
Durrett, R. (2019).
\newblock {\em Probability: Theory and Examples}, volume~49 of {\em Cambridge
  Series in Statistical and Probabilistic Mathematics}.
\newblock Cambridge University Press.

\bibitem[Dutilleul et~al., 1993]{dutilleul1993modifying}
Dutilleul, P., Clifford, P., Richardson, S., and H{\'e}mon, D. (1993).
\newblock Modifying the {$t$} test for assessing the correlation between two
  spatial processes.
\newblock {\em Biometrics}, 49(1):305--314.

\bibitem[Folland, 1999]{folland1999real}
Folland, G.~B. (1999).
\newblock {\em Real Analysis: Modern Techniques and Their Applications}.
\newblock Wiley, New York, 2nd edition.

\bibitem[Fox et~al., 2005]{fox2005human}
Fox, M.~D., Snyder, A.~Z., Vincent, J.~L., Corbetta, M., Van~Essen, D.~C., and
  Raichle, M.~E. (2005).
\newblock The human brain is intrinsically organized into dynamic,
  anticorrelated functional networks.
\newblock {\em Proceedings of the National Academy of Sciences},
  102(27):9673--9678.

\bibitem[Hannan, 1976]{hannan1976asymptotic}
Hannan, E.~J. (1976).
\newblock The asymptotic distribution of serial covariances.
\newblock {\em The Annals of Statistics}, 4(2):396--399.

\bibitem[Hansen, 2022]{hansen2022econometrics}
Hansen, B. (2022).
\newblock {\em Econometrics}.
\newblock Princeton University Press.

\bibitem[H{\"o}rmann and Kokoszka, 2013]{hormann2013consistency}
H{\"o}rmann, S. and Kokoszka, P. (2013).
\newblock Consistency of the mean and the principal components of spatially
  distributed functional data.
\newblock {\em Bernoulli}, 19(5A):1535--1558.

\bibitem[Jenish and Prucha, 2009]{jenish2009central}
Jenish, N. and Prucha, I.~R. (2009).
\newblock Central limit theorems and uniform laws of large numbers for arrays
  of random fields.
\newblock {\em Journal of Econometrics}, 150(1):86--98.

\bibitem[Kechris, 2012]{kechris2012classical}
Kechris, A. (2012).
\newblock {\em Classical Descriptive Set Theory}, volume 156 of {\em Graduate
  Texts in Mathematics}.
\newblock Springer.

\bibitem[Kelejian and Prucha, 2007]{kelejian2007hac}
Kelejian, H.~H. and Prucha, I.~R. (2007).
\newblock {HAC} estimation in a spatial framework.
\newblock {\em Journal of Econometrics}, 140(1):131--154.

\bibitem[Kim and Sun, 2011]{kim2011spatial}
Kim, M.~S. and Sun, Y. (2011).
\newblock Spatial heteroskedasticity and autocorrelation consistent estimation
  of covariance matrix.
\newblock {\em Journal of Econometrics}, 160(2):349--371.

\bibitem[Li et~al., 2008a]{li2008asymptotic}
Li, B., Genton, M.~G., and Sherman, M. (2008a).
\newblock On the asymptotic joint distribution of sample space-time covariance
  estimators.
\newblock {\em Bernoulli}, 14(1):228--248.

\bibitem[Li et~al., 2008b]{li2008testing}
Li, B., Genton, M.~G., and Sherman, M. (2008b).
\newblock Testing the covariance structure of multivariate random fields.
\newblock {\em Biometrika}, 95(4):813--829.

\bibitem[Markello and Misic, 2021]{markello2021comparing}
Markello, R.~D. and Misic, B. (2021).
\newblock Comparing spatial null models for brain maps.
\newblock {\em NeuroImage}, 236:118052.

\bibitem[Martens and Meesters, 1998]{martens1998image}
Martens, J.-B. and Meesters, L. (1998).
\newblock Image dissimilarity.
\newblock {\em Signal Processing}, 70(3):155--176.

\bibitem[Medrano et~al., 2024]{medrano2024assessing}
Medrano, J., Kheddar, A., and Ramdani, S. (2024).
\newblock Assessing time series correlation significance: A parametric approach
  with application to physiological signals.
\newblock {\em Biomedical Signal Processing and Control}, 94:106235.

\bibitem[Newey and West, 1994]{newey1994automatic}
Newey, W.~K. and West, K.~D. (1994).
\newblock Automatic lag selection in covariance matrix estimation.
\newblock {\em The Review of Economic Studies}, 61(4):631--653.

\bibitem[Politis and Romano, 1994]{politis1994large}
Politis, D.~N. and Romano, J.~P. (1994).
\newblock Large sample confidence regions based on subsamples under minimal
  assumptions.
\newblock {\em The Annals of Statistics}, 22(4):2031--2050.

\bibitem[Prause and Steland, 2018]{prause2018estimation}
Prause, A. and Steland, A. (2018).
\newblock Estimation of the asymptotic variance of univariate and multivariate
  random fields and statistical inference.
\newblock {\em Electronic Journal of Statistics}, 12(1):890--940.

\bibitem[Qin, 2026]{RFCorr}
Qin, J. (2026).
\newblock {RFCorr}: A {P}ython package for inference on the correlation between
  random fields.
\newblock \url{https://github.com/JiyueQin/RFcorr}.
\newblock Python package.

\bibitem[Qin et~al., 2025]{qin2025simultaneous}
Qin, J., Davenport, S., and Schwartzman, A. (2025).
\newblock Simultaneous confidence regions for image excursion sets: A
  validation study with applications in f{MRI}.
\newblock {\em Imaging Neuroscience}, 3:IMAG.a.1044.

\bibitem[Quenouille, 1947]{quenouille1947notes}
Quenouille, M.~H. (1947).
\newblock Notes on the calculation of autocorrelations of linear autoregressive
  schemes.
\newblock {\em Biometrika}, 34(3/4):365--367.

\bibitem[Roy, 1989]{roy1989asymptotic}
Roy, R. (1989).
\newblock Asymptotic covariance structure of serial correlations in
  multivariate time series.
\newblock {\em Biometrika}, 76(4):824--827.

\bibitem[Schmidt, 2014]{schmidt2014stochastic}
Schmidt, V., editor (2014).
\newblock {\em Stochastic Geometry, Spatial Statistics and Random Fields}.
\newblock Lecture Notes in Mathematics. Springer.

\bibitem[Shashkin, 2002]{shashkin2002quasi}
Shashkin, A.~P. (2002).
\newblock Quasi-associatedness of a {G}aussian system of random vectors.
\newblock {\em Russian Mathematical Surveys}, 57(6):1243--1244.

\bibitem[Sommerfeld et~al., 2018]{sommerfeld2018confidence}
Sommerfeld, M., Sain, S., and Schwartzman, A. (2018).
\newblock Confidence regions for spatial excursion sets from repeated random
  field observations, with an application to climate.
\newblock {\em Journal of the American Statistical Association},
  113(523):1327--1340.

\bibitem[Vallejos et~al., 2020]{vallejos2020spatial}
Vallejos, R., Osorio, F., and Bevilacqua, M. (2020).
\newblock {\em Spatial Relationships between Two Georeferenced Variables: With
  Applications in {R}}.
\newblock Springer.

\bibitem[Van~der Vaart, 2000]{van2000asymptotic}
Van~der Vaart, A.~W. (2000).
\newblock {\em Asymptotic Statistics}, volume~3 of {\em Cambridge Series in
  Statistical and Probabilistic Mathematics}.
\newblock Cambridge University Press.

\bibitem[Viladomat et~al., 2014]{viladomat2014assessing}
Viladomat, J., Mazumder, R., McInturff, A., McCauley, D.~J., and Hastie, T.
  (2014).
\newblock Assessing the significance of global and local correlations under
  spatial autocorrelation: a nonparametric approach.
\newblock {\em Biometrics}, 70(2):409--418.

\bibitem[Weir, 1996]{weir1996genetic}
Weir, B.~S. (1996).
\newblock {\em Genetic Data Analysis {II}}.
\newblock Sinauer Associates, Sunderland, MA.

\bibitem[Zhang, 2004]{zhang2004inconsistent}
Zhang, H. (2004).
\newblock Inconsistent estimation and asymptotically equal interpolations in
  model-based geostatistics.
\newblock {\em Journal of the American Statistical Association},
  99(465):250--261.

\end{thebibliography}

\clearpage
\appendix
\section{Appendix}
\etocsettocstyle{\subsection*{Contents}}{}
\localtableofcontents
\subsection{Results under the generalized \texorpdfstring{$(BL,\theta)$}{(BL,theta)} definition}\label{append:clt}
{This section justifies Remark~\ref{BL-general}, which allows $\theta_r=\infty$ for
finitely many $r$ in the definition of $(BL,\theta)$-dependence.
The generalization affects two results of \citet{bulinski2007limit} that the paper uses.
The first is the absolute summability of the covariances of a weakly stationary $(BL,\theta)$
field (Remark~\ref{remark:finite_sus}), which Lemmas~\ref{lemma:bl-cov} and \ref{lemma:abs_sum}
establish.
The second is the CLT restated in Section~\ref{sec:pre}, which Lemma~\ref{lemma:clt_general}
establishes.
Throughout this section, $(BL,\theta)$ refers to the generalized definition.

We first bound sums of covariances between finite sets that are far apart.}

\begin{lemma}\label{lemma:bl-cov}
If a {vector-valued} random field $\{\bm{V}_s \in \mathbb{R}^p:s\in \mathbb{Z}^d\}$ is
$(BL,\theta)$ {and square integrable}, then for
all finite disjoint sets $I, J \subset \mathbb{Z}^d$ with $\operatorname{dist}(I, J)=r \in \mathbb{N}$
{and all $l,m\in\{1,\dots,p\}$}, we have
\[
\sum_{i \in I} \sum_{j \in J} |\cov(V_i^{(l)}, V_j^{(m)})|
\le (|I|\wedge |J|)\theta_r,
\]
{where $V_s^{(l)}$ denotes the $l$th component of $\bm{V}_s$.}
\end{lemma}
\begin{proof}[Proof of Lemma \ref{lemma:bl-cov}]
{
The linear functions needed here are unbounded, so we first truncate.
For $M>0$, let $c_M(v)=\max\{-M,\min(M,v)\}$ be the clipping map at level $M$, which is bounded by $M$ and $1$-Lipschitz.
Fix $j \in J$, write $\epsilon_i = \operatorname{sgn}(\cov(V_i^{(l)},V_j^{(m)}))$ for $i\in I$, and
define
\[
f_j^M(\bm{v}_I) = \sum_{i \in I} \epsilon_i\, c_M(v_i^{(l)}), \qquad g^M(\bm{v}_j) =
c_M(v_j^{(m)}),
\]
where $v_i^{(l)}$ is the $l$th component of the block $\bm{v}_i\in\mathbb{R}^p$ of $\bm{v}_I$.
These functions read only component $l$ at the grid points of $I$ and component $m$ at $j$.
Then $f_j^M\in BL(p|I|)$ and $g^M\in BL(p)$, with $\operatorname{Lip}(f_j^M)\le 1$ since
$|f_j^M(\bm{v}_I)-f_j^M(\bm{w}_I)|\le\sum_{i\in I}|v_i^{(l)}-w_i^{(l)}| \le \|\bm{v}_I-\bm{w}_I\|_1$,
and similarly $\operatorname{Lip}(g^M)\le 1$.
Since $\operatorname{dist}(I,\{j\})\ge r$ and $\theta$ is nonincreasing, $(BL,\theta)$-dependence applied to the
sets $I$ and $\{j\}$ gives
\[
\big|\cov(f_j^M(\bm{V}_I), g^M(\bm{V}_j))\big| = \Big|\sum_{i\in I}\epsilon_i\cov(c_M(V_i^{(l)}),
c_M(V_j^{(m)}))\Big| \le \theta_r
\]
for every $M>0$.
As $M\to\infty$, $c_M(V_i^{(l)})\to V_i^{(l)}$ pointwise, and $|c_M(V_i^{(l)})c_M(V_j^{(m)})|\le
|V_i^{(l)}V_j^{(m)}|$ and $|c_M(V_i^{(l)})|\le|V_i^{(l)}|$, which are integrable by square
integrability and the Cauchy--Schwarz inequality.
By dominated convergence, $\cov(c_M(V_i^{(l)}), c_M(V_j^{(m)}))\to\cov(V_i^{(l)},V_j^{(m)})$ for
each $i\in I$, and since $I$ is finite,
\[
\sum_{i \in I} |\cov(V_i^{(l)}, V_j^{(m)})| = \lim_{M\to\infty}\sum_{i\in
I}\epsilon_i\cov(c_M(V_i^{(l)}), c_M(V_j^{(m)})) \le \theta_r.
\]
Summing over $j \in J$ gives
\[
\sum_{j \in J}\sum_{i \in I} |\cov(V_i^{(l)}, V_j^{(m)})|\leq |J|\theta_r.
\]
Exchanging the roles of $(I,l)$ and $(J,m)$, that is, fixing $i\in I$ and letting the test
function on $J$ read component $m$, gives the same sum bounded by $|I|\theta_r$.
The result follows.}
\end{proof}

{Lemma~\ref{lemma:bl-cov} controls covariances only at distances where $\theta_r$ is
finite.
Combining it with a direct bound at short distances gives absolute summability.}

\begin{lemma}\label{lemma:abs_sum}
{If a random field $V=\{\bm{V}_s\in\mathbb{R}^p, s\in\mathbb{Z}^d\}$ is weakly
stationary and $(BL,\theta)$, then $\sum_{k\in\mathbb{Z}^d}|\cov(V_0^{(l)}, V_k^{(m)})|<\infty$ for
all $l,m\in\{1,\dots,p\}$, where $V_s^{(l)}$ denotes the $l$th component of $\bm{V}_s$.
That is, $V$ has summable covariance, so Remark~\ref{remark:finite_sus} holds
under the generalized definition.}
\end{lemma}

\begin{proof}[Proof of Lemma \ref{lemma:abs_sum}]
{Fix $l,m\in\{1,\dots,p\}$.
Since $\theta_r\to 0$, there is an integer $r_0\geq 1$ with $\theta_{r_0}<\infty$.
Weak stationarity implies square integrability, so Lemma~\ref{lemma:bl-cov} applies to $V$.
We split the sum $\sum_{k\in\mathbb{Z}^d}|\cov(V_0^{(l)}, V_k^{(m)})|$ into far grid points with
$\|k\|_\infty\geq r_0$, which Lemma~\ref{lemma:bl-cov} controls, and near grid points with
$\|k\|_\infty<r_0$, where $\theta_r$ may be infinite and we bound each term directly.

For the far grid points, let $J$ be any finite set of grid points $j$ with $\|j\|_\infty\geq r_0$.
Then $0\notin J$ and $r = \operatorname{dist}(\{0\}, J)\geq r_0$.
Lemma~\ref{lemma:bl-cov} with $I=\{0\}$ gives $\sum_{j\in J}|\cov(V_0^{(l)}, V_j^{(m)})|\leq
\theta_r\leq \theta_{r_0}$, where the second inequality holds because $\theta$ is nonincreasing.
This bound is the same for every such $J$.
Lemma~\ref{lemma:bl-cov} applies only to finite sets, so we reach the infinite sum through the
increasing finite sets $J_N=\{j\in\mathbb{Z}^d: r_0\leq\|j\|_\infty\leq N\}$, $N\geq r_0$, whose
union is the set of all far grid points.
The partial sums $\sum_{j\in J_N}|\cov(V_0^{(l)}, V_j^{(m)})|$ are nondecreasing in $N$, because
every term is nonnegative, and they are bounded by $\theta_{r_0}$.
Hence they converge, and $\sum_{\|k\|_\infty\geq r_0}|\cov(V_0^{(l)},
V_k^{(m)})|\leq\theta_{r_0}$.

For the near grid points, $\|k\|_\infty<r_0$ means that every coordinate of $k$ is an integer in
$\{-(r_0-1),\dots,r_0-1\}$.
Each coordinate takes $2r_0-1$ values, so there are $(2r_0-1)^d$ near grid points, including $k=0$.
For each of them, the Cauchy--Schwarz inequality and weak stationarity give
$|\cov(V_0^{(l)}, V_k^{(m)})|\leq \{\operatorname{Var}(V_0^{(l)})\operatorname{Var}(V_k^{(m)})\}^{1/2}
= \{\operatorname{Var}(V_0^{(l)})\operatorname{Var}(V_0^{(m)})\}^{1/2}$.

Combining the two parts,
\[
\sum_{k\in\mathbb{Z}^d}|\cov(V_0^{(l)}, V_k^{(m)})|\leq
(2r_0-1)^d\{\operatorname{Var}(V_0^{(l)})\operatorname{Var}(V_0^{(m)})\}^{1/2}+\theta_{r_0}<\infty.
\]}
\end{proof}

{With absolute summability in hand, the CLT carries over to the generalized
definition.}

\begin{lemma}\label{lemma:clt_general}
{The CLT restated in Section~\ref{sec:pre} holds under the generalized definition
(Remark~\ref{BL-general}).}
\end{lemma}

\begin{proof}[Proof of Lemma \ref{lemma:clt_general}]
{We first treat the real-valued case $p=1$.
Let $X=\{X_s, s\in\mathbb{Z}^d\}$ be a real-valued, strictly stationary, centered and
square-integrable $(BL,\theta)$ field, and let $S(U)=\sum_{j\in U}X_j$ for finite
$U\subset\mathbb{Z}^d$.
It suffices to check the steps of the proof of \citet[Chapter 3, Theorem 1.12]{bulinski2007limit}
that involve the sequence $\theta=(\theta_r)_{r\in\mathbb{N}}$, and that proof uses $\theta$ in
three ways.

First, it bounds by a multiple of $\theta_q$ the difference between the characteristic function of a
sum over blocks at distance at least $q$ from each other and the product of the characteristic
functions of the block sums, where $q$ is taken large and then $q\to\infty$.
This involves $\theta_q$ only for large $q$, so it is unaffected when $\theta_r=\infty$ for finitely
many $r$.

Second, it uses $\theta_1$ through the variance bound $E\{S(U)^2\}\le\{E(X_0^2)+\theta_1\}|U|$, and
every use of $\theta_1$ applies this bound to a finite set $U$.
It therefore suffices that the bound holds with $\theta_1$ replaced by a finite constant.
Let $C_X=\sum_{s \in \mathbb{Z}^d}|\cov(X_0,X_s)|$, which is finite by Lemma~\ref{lemma:abs_sum}.
For any finite $U\subset\mathbb{Z}^d$,
\begin{align*}
E\{S(U)^2\} &=\sum_{j\in U}E(X_j^2) + \sum_{i\in U}\sum_{j \in U, j\neq
i}\cov(X_i,X_j)
\quad\text{(field is centered)}\\
&\leq \sum_{j\in U}E(X_j^2) + \sum_{i\in U}\sum_{j \in \mathbb{Z}^d, j\neq
i}|\cov(X_i,X_j)|\\
&=|U|E(X_0^2) + |U|\sum_{j \in \mathbb{Z}^d, j\neq 0}|\cov(X_0,X_j)|
\quad\text{(weak stationarity)}\\
&\leq\{E(X_0^2)+C_X\}|U|.
\end{align*}

Third, the proof uses the convergence of the series $\sum_{s\in\mathbb{Z}^d}\cov(X_0,X_s)$ that
defines the limit variance, through \citet[Chapter 3, Theorem 1.8]{bulinski2007limit}.
Under the generalized definition, this series converges absolutely by Lemma~\ref{lemma:abs_sum}.
Hence the proof of \citet[Chapter 3, Theorem 1.12]{bulinski2007limit} goes through, which gives
the case $p=1$.

For $p\geq 2$, the proof of \citet[Chapter 3, Corollary 1.13]{bulinski2007limit} reduces the
result to the case $p=1$ by the Cramér--Wold device, applied to the projected fields
$\{\bm{v}'\bm{X}_s, s\in\mathbb{Z}^d\}$ for $\bm{v}\in\mathbb{R}^p$.
Each projected field is strictly stationary, centered, square integrable and $(BL,\theta)$ by
Lemma~\ref{lemma:bl}, so the case $p=1$ applies to it.
Hence their proof goes through.}
\end{proof}

\subsection{Proofs for \texorpdfstring{Section~\ref{sec:pre}}{Section 2}}\label{append:proofs_pre}

\begin{proof}[Proof of Lemma \ref{lemma:bernstein_wlln}]
Write $\bar V_n = |U_n|^{-1}\sum_{s\in U_n} V_s$ and
$\psi_r=\sup\{|\cov(V_s, V_t)| : s,t\in\mathbb{Z}^d,\ \|s-t\|_\infty\geq r\}$ for $r\in\mathbb{N}$.
Let $\varepsilon>0$. Since $\psi_r\to0$, there is $r\in\mathbb{N}$ with $\psi_r<\varepsilon$, so that
$|\cov(V_s, V_t)|<\varepsilon$ whenever $\|s-t\|_\infty\geq r$.
Splitting the pairs of grid points at distance $r$,
$$
\operatorname{Var}\!\left(\sum_{s\in U_n} V_s\right) = \sum_{\substack{s,t\in U_n\\ \|s-t\|_\infty<r}}
\cov(V_s, V_t) + \sum_{\substack{s,t\in U_n\\ \|s-t\|_\infty\geq r}} \cov(V_s, V_t).
$$
By the Cauchy--Schwarz inequality, $|\cov(V_s, V_t)|\leq
\sqrt{\operatorname{Var}(V_s)\operatorname{Var}(V_t)}\leq c$ for all $s,t$.
For each $s\in U_n$ there are at most $(2r-1)^d$ grid points $t\in\mathbb{Z}^d$ with $\|s-t\|_\infty<r$, so
the first sum has at most $(2r-1)^d|U_n|$ terms and is bounded in absolute value by
$c\,(2r-1)^d|U_n|$.
The second sum has at most $|U_n|^2$ terms, each smaller than $\varepsilon$ in absolute value, so it is
bounded in absolute value by $\varepsilon|U_n|^2$.
Dividing by $|U_n|^2$ gives
$$
\operatorname{Var}(\bar V_n)\leq \frac{c\,(2r-1)^d}{|U_n|}+\varepsilon.
$$
Here $r$ depends on $\varepsilon$ but not on $n$, and $|U_n|\to\infty$, so
$\limsup_{n\to\infty}\operatorname{Var}(\bar V_n)\leq\varepsilon$.
Since $\varepsilon>0$ is arbitrary and variances are nonnegative, $\operatorname{Var}(\bar V_n)\to0$.
By Chebyshev's inequality, for every $\delta>0$,
$$
P\left(|\bar V_n-E(\bar V_n)|>\delta\right)\leq \frac{\operatorname{Var}(\bar V_n)}{\delta^2}\to0
\text{ as } n\to\infty,
$$
that is, $\bar V_n-E(\bar V_n)\xrightarrow{p}0$ as $n\to\infty$.
Since $\bar V_n-E(\bar V_n)=|U_n|^{-1}\sum_{s\in U_n}(V_s-E(V_s))$, this proves
$$
\frac{1}{|U_n|}\sum_{s\in U_n} \left(V_s-E(V_s)\right) \xrightarrow{p} 0 \text{ as } n\to \infty.
$$
If, in addition, $E(V_s)=\mu$ for all $s$, then $E(\bar V_n)=\mu$, and hence
$|U_n|^{-1}\sum_{s\in U_n}V_s=\bar V_n\xrightarrow{p}\mu$ as $n\to\infty$.
\end{proof}

\begin{proof}[Proof of Corollary \ref{cor:bernstein_wlln}]
By weak stationarity,
$$
\operatorname{Var}(V_s)=\operatorname{Var}(V_0),\qquad E(V_s)=E(V_0),\qquad
\cov(V_s, V_t)=\cov(V_0, V_{t-s})
$$
for all $s,t\in\mathbb{Z}^d$.
Every lag $k\in\mathbb{Z}^d$ arises from the pair $s=0$, $t=k$, so
$$
\sup\left\{|\cov(V_s, V_t)| : s,t\in\mathbb{Z}^d,\ \|s-t\|_\infty\geq r\right\}
=\sup_{\|k\|_\infty\geq r}|\cov(V_0, V_k)|,
$$
which tends to $0$ as $r\to\infty$ by vanishing covariance (Definition~\ref{def:dying_cov}).
The result follows from Lemma~\ref{lemma:bernstein_wlln} with $c=\operatorname{Var}(V_0)$ and
$\mu=E(V_0)$.
\end{proof}

\begin{proof}[Proof of Lemma \ref{lemma:finite_sus_rate}]
{Write $\bar V_n = |U_n|^{-1}\sum_{s\in U_n} V_s$ and $C_V = \sum_{k\in\mathbb{Z}^d}
|\cov(V_0, V_k)|$, which is finite by summable covariance.
For any finite $U\subset\mathbb{Z}^d$, weak stationarity gives $\cov(V_s, V_t) = \cov(V_0,
V_{t-s})$, so
$$
\operatorname{Var}\!\left(\sum_{s\in U} V_s\right) = \sum_{s,t\in U}\cov(V_s, V_t) \leq \sum_{s\in
U}\sum_{k\in\mathbb{Z}^d}|\cov(V_0, V_k)| = C_V|U|.
$$
Hence $\operatorname{Var}(\bar V_n) \leq C_V/|U_n|$, and Chebyshev's inequality gives, for every
$M>0$,
$$
P\!\left(\sqrt{|U_n|}\,|\bar V_n - E(V_0)| > M\right) \leq \frac{C_V}{M^2}.
$$
For any $\varepsilon>0$, choosing $M = \max\{1, \sqrt{C_V/\varepsilon}\}$ makes the right-hand side
at most $\varepsilon$, uniformly in $n$.
Thus $\sqrt{|U_n|}(\bar V_n - E(V_0)) = O_p(1)$, i.e., $\bar V_n = E(V_0) + O_p(1/\sqrt{|U_n|})$.
If $|U_n|\to\infty$, then $O_p(1/\sqrt{|U_n|}) = o_p(1)$, which is the WLLN.}
\end{proof}

\begin{proof}[Proof of Lemma \ref{lemma:bl}]
Let $l = l_1+\dots+l_q< \infty$.
Let $k$ be a positive integer.
Let
\[
\bm{\eta} = [\eta_1^{(1)}, \dots, \eta_1^{(p)},\dots, \eta_k^{(1)}, \dots, \eta_k^{(p)}]'\in\mathbb{R}^{kp},
\qquad
\bm{\eta}_m = [\eta_m^{(1)}, \dots, \eta_m^{(p)}]',\quad m = 1,\dots, k.
\]
For any $f\in BL(kq)$, define $\tilde{f}: \mathbb{R}^{kp}\rightarrow \mathbb{R}$ by
\[
\tilde{f}(\bm{\eta}) = f(h_1(\bm{\eta}_1), \dots, h_q(\bm{\eta}_1), \dots, h_1(\bm{\eta}_k), \dots,
h_q(\bm{\eta}_k)).
\]
Clearly, $\tilde{f}$ is bounded since $f$ is bounded.
Note that $\tilde{f}$ is $[\operatorname{Lip}(f)l]$-Lipschitz since for all $\tilde{\bm{\eta}}$, ${\bm{\eta}}\in
\mathbb{R}^{kp}$, $\tilde{\bm{\eta}}\ne{\bm{\eta}}$, we have:
{\small
\begin{align*}
&\quad\left|\tilde{f}(\tilde{\bm{\eta}}) -\tilde{f}(\bm{\eta})\right| \\
&\leq \operatorname{Lip}(f)\big\|[h_1(\tilde{\bm{\eta}}_1), \cdots, h_q(\tilde{\bm{\eta}}_1), \cdots,
h_1(\tilde{\bm{\eta}}_k), \cdots, h_q(\tilde{\bm{\eta}}_k)]'\\
&\qquad\qquad -[h_1(\bm{\eta}_1) ,\cdots,
h_q(\bm{\eta}_1) ,\cdots, h_1(\bm{\eta}_k) ,\cdots, h_q(\bm{\eta}_k)]'\big\|_1\\
& \leq \operatorname{Lip}(f)\left\{l_1\left\|\tilde{\bm{\eta}}_1-\bm{\eta}_1\right\|_1 + \cdots+
l_q\left\|\tilde{\bm{\eta}}_1-\bm{\eta}_1\right\|_1+\cdots+
l_1\left\|\tilde{\bm{\eta}}_k-\bm{\eta}_k\right\|_1+ \cdots+
l_q\left\|\tilde{\bm{\eta}}_k-\bm{\eta}_k\right\|_1\right\}\\
&=\operatorname{Lip}(f)\times l \times{\|\tilde{\bm{\eta}}-\bm{\eta}\|_1}.
\end{align*}}

Since $\{\bm{V}_s \in \mathbb{R}^p:s\in \mathbb{Z}^d\}$ is $(BL,\theta)$, there exists
$\left(\theta_r\right)_{r \in \mathbb{N}}$ that satisfies inequality \eqref{eq:bl}.
Let $\tilde{\theta}_r = l^2\theta_r$.
Then for all finite disjoint sets $I, J \subset \mathbb{Z}^d$ with $\operatorname{dist}(I, J)=r \in \mathbb{N}$,
and any functions $f \in BL(q|I|), g \in BL(q|J|)$, we have:
\begin{align*}
\left|\cov\left(f\left(h(\bm{V}_I)\right), g\left(h(\bm{V}_J)\right)\right)\right| &=
\left|\cov\left(\tilde{f}(\bm{V}_I), \tilde{g}\left(\bm{V}_J\right)\right)\right| \\
& \leq \operatorname{Lip}(\tilde{f}) \operatorname{Lip}(\tilde{g})(|I| \wedge|J|) \theta_r\\
& \leq \operatorname{Lip}(f) \operatorname{Lip}(g)l^2(|I| \wedge|J|) \theta_r\\
& = \operatorname{Lip}(f) \operatorname{Lip}(g)(|I| \wedge|J|) \tilde{\theta}_r.
\end{align*}
Therefore, $\{h(\bm{V}_s)\in \mathbb{R}^q:s\in \mathbb{Z}^d\}$ is $(BL,\theta)$.
\end{proof}

\begin{proof}[Proof of Lemma \ref{lemma:pa}] We prove the case for PA random fields.
The proof for NA random fields follows analogously.
Let $k$ be a positive integer.
Let
\[
\bm{\eta} = [\eta_1^{(1)}, \dots, \eta_1^{(p)},\dots, \eta_k^{(1)}, \dots, \eta_k^{(p)}]'\in\mathbb{R}^{kp},
\qquad
\bm{\eta}_m = [\eta_m^{(1)}, \dots, \eta_m^{(p)}]',\quad m = 1,\dots, k.
\]
For any $f\in M(kq)$, define $\tilde{f}: \mathbb{R}^{kp}\rightarrow \mathbb{R}$ by
\[
\tilde{f}(\bm{\eta}) = f(h_1(\bm{\eta}_1), \dots, h_q(\bm{\eta}_1), \dots, h_1(\bm{\eta}_k), \dots,
h_q(\bm{\eta}_k)).
\]
Clearly, $\tilde{f}$ is bounded since $f$ is bounded.
Note that $\tilde{f}$ is coordinate-wise nondecreasing Borel since each $h_i$ is coordinate-wise
nondecreasing Borel for $i = 1, \dots, q$.
Thus, $\tilde f \in M(kp)$.
Then for all finite disjoint sets $I, J \subset \mathbb{Z}^d$, and any functions $f \in M(q|I|), g
\in M(q|J|)$, we have:
$$
\cov\left(f\left(h(\bm{V}_I)\right), g\left(h(\bm{V}_J)\right)\right) =
\cov\left(\tilde{f}(\bm{V}_I), \tilde{g}\left(\bm{V}_J\right)\right) \geq 0.
$$
Therefore, $\left\{h(\bm{V}_s)\in \mathbb{R}^q:s\in \mathbb{Z}^d\right\}$ is PA.
\end{proof}

\begin{proof}[Proof of Lemma \ref{lemma:cox}]
Fix $l $ and $m$. Then we have:
\begin{align*}
    \sum_{j \in \mathbb{Z}^d} \left|\cov\left(V_0^{(l)},
V_j^{(m)}\right)\right| \leq C\sum_{j \in \mathbb{Z}^d}
\exp(-\lambda\|j\|_\infty).
\end{align*}
For each integer $k\geq 1$, let $E_k = \{j\in \mathbb{Z}^d: \|j\|_{\infty}=k\}.$ Then $|E_k| =
(2k+1)^d-(2k-1)^d= O(k^{d-1})$ and
\begin{align*}
\sum_{j \in \mathbb{Z}^d} \left|\cov\left(V_0^{(l)},V_j^{(m)}\right)\right|
&\leq C\Big(1+\sum_{k\geq 1} [(2k+1)^d-(2k-1)^{d}]\exp(-\lambda k)\Big).
\end{align*}
The series on the right converges, since its terms are $O(k^{d-1}e^{-\lambda k})$.
This holds for all $l$ and $m$, so $\bm{V}$ has summable covariance.

For the squared-exponential case, let $\lambda_{\min}>0$ be the smallest eigenvalue of $A$.
For all $j\in\mathbb{Z}^d$, we have $j'Aj \ge \lambda_{\min}\|j\|_2^2 \ge
\lambda_{\min}\|j\|_\infty^2 \ge \lambda_{\min}\|j\|_\infty$, where the last inequality holds because
$\|j\|_\infty$ is a nonnegative integer.
Hence $\exp(-(t-s)'A(t-s)) \le \exp(-\lambda_{\min}\|t-s\|_\infty)$, and \eqref{eq:exp_decay} holds
with $\lambda = \lambda_{\min}$.
\end{proof}

\subsection{Proofs for \texorpdfstring{Section~\ref{sec:theorem}}{Section 3}}\label{append:proofs_theorem}
\subsubsection{Asymptotic normality of the sample correlation}
\begin{proof}[Proof of Theorem \ref{theorem:zd}]
Recall $\tilde{\bm{T}}_s = [\tilde{X}_s^2, \tilde{Y}_s^2, \tilde{X}_s\tilde{Y}_s]'$ from Assumption
\ref{ass:bl}.
The standardized fields $\tilde X$ and $\tilde Y$ are well defined by Assumption
\ref{ass:zd_squareint}.
Note that $\tilde T=\{\tilde{\bm{T}}_s:s\in \mathbb{Z}^d\}$ is strictly stationary by Assumption
\ref{ass:zd_stationary}, since $\tilde{\bm{T}}_s$ is a function of $\bm{Z}_s$ and strict
stationarity is preserved.
Then $E[\tilde{\bm{T}}_s]$ is constant for all $s$ and let $\bm{\mu} =E[\tilde{\bm{T}}_s] = [1, 1,
\rho]'$.
Let $$\bm{w} = [w_1, w_2, w_3]' = \left[\frac{1}{|U_n|}\sum_{s \in U_n}\tilde{X}_s^2,
\frac{1}{|U_n|}\sum_{s \in U_n}\tilde{Y}_s^2, \frac{1}{|U_n|}\sum_{s \in
U_n}\tilde{X}_s\tilde{Y}_s\right]{'}.$$ Note that $\bm{w} =|U_n|^{-1}\sum_{s\in U_n}\tilde{\bm{T}}_s$.
Let $\bm{B}_s= \tilde{\bm{T}}_s-\bm{\mu}$.
Then the field $B=\{\bm{B}_s, s\in \mathbb{Z}^d\}$ is centered, strictly stationary, square
integrable and $(BL,\theta)$. It is square integrable because $E\tilde X_0^4$ and
$E\tilde Y_0^4$ are finite by Assumption \ref{ass:zd_squareint} and $\tilde X_0^2\tilde Y_0^2\le
(\tilde X_0^4+\tilde Y_0^4)/2$, and it is $(BL,\theta)$ by Assumption \ref{ass:bl}, since
subtracting a constant preserves this property.
By applying the CLT {restated in Section~\ref{sec:pre}, with $p=3$,} on
the field $B$ over the regularly growing domains $U_n$ (Assumption \ref{ass:regular}), we have
$$\sqrt{|U_n|}(\bm{w}-\bm{\mu})= |U_n|^{-1/2}\sum_{s\in U_n}\bm{B}_s \xrightarrow{d} \mathcal{N}(0,
\Omega),$$ where $\Omega = \sum_{k\in \mathbb{Z}^d}\cov(\bm{B}_0,\bm{B}_{k})
{=\sum_{k\in \mathbb{Z}^d}\cov(\tilde{\bm T}_0,\tilde{\bm T}_{k})}$.
{By the same theorem, every entry of {this series} converges
absolutely.}
Note that
\begin{align*}
\widehat{\cov}(X_0, Y_0) &\stackrel{\text{def}}{=} \frac{1}{|U_n|}\sum_{s \in
U_n}(X_s-\bar{X}_n)(Y_s-\bar{Y}_n) \\
&= \frac{1}{|U_n|}\sum_{s \in U_n}(X_s-\mu_x+\mu_x-\bar{X}_n)(Y_s-\mu_y+\mu_y-\bar{Y}_n)\\
&=\frac{1}{|U_n|}\sum_{s \in U_n}(X_s-\mu_x)(Y_s-\mu_y) + (\mu_x-\bar{X}_n)(\bar{Y}_n-\mu_y)\\
&=\frac{\sigma_x\sigma_y}{|U_n|}\sum_{s \in U_n}\tilde X_s\tilde Y_s +
(\mu_x-\bar{X}_n)(\bar{Y}_n-\mu_y)\\
&=\sigma_x\sigma_yw_3 - \sigma_x\sigma_y\frac{1}{|U_n|}\sum_{s \in U_n}\tilde
X_s\frac{1}{|U_n|}\sum_{s \in U_n}\tilde Y_s\\
&=\sigma_x\sigma_yw_3 + o_p(|U_n|^{-1/2}),
\end{align*}
where the last equal sign follows because by
Lemma \ref{lemma:finite_sus_rate}, which applies since $\tilde X$ and $\tilde Y$ are weakly
stationary by Assumptions \ref{ass:zd_stationary} and \ref{ass:zd_squareint}, each have summable covariance by Assumption \ref{ass:bl}, {and $|U_n|\to\infty$ by Assumption
\ref{ass:regular},} we have $|U_n|^{-1}\sum_{s \in U_n}\tilde
X_s=O_p(1/\sqrt{|U_n|})$, and $|U_n|^{-1}\sum_{s \in U_n}\tilde Y_s=O_p(1/\sqrt{|U_n|})$,
which implies $$\frac{1}{|U_n|}\sum_{s \in U_n}\tilde X_s\frac{1}{|U_n|}\sum_{s \in U_n}\tilde Y_s =
o_p(|U_n|^{-1/2}).$$
Similarly, we have $$\hat \sigma^2_x \equiv \widehat{\operatorname{Var}}(X_0) = \sigma^2_xw_1 +
o_p(|U_n|^{-1/2}) \quad[\text{by Lemma \ref{lemma:finite_sus_rate} on }\tilde X ]$$ and
$$\hat \sigma^2_y \equiv \widehat{\operatorname{Var}}(Y_0) = \sigma^2_yw_2 + o_p(|U_n|^{-1/2})
\quad[\text{by Lemma \ref{lemma:finite_sus_rate} on } \tilde Y].$$
Let $$\tilde{\bm{w}} = [\tilde{w}_1, \tilde{w}_2, \tilde{w}_3]' =\left[\frac{\hat
\sigma_x^2}{\sigma_x^2}, \frac{\hat \sigma_y^2}{\sigma_y^2}, \frac{\widehat{\cov}(X_0,
Y_0)}{\sigma_x \sigma_y}\right]',$$
then we have
\begin{equation}\label{eq:tilde_w}
\tilde{\bm{w}} = \bm{w}+o_p(|U_n|^{-1/2}).
\end{equation}
Since $\sqrt{|U_n|}(\bm{w}-\bm{\mu}) \xrightarrow{d} \mathcal{N}(0, \Omega)$, we have
\begin{align*}
\sqrt{|U_n|}(\tilde{\bm{w}}-\bm{\mu}) &= \sqrt{|U_n|}\{\bm{w}+o_p(|U_n|^{-1/2})-\bm{\mu}\}
\\
&=  \sqrt{|U_n|}(\bm{w}-\bm{\mu}) + o_p(1)\\
&\xrightarrow{d} \mathcal{N}(0, \Omega).
\end{align*}
{Define $f:\mathbb{R}^3 \rightarrow\mathbb{R}$ by $f(\bm{v}) = (v_1v_2)^{-1/2}v_3$ if
$v_1v_2>0$ and $f(\bm{v})=0$ otherwise.
Since $\tilde w_1\tilde w_2 = \hat\sigma_x^2\hat\sigma_y^2/(\sigma_x^2\sigma_y^2)$, we have $\hat
\rho = f(\tilde{\bm{w}})$ whenever both sample variances are positive.
When a sample variance is $0$, $\hat\rho$ is undefined, and we set $\hat\rho = 0 =
f(\tilde{\bm{w}})$.
The function $f$ is differentiable at $\bm{\mu} = [1, 1, \rho]'$, because $v_1, v_2>0$ in a
neighborhood of $\bm{\mu}$.}
By the delta method \cite[Theorem 3.1]{van2000asymptotic}, we have $$\sqrt{|U_n|}(\hat \rho-\rho) =
\sqrt{|U_n|}(f(\tilde{\bm{w}})-f(\bm{\mu}))\xrightarrow{d}\mathcal{N}(0, \sigma^2),$$ where
$\sigma^2 = \nabla f(\bm{\mu})'\Omega \nabla f(\bm{\mu})$ and $\nabla f(\bm{\mu})$ denotes the
gradient of $f$ evaluated at $\bm{\mu}$.
Note that $\sigma^2$ is nonnegative since $\Omega$ is a limiting variance-covariance matrix and is
positive semi-definite.
{Note that $\nabla f(\bm{\mu}) = [-0.5\rho, -0.5\rho, 1]'$, which we denote by
$\bm{\alpha}(\rho)$.
Then $\bm\alpha(\rho)'\tilde{\bm T}_s=L_s$ with $L_s$ in \eqref{eq:variance_v}, so
$$
\begin{aligned}
\sigma^2=\bm{\alpha}(\rho)'\Omega\,\bm{\alpha}(\rho)
&=\sum_{k\in\mathbb{Z}^d}\bm\alpha(\rho)'\cov(\tilde{\bm T}_0,\tilde{\bm T}_k)\,\bm\alpha(\rho)\\
&=\sum_{k\in\mathbb{Z}^d}\cov\big(\bm\alpha(\rho)'\tilde{\bm T}_0,\bm\alpha(\rho)'\tilde{\bm T}_k\big)
=\sum_{k\in\mathbb{Z}^d}\cov(L_0,L_k),
\end{aligned}
$$
where the series converges absolutely since every entry of $\Omega$ does.
This gives \eqref{eq:variance_v}.
Finally, $E[L_s]=\bm\alpha(\rho)'\bm\mu=-0.5\rho-0.5\rho+\rho=0$, and expanding
$\cov(L_0,L_k)=E[L_0L_k]$ gives $h_k$ in
\eqref{eq:variance}, so $\sigma^2=\sum_{k\in\mathbb{Z}^d}h_k$.}

\end{proof}

\begin{proof}[Proof of Corollary \ref{cor:zd_gaussian}]
If $Z$ is a Gaussian random field, then by Isserlis' Theorem, we have $E[\tilde{X}_0^2
\tilde{X}_{k}^2]=2\rho_{xx,k}^2+1$, $E[\tilde{Y}_0^2 \tilde{Y}_{k}^2]=2\rho_{yy,k}^2+1$,
$E[\tilde{X}_0^2 \tilde{Y}_{k}^2]=2\rho_{xy,k}^2+1$, $E[\tilde{X}_{k}^2
\tilde{Y}_{0}^2]=2\rho_{xy,-k}^2+1$, $E[\tilde{X}_0^2\tilde{X}_{k}\tilde{Y}_{k}] = \rho +
2\rho_{xx,k}\rho_{xy,k}$, $E[\tilde{X}_0\tilde{X}_{k}^2\tilde{Y}_0] = \rho +
2\rho_{xx,k}\rho_{xy,-k}$, $E[\tilde{Y}_0^2\tilde{Y}_{k}\tilde{X}_{k}] = \rho +
2\rho_{yy,k}\rho_{xy,-k}$, $E[\tilde{Y}_0\tilde{Y}_{k}^2\tilde{X}_0] = \rho +
2\rho_{yy,k}\rho_{xy,k}$, $E[\tilde{X}_0\tilde{X}_{k}\tilde{Y}_0\tilde{Y}_{k}] =
\rho^2+\rho_{xx,k}\rho_{yy,k}+\rho_{xy,k}\rho_{xy,-k}$.
{Substituting these moments into $h_k$ in \eqref{eq:variance} of
Theorem~\ref{theorem:zd} gives \eqref{eq:variance_gaussian}.}

\end{proof}

\begin{proof}[Proof of Theorem \ref{theorem:zd_rho0}]
Note $\tilde X^2$,
$\tilde Y^2$ and $\tilde X\tilde Y$ are all strictly stationary
by Assumption \ref{ass:zd_stationary}, since they are functions of $\bm{Z}_s$ and strict
stationarity is preserved.
Let $$\bm{w} = \left[w_1, w_2, w_3\right]' = \left[\frac{1}{|U_n|}\sum_{s \in U_n}\tilde{X}_s^2,
\frac{1}{|U_n|}\sum_{s \in U_n}\tilde{Y}_s^2, \frac{1}{|U_n|}\sum_{s \in
U_n}\tilde{X}_s\tilde{Y}_s\right]'.$$
The fields $\tilde X^2$ and $\tilde Y^2$ are weakly stationary by Assumptions
\ref{ass:zd_stationary} and \ref{ass:zd_squareint}, and each have vanishing covariance by Assumption
\ref{ass:finite}.
Since $|U_n|\to\infty$ (Assumption \ref{ass:regular}), Corollary~\ref{cor:bernstein_wlln} gives
$w_1 \xrightarrow{p} 1$ and $w_2\xrightarrow{p} 1$.
Note that $E[\tilde X_s\tilde Y_s]=0$ for all $s$ when $\rho=0$.
Since $\tilde X \tilde Y$ is $(BL,\theta)$ by Assumption \ref{ass:finite}, and $\tilde X \tilde
Y$ is square-integrable because $E\tilde X_0^2\tilde Y_0^2\le(E\tilde X_0^4+E\tilde
Y_0^4)/2<\infty$ by Assumption \ref{ass:zd_squareint}, by the CLT {restated in Section~\ref{sec:pre}, with $p=1$,} applied over the regularly growing
domains $U_n$ (Assumption \ref{ass:regular}), we have $$\sqrt{|U_n|}w_3= \frac{1}{\sqrt{|U_n|}}\sum_{s\in
U_n} \tilde{X}_s\tilde{Y}_s \xrightarrow{d} \mathcal{N}(0, \sigma^2),$$ where $\sigma^2 = \sum_{k\in
\mathbb{Z}^d}h_k$ and $h_k = \cov(\tilde{X}_0\tilde{Y}_0,
\tilde{X}_{k}\tilde{Y}_{k})=E[\tilde{X}_0\tilde{Y}_0\tilde{X}_{k}\tilde{Y}_{k}]$ as in
\eqref{eq:variance_0}.
{By the same theorem, this series converges absolutely.}
As in the proof of Theorem \ref{theorem:zd}, since $\tilde X$ and $\tilde Y$ each have summable covariance (Assumption \ref{ass:finite}), we have the following as in \eqref{eq:tilde_w}:
$$\tilde{\bm{w}} = \bm{w}+o_p(|U_n|^{-1/2}).$$ Note that this yields $\tilde{w}_3
=w_3+o_p(1/\sqrt{|U_n|})$.
In addition, it implies that $\tilde{w}_1 =w_1+o_p(1)=1+o_p(1)$ and $\tilde{w}_2
=w_2+o_p(1)=1+o_p(1)$.
{As in the proof of Theorem~\ref{theorem:zd}, $\hat\rho = f(\tilde{\bm{w}}) =
\tilde{w}_3\, g(\tilde{w}_1, \tilde{w}_2)$, where $g(v_1, v_2) = (v_1v_2)^{-1/2}$ if $v_1v_2>0$ and
$g(v_1, v_2)=0$ otherwise.
The function $g$ is continuous at $(1,1)$, so $g(\tilde{w}_1, \tilde{w}_2) = 1+o_p(1)$ by the
continuous mapping theorem.}
Since $\sqrt{|U_n|}w_3 \xrightarrow{d} \mathcal{N}(0, \sigma^2)$, we have
\begin{align*}
\sqrt{|U_n|} \hat \rho
&={\sqrt{|U_n|}\tilde{w}_3\, g(\tilde{w}_1, \tilde{w}_2)}\\
&={\{\sqrt{|U_n|}w_3+o_p(1)\}\{1+o_p(1)\}}
\xrightarrow{d}\mathcal{N}(0, \sigma^2).
\end{align*}
\end{proof}

\begin{proof}[Proof of Corollary \ref{cor:zd_rho0_gaussian}]
Since $Z$ is a Gaussian field, by Isserlis' Theorem, we have
$E[\tilde{X}_0\tilde{X}_{k}\tilde{Y}_0\tilde{Y}_{k}] =
\rho_{xx,k}\rho_{yy,k}+\rho_{xy,k}\rho_{xy,-k}$ when $\rho=0$ and we obtain the desired expression
as in \eqref{eq:variance_gaussian_0}.
\end{proof}

\begin{proof}[Proof of Corollary \ref{cor:zd_indep}]
Since $X$ and $Y$ are independent, $\tilde X$ and $\tilde Y$ are independent.
Also recall that  $\tilde X$ and $\tilde Y$ are standardized fields with mean 0 and SD 1.
Thus,  we can express \eqref{eq:variance_0} as
$$
E[\tilde{X}_0\tilde{X}_{k}\tilde{Y}_0\tilde{Y}_{k}] =
E[\tilde{X}_0\tilde{X}_{k}]E[\tilde{Y}_0\tilde{Y}_{k}]= \rho_{xx,k}\rho_{yy,k},$$
which yields the desired expression as in \eqref{eq:variance_indep}.
\end{proof}

\begingroup
\subsubsection{Asymptotic equivalence with Afyouni et al.\ on \texorpdfstring{$\mathbb{Z}$}{Z}}\label{append:afyouni}
{We show that the finite-sample approximation of \citet{afyouni2019effective} in
Remark~\ref{remark:1d} is asymptotically equivalent to the asymptotic variance
\eqref{eq:variance_gaussian}, in the sense that $n$ times their approximation converges to
\eqref{eq:variance_gaussian}.}

Let $d=1$ and $U_n=\{1,\dots,n\}$.
{By Corollary~\ref{cor:zd_gaussian}, $\sigma^2=\sum_{k\in\mathbb{Z}}h_k$, and we write
$h_k$ in \eqref{eq:variance_gaussian} as $h_k=A_k+B_k+C_k$}, where
\begin{align*}
A_k&=\frac{1}{2}\rho^2\left(\rho_{xx,k}^2+\rho_{yy,k}^2+\rho_{xy,k}^2+\rho_{xy,-k}^2\right),\qquad
B_k=-\rho\left(\rho_{xx,k}+\rho_{yy,k}\right)\left(\rho_{xy,k}+\rho_{xy,-k}\right),\\
C_k&=\rho_{xx,k}\rho_{yy,k}+\rho_{xy,k}\rho_{xy,-k}.
\end{align*}
At lag $0$, $\rho_{xx,0}=\rho_{yy,0}=1$ and $\rho_{xy,0}=\rho$, so
\[
h_0=\tfrac12\rho^2\left(2+2\rho^2\right)-4\rho^2+\left(1+\rho^2\right)=1-2\rho^2+\rho^4=(1-\rho^2)^2.
\]

From a first-order Taylor expansion of $\hat\rho$ for Gaussian series,
\citet[eq.~2]{afyouni2019effective} approximated $\operatorname{Var}(\hat\rho)$ by
\begin{align*}
V_n^{\mathrm{A}}=\frac{1}{n^2}\Big\{&(n-2)(1-\rho^2)^2
+\rho^2\sum_{k=1}^{M}w_k\left(\rho_{xx,k}^2+\rho_{yy,k}^2+\rho_{xy,k}^2+\rho_{xy,-k}^2\right)\\
&-2\rho\sum_{k=1}^{M}w_k\left(\rho_{xx,k}+\rho_{yy,k}\right)\left(\rho_{xy,k}+\rho_{xy,-k}\right)\\
&+2\sum_{k=1}^{M}w_k\left(\rho_{xx,k}\rho_{yy,k}+\rho_{xy,k}\rho_{xy,-k}\right)\Big\},
\end{align*}
where $w_k=n-2-k$ and $M$ is the largest lag used, with $M=n-2$ when no lags are truncated.

With $w_0=n-2$, so that $(n-2)(1-\rho^2)^2=w_0h_0$, and using $h_{-k}=h_k$,
\begin{align}
n\,V_n^{\mathrm{A}}
&=\frac{1}{n}\Big\{w_0h_0+\sum_{k=1}^{n-2}\left(2w_kA_k+2w_kB_k+2w_kC_k\right)\Big\}
=\frac{1}{n}\Big\{w_0h_0+\sum_{k=1}^{n-2}2w_kh_k\Big\}\nonumber\\
&=\frac{1}{n}\sum_{|k|\le n-2}w_{|k|}\,h_k
=\frac{1}{n}\sum_{|k|\le n-2}(n-2-|k|)\,h_k .\label{eq:afyouni_lag}
\end{align}

\begin{proposition}\label{prop:afyouni}
Under the assumptions of Corollary~\ref{cor:zd_gaussian} with $d=1$, $n\,V_n^{\mathrm{A}}\to\sigma^2$
as $n\to\infty$, where $\sigma^2$ is the asymptotic variance \eqref{eq:variance_gaussian}.
\end{proposition}
\begin{proof}[Proof of Proposition~\ref{prop:afyouni}]
Let $S_m=\sum_{|k|\le m}h_k$ be the partial sums of the series $\sigma^2$.
A partial sum $S_m$ contains $h_k$ if and only if $m\ge|k|$.
So among $S_0,\dots,S_{n-3}$, the term $h_k$ appears in the $n-2-|k|$ partial sums
$S_{|k|},\dots,S_{n-3}$, and
\[
\sum_{m=0}^{n-3}S_m=\sum_{|k|\le n-2}(n-2-|k|)\,h_k ,
\]
where the terms with $|k|=n-2$ have weight $0$.
Hence, by \eqref{eq:afyouni_lag},
\[
n\,V_n^{\mathrm{A}}=\frac{n-2}{n}\cdot\frac{1}{n-2}\sum_{m=0}^{n-3}S_m .
\]
Under the assumptions of Corollary~\ref{cor:zd_gaussian}, the series $\sigma^2$ converges, so
$S_m\to\sigma^2$.
By the Ces\`aro means theorem, $(n-2)^{-1}\sum_{m=0}^{n-3}S_m\to\sigma^2$, and since
$(n-2)/n\to1$, $n\,V_n^{\mathrm{A}}\to\sigma^2$.
\end{proof}
\endgroup

\begingroup
\subsubsection{Asymptotic equivalence with Clifford and Dutilleul on $\mathbb{Z}^d$}\label{append:dutilleul_zd}
\begin{proof}[Proof of Theorem~\ref{theorem:clifford}]
Write $N=|U_n|$ and let $\bm{1}_N$ be the vector of ones of length $N$, so that $I_N$ and $\bm{1}_N$ match $I_{|U_n|}$ and $\bm{1}_{|U_n|}$ in Section~\ref{sec:theorem_indep}.
For $k\in\mathbb{Z}^d$, let $C_{xx,k}=\cov(X_s,X_{s+k})$ and $C_{yy,k}=\cov(Y_s,Y_{s+k})$, so that
$\rho_{xx,k}=C_{xx,k}/C_{xx,0}$ and $\rho_{yy,k}=C_{yy,k}/C_{yy,0}$.
By weak stationarity, the entry of $\Sigma_X$ indexed by $(i,j)\in U_n\times U_n$
is $C_{xx,j-i}$, and similarly for $\Sigma_Y$.
Let $S_x = \sum_{k\in\mathbb{Z}^d}|C_{xx,k}|$ and $S_y = \sum_{k\in\mathbb{Z}^d}|C_{yy,k}|$, which
are finite since $X$ and $Y$ each have summable covariance.

We first consider the leading term $\operatorname{tr}(\Sigma_X\Sigma_Y)$.
By weak stationarity,
$$
C_{yy,-k} =\cov(Y_{s}, Y_{s-k}) = \cov(Y_{s-k}, Y_{s}) = C_{yy,k},
$$
so $C_{yy,i-j}=C_{yy,j-i}$ and
\begin{align*}
\operatorname{tr}(\Sigma_X\Sigma_Y) &= \sum_{i,j\in U_n} C_{xx,j-i}C_{yy,i-j} = \sum_{i,j\in U_n}
C_{xx,j-i}C_{yy,j-i}\\
&= C_{xx,0}C_{yy,0}\sum_{k\in\mathbb{Z}^d} |A(k)|\,\rho_{xx,k}\rho_{yy,k},
\end{align*}
where $|A(k)|$, with $A(k) = \{s\in U_n: s+k\in U_n\}$ as in Section~\ref{sec:var_special}, is the
number of pairs in $U_n$ at lag $k$.
We claim that $|A(k)|/N\to1$ for every fixed $k$.
Fix a lattice path $0=p_0, p_1, \dots, p_{\|k\|_1}=k$ with $\|p_q-p_{q-1}\|_1=1$, for example moving
one unit at a time along each coordinate.
The path depends only on $k$.
If $i\in U_n$ but $i+k\notin U_n$, the shifted path $i+p_0,\dots,i+p_{\|k\|_1}$ starts in $U_n$ and
ends outside it, so let $\tau(i)$ be the first step with $i+p_{\tau(i)}\notin U_n$.
Then $i+p_{\tau(i)-1}\in U_n$ and $\|p_{\tau(i)}-p_{\tau(i)-1}\|_\infty=1$, so $i+p_{\tau(i)}\in\delta U_n$.
For each $q=1,\dots,\|k\|_1$, let $B_q$ be the set of such $i$ with $\tau(i)=q$.
On $B_q$, the map $i\mapsto i+p_q$ is a translation by the fixed vector $p_q$, hence injective, so
$|B_q|\le|\delta U_n|$.
Summing over $q=1,\dots,\|k\|_1$ gives
$$
0\le N - |A(k)| \le \|k\|_1\,|\delta U_n|.
$$
Since $U_n$ is regularly growing, $|\delta U_n|/N\to 0$, and therefore $|A(k)|/N\to1$.
Moreover,
$$
|A(k)|\,|\rho_{xx,k}\rho_{yy,k}|/N\le|\rho_{xx,k}|,
$$
which is summable over $k$, so by dominated convergence,
$$
\frac{\operatorname{tr}(\Sigma_X\Sigma_Y)}{N C_{xx,0}C_{yy,0}} \to
\sum_{k\in\mathbb{Z}^d}\rho_{xx,k}\rho_{yy,k}.
$$
Since $\operatorname{tr}(\Sigma_X) = NC_{xx,0}$ and $\operatorname{tr}(\Sigma_Y) = NC_{yy,0}$, the
left-hand side equals $N V_n^{\mathrm{C}}$, which proves the limit for the formula of
\citet{clifford1989assessing}.

Next, we show that the centering matrix $\Pi$ changes each trace by only $O(1)$.
Using $\Pi = I_N - N^{-1}\bm{1}_N\bm{1}_N'$, the symmetry of $\Sigma_X$ and $\Sigma_Y$, and
$\operatorname{tr}(\bm{1}_N\bm{1}_N'A)=\operatorname{tr}(\bm{1}_N'A\bm{1}_N)=\bm{1}_N'A\bm{1}_N$ for any
$N\times N$ matrix $A$ (by the cyclic property of the trace, since $\bm{1}_N'A\bm{1}_N$ is a scalar),
we have
\begin{align*}
\operatorname{tr}(\Pi\Sigma_X) &= \operatorname{tr}(\Sigma_X) - N^{-1}\bm{1}_N'\Sigma_X\bm{1}_N,\\
\operatorname{tr}(\Pi\Sigma_X\Pi\Sigma_Y) &= \operatorname{tr}(\Sigma_X\Sigma_Y) -
2N^{-1}\bm{1}_N'\Sigma_X\Sigma_Y\bm{1}_N + N^{-2}(\bm{1}_N'\Sigma_X\bm{1}_N)(\bm{1}_N'\Sigma_Y\bm{1}_N).
\end{align*}
Here $\operatorname{tr}(\Sigma_X) = NC_{xx,0}$,
\[
|\bm{1}_N'\Sigma_X\bm{1}_N| \le \sum_{i,j\in U_n}|C_{xx,j-i}| \le NS_x,
\]
and similarly $|\bm{1}_N'\Sigma_Y\bm{1}_N|\le NS_y$.
In addition,
$$
|\bm{1}_N'\Sigma_X\Sigma_Y\bm{1}_N| \le \sum_{j\in U_n}\Big(\sum_{i\in U_n}|C_{xx,j-i}|\Big)\Big(\sum_{l\in
U_n}|C_{yy,l-j}|\Big) \le NS_xS_y.
$$
Hence $\operatorname{tr}(\Pi\Sigma_X) = NC_{xx,0} + O(1)$, $\operatorname{tr}(\Pi\Sigma_Y) = NC_{yy,0} +
O(1)$ and $\operatorname{tr}(\Pi\Sigma_X\Pi\Sigma_Y) = \operatorname{tr}(\Sigma_X\Sigma_Y) + O(1)$.

Combining these results,
\begin{align*}
N V_n^{\mathrm{D}} &= \frac{N\{\operatorname{tr}(\Sigma_X\Sigma_Y) + O(1)\}}{\{NC_{xx,0}+O(1)\}
\{NC_{yy,0}+O(1)\}}\\
&= \frac{\operatorname{tr}(\Sigma_X\Sigma_Y)/(NC_{xx,0}C_{yy,0}) +
O(N^{-1})}{1+O(N^{-1})} \to \sum_{k\in\mathbb{Z}^d}\rho_{xx,k}\rho_{yy,k},
\end{align*}
since $N\to\infty$.
Under the assumptions of Corollary~\ref{cor:zd_indep}, this limit is $\sigma^2$ in
\eqref{eq:variance_indep} by that corollary.
\end{proof}
\endgroup

\subsection{Proofs for \texorpdfstring{Section~\ref{sec:condition}}{Section 4}}\label{append:proofs_condition}

\begin{proof}[Proof of Lemma \ref{lemma:bounded_bl}]
Suppose $|X_s|\le M, |Y_s|\le M$ for all $s\in \mathbb{Z}^d$, where $M\in\mathbb{R}^+$.
{Lemma~\ref{lemma:bl} needs a function that is Lipschitz on all of $\mathbb{R}^2$, and
$x^2$, $y^2$ and $xy$ are not.
We therefore clip the arguments first.
Let $c_M(v)=\max\{-M,\min(M,v)\}$ be the clipping map from the proof of
Lemma~\ref{lemma:bl-cov}, which is $1$-Lipschitz and bounded by $M$, and define
$h:\mathbb{R}^2\rightarrow \mathbb{R}^3$ by
\[
h(x, y) = [h_1(x,y), h_2(x,y), h_3(x,y)]' = [c_M(x)^2, c_M(y)^2, c_M(x)c_M(y)]'.
\]
For all $(x,y), (\tilde x, \tilde y)\in\mathbb{R}^2$,
\begin{align*}
|h_1(\tilde x, \tilde y) - h_1(x, y)| &= |c_M(\tilde x)+c_M(x)|\,|c_M(\tilde x)-c_M(x)|\\
&\leq 2M|\tilde x - x| \leq 2M(|\tilde x - x| + |\tilde y - y|),
\end{align*}
so $h_1$ is $2M$-Lipschitz on $\mathbb{R}^2$, and similarly $h_2$ is $2M$-Lipschitz.
For $h_3$,
\begin{align*}
|h_3(\tilde x, \tilde y) - h_3(x, y)| &\leq |c_M(\tilde x)|\,|c_M(\tilde y)-c_M(y)| +
|c_M(\tilde x)-c_M(x)|\,|c_M(y)|\\
&\leq M(|\tilde x - x| + |\tilde y - y|),
\end{align*}
so $h_3$ is $M$-Lipschitz on $\mathbb{R}^2$.
Since $|X_s|\le M$ and $|Y_s|\le M$, we have $c_M(X_s)=X_s$ and $c_M(Y_s)=Y_s$, so $h(X_s,
Y_s)=[X_s^2, Y_s^2, X_sY_s]'$ for every $s\in\mathbb{Z}^d$.
The result follows by applying Lemma \ref{lemma:bl} with $h$.}
\end{proof}

\begin{proof}[Proof of Proposition \ref{prop:bounded}]
First, we show Assumption \ref{ass:bl} is met under Condition \ref{ass:B1} when $X$ and $Y$ are
bounded.
Since $X$ and $Y$ are bounded,  $\tilde T$ is square-integrable.
By Lemma \ref{lemma:bl}, the standardized field $\{[\tilde X_s, \tilde Y_s]'\}$ is bounded and $(BL,\theta)$.
Then by Lemma \ref{lemma:bounded_bl}, $\tilde T$ is $(BL,\theta)$.
Since $\tilde X$ and $\tilde Y$ are bounded and strictly stationary by Assumption
\ref{ass:zd_stationary}, they are weakly stationary, so by Remark \ref{remark:finite_sus}, $\tilde
X$ and $\tilde Y$ each have summable covariance.
And we have proved the result.
Second, Condition \ref{ass:B2} implies Condition \ref{ass:B1} by \citet[Chapter 1, Theorem
5.17]{bulinski2007limit}.
Finally, Condition \ref{ass:B3} implies  Condition \ref{ass:B2} by Lemma \ref{lemma:cox}.
\end{proof}

The proof of Lemma~\ref{lemma:gauss_square} uses the following Gaussian covariance bound. Write
$C^\infty(\mathbb{R}^n)$ for the class of functions $F:\mathbb{R}^n\to\mathbb{R}$ whose partial
derivatives of every order exist and are continuous, and $\|\xi\|_{L^2}=(E\xi^2)^{1/2}$ for a random
variable $\xi$. A function $h$ on $\mathbb{R}^n$ grows at most polynomially if
$|h(x)|\le c(1+\|x\|_2^k)$ for all $x$, for some constants $c$ and $k$.

\begin{lemma}\label{lemma:gauss_cov}
Let $\bm{U}\in\mathbb{R}^n$ and $\bm{V}\in\mathbb{R}^m$ be jointly Gaussian, with a possibly
singular joint covariance matrix. If $F\in C^\infty(\mathbb{R}^n)$ and $H\in C^\infty(\mathbb{R}^m)$,
and all their partial derivatives grow at most polynomially, then
\begin{equation}\label{eq:gauss_cov}
\big|\cov\big(F(\bm{U}),H(\bm{V})\big)\big|
\le\sum_{i=1}^n\sum_{j=1}^m\big|\cov(U_i,V_j)\big|\,
\|\partial_iF(\bm{U})\|_{L^2}\,\|\partial_jH(\bm{V})\|_{L^2}.
\end{equation}
\end{lemma}

\begin{proof}[Proof of Lemma~\ref{lemma:gauss_cov}]
By translation, assume $\bm{U}$ and $\bm{V}$ are centered. Let $\bm{W}=[\bm{U}',\bm{V}']'$ with
covariance matrix $\Sigma_W$, and set $\bar F(\bm{W})=F(\bm{U})$ and $\bar H(\bm{W})=H(\bm{V})$. Let
$\bm{W}^*$ be an independent copy of $\bm{W}$, and
$\bm{W}_t=t\bm{W}+\sqrt{1-t^2}\,\bm{W}^*$ for $t\in[0,1]$. The Gaussian covariance identity
\citep[Corollary~4.12]{bandeira2023matrix} gives
\[
\cov\big(F(\bm{U}),H(\bm{V})\big)
=\int_0^1E\big[\nabla\bar F(\bm{W})'\,\Sigma_W\,\nabla\bar H(\bm{W}_t)\big]\,dt.
\]
Since $\nabla\bar F(\bm{W})=[\nabla F(\bm{U})',\bm{0}']'$ and
$\nabla\bar H(\bm{W})=[\bm{0}',\nabla H(\bm{V})']'$, only the block of $\Sigma_W$ formed by the
covariances between $\bm{U}$ and $\bm{V}$ enters, and
\[
\nabla\bar F(\bm{W})'\,\Sigma_W\,\nabla\bar H(\bm{W}_t)
=\sum_{i=1}^n\sum_{j=1}^m\cov(U_i,V_j)\,\partial_iF(\bm{U})\,\partial_jH(\bm{V}_t),
\]
where $\bm{V}_t=t\bm{V}+\sqrt{1-t^2}\,\bm{V}^*$ is the $\bm{V}$ block of $\bm{W}_t$. Since $\bm{V}_t$ has the same distribution as
$\bm{V}$, the Cauchy--Schwarz inequality gives
$|E[\partial_iF(\bm{U})\,\partial_jH(\bm{V}_t)]|\le\|\partial_iF(\bm{U})\|_{L^2}
\|\partial_jH(\bm{V})\|_{L^2}$ for every $t$. Integrating over $t\in[0,1]$ yields
\eqref{eq:gauss_cov}.
\end{proof}

\begin{proof}[Proof of Lemma~\ref{lemma:gauss_square}]
Let $\bm{Z}_s=[X_s,Y_s]'$ and $\bm{T}_s=[X_s^2,Y_s^2,X_sY_s]'$, and write $Z_s^{(1)}=X_s$ and
$Z_s^{(2)}=Y_s$ for the components. Write
$\gamma^{ab}(k)=\cov(Z_0^{(a)},Z_k^{(b)})$ and set
\[
v_r=\sum_{\substack{k\in\mathbb{Z}^d\\ \|k\|_\infty\ge r}}\ \sum_{a,b=1}^{2}\big|\gamma^{ab}(k)\big|.
\]
Since $Z$ is weakly stationary and $(BL,\theta)$, it has summable covariance by
Lemma~\ref{lemma:abs_sum}, so $v_0<\infty$ and $v_r\downarrow0$.

Let $I,J\subset\mathbb{Z}^d$ be nonempty finite disjoint sets with $\operatorname{dist}(I,J)=r$, and take
$f\in BL(3|I|)$ and $g\in BL(3|J|)$. The goal is to show that
\begin{equation}\label{eq:gauss_goal}
\big|\cov\big(f(\bm{T}_I),g(\bm{T}_J)\big)\big|
\le K\operatorname{Lip}(f)\operatorname{Lip}(g)\,(|I|\wedge|J|)\,v_r
\end{equation}
for a constant $K$ that does not depend on $f$, $g$, $I$, $J$ or $r$.
Lemma~\ref{lemma:gauss_cov} requires $C^\infty$ functions, but $f$ and $g$ are only Lipschitz and
need not be differentiable. We therefore approximate them by smooth functions, using a convolution
approximation similar to the one in \citet{shashkin2002quasi}.
Let $f_\varepsilon=f*\psi_\varepsilon$ and $g_\varepsilon=g*\psi_\varepsilon$, where
$\psi_\varepsilon\ge0$ is a $C^\infty$ function with $\int\psi_\varepsilon=1$ that vanishes
outside the $\varepsilon$-ball $\{\bm{y}:\|\bm{y}\|_1\le\varepsilon\}$, taken on
$\mathbb{R}^{3|I|}$ for $f$ and on $\mathbb{R}^{3|J|}$ for $g$. By
\citet[Lemma]{shashkin2002quasi}, $f_\varepsilon$ and $g_\varepsilon$ are $C^\infty$ and bounded,
converge uniformly to $f$ and $g$ as $\varepsilon\to0$ \citep[see also][Theorem
8.14]{folland1999real}, and satisfy $\operatorname{Lip}(f_\varepsilon)\le\operatorname{Lip}(f)$
and $\operatorname{Lip}(g_\varepsilon)\le\operatorname{Lip}(g)$.
Their partial derivatives of every order are bounded. Indeed, since $f$ is bounded, hence
locally integrable, and $\psi_\varepsilon$ is smooth with compact support, differentiation under
the integral moves each derivative onto $\psi_\varepsilon$ \citep[Section 8.2, Exercise
7]{folland1999real}, for example $\partial_mf_\varepsilon=f*\partial_m\psi_\varepsilon$ and
$\partial_m\partial_{m'}f_\varepsilon=f*\partial_m\partial_{m'}\psi_\varepsilon$, where
$\partial_m$ is the partial derivative in the $m$th of the $3|I|$ coordinates. Hence, for
instance,
\begin{align*}
|\partial_mf_\varepsilon(\bm{x})|
&=\Big|\int f(\bm{y})\,\partial_m\psi_\varepsilon(\bm{x}-\bm{y})\,d\bm{y}\Big|
\le\|f\|_\infty\int|\partial_m\psi_\varepsilon(\bm{x}-\bm{y})|\,d\bm{y}\\
&=\|f\|_\infty\|\partial_m\psi_\varepsilon\|_{L^1}<\infty
\end{align*}
for every $\bm{x}$, where $\|f\|_\infty=\sup_{\bm{x}}|f(\bm{x})|$ and
$\|h\|_{L^1}=\int|h(\bm{y})|\,d\bm{y}$. The last equality shifts the variable of integration, and
the bound is finite because $f$ is bounded and $\partial_m\psi_\varepsilon$ is continuous and
vanishes outside the $\varepsilon$-ball. Higher derivatives are bounded in the same way, and
similarly for $g_\varepsilon$.

Suppose \eqref{eq:gauss_goal} holds for $f_\varepsilon$ and $g_\varepsilon$ with a constant $K$.
Since $\operatorname{Lip}(f_\varepsilon)\le\operatorname{Lip}(f)$ and
$\operatorname{Lip}(g_\varepsilon)\le\operatorname{Lip}(g)$, the right-hand side is at most
$K\operatorname{Lip}(f)\operatorname{Lip}(g)(|I|\wedge|J|)v_r$ for every $\varepsilon$. Since
$f_\varepsilon\to f$ and $g_\varepsilon\to g$ uniformly, and $f$, $g$, $f_\varepsilon$,
$g_\varepsilon$ are bounded, the left-hand side converges as $\varepsilon\to0$ to
$|\cov(f(\bm{T}_I),g(\bm{T}_J))|$. Indeed, adding and subtracting
$f(\bm{T}_I)g_\varepsilon(\bm{T}_J)$ in the product term, and
$E f(\bm{T}_I)\,E g_\varepsilon(\bm{T}_J)$ in the product of means, gives
\[
\begin{aligned}
&\cov\big(f_\varepsilon(\bm{T}_I),g_\varepsilon(\bm{T}_J)\big)
-\cov\big(f(\bm{T}_I),g(\bm{T}_J)\big)\\
&\quad=E\big[(f_\varepsilon-f)(\bm{T}_I)\,g_\varepsilon(\bm{T}_J)\big]
+E\big[f(\bm{T}_I)\,(g_\varepsilon-g)(\bm{T}_J)\big]\\
&\qquad-E\big[(f_\varepsilon-f)(\bm{T}_I)\big]\,E\big[g_\varepsilon(\bm{T}_J)\big]
-E\big[f(\bm{T}_I)\big]\,E\big[(g_\varepsilon-g)(\bm{T}_J)\big].
\end{aligned}
\]
Since $|f_\varepsilon-f|\le\|f_\varepsilon-f\|_\infty$, $|g_\varepsilon|\le\|g_\varepsilon\|_\infty$,
$|f|\le\|f\|_\infty$ and $|g_\varepsilon-g|\le\|g_\varepsilon-g\|_\infty$ at every point, the
first and third terms are each at most $\|f_\varepsilon-f\|_\infty\|g_\varepsilon\|_\infty$ in
absolute value, and the second and fourth each at most $\|f\|_\infty\|g_\varepsilon-g\|_\infty$, so
\[
\big|\cov\big(f_\varepsilon(\bm{T}_I),g_\varepsilon(\bm{T}_J)\big)
-\cov\big(f(\bm{T}_I),g(\bm{T}_J)\big)\big|
\le2\big(\|f_\varepsilon-f\|_\infty\|g_\varepsilon\|_\infty+\|f\|_\infty\|g_\varepsilon-g\|_\infty\big),
\]
which tends to $0$ because $\|f_\varepsilon-f\|_\infty\to0$ and $\|g_\varepsilon-g\|_\infty\to0$,
while the sup norms $\|f\|_\infty<\infty$ and $\|g_\varepsilon\|_\infty\le\|g\|_\infty<\infty$
stay bounded uniformly in $\varepsilon$. Hence \eqref{eq:gauss_goal} holds for $f$ and $g$ with the
same $K$. It therefore suffices to prove \eqref{eq:gauss_goal} for smooth $f,g$ with bounded
derivatives, and from now on we assume that $f$ and $g$ are such functions.
Define
\[
q(x,y)=[x^2,y^2,xy]',\qquad
F(\bm{z}_I)=f\big((q(\bm{z}_i))_{i\in I}\big),\qquad
H(\bm{z}_J)=g\big((q(\bm{z}_j))_{j\in J}\big),
\]
where $\bm{z}_I=(\bm{z}_i)_{i\in I}$ with $\bm{z}_i=[x_i,y_i]'\in\mathbb{R}^2$, and similarly for
$\bm{z}_J$.
We check that $F$ and $H$ satisfy the growth assumptions of Lemma~\ref{lemma:gauss_cov}. First,
$|F|\le\|f\|_\infty$. For the derivatives, consider first a single grid point, $|I|=1$, so that
$F(x,y)=f(x^2,y^2,xy)$ with $f=f(a,b,c)$. Write $f_a,f_{ac},\dots$ for the partial derivatives of
$f$, evaluated at $(a,b,c)=(x^2,y^2,xy)$. The chain rule gives
\[
\partial_xF=2x\,f_a+y\,f_c,\qquad \partial_yF=2y\,f_b+x\,f_c,
\]
and, differentiating once more,
\[
\begin{aligned}
\partial_x^2F&=4x^2f_{aa}+4xy\,f_{ac}+y^2f_{cc}+2f_a,\\
\partial_y^2F&=4y^2f_{bb}+4xy\,f_{bc}+x^2f_{cc}+2f_b,\\
\partial_x\partial_yF&=4xy\,f_{ab}+2x^2f_{ac}+2y^2f_{bc}+xy\,f_{cc}+f_c.
\end{aligned}
\]
By the reduction above, $f$ has bounded derivatives. Since Lipschitz constants are taken with
respect to the $\ell_1$ norm, its first derivatives are bounded by $\operatorname{Lip}(f)$, and its
second derivatives are bounded by some $M_2$. With $m=\max(|x|,|y|)$, this gives
$|\partial F|\le3\operatorname{Lip}(f)\,m$ and $|\partial^2F|\le9M_2m^2+2\operatorname{Lip}(f)$
for every first and second partial derivative.

The same holds for general $I$ and for derivatives of any order. By the chain rule, every partial
derivative of $F$ is a sum of terms. Each term is a derivative of $f$, which is bounded, times a
product of derivatives of $q$, which is a polynomial. So every partial derivative of $F$ grows at
most polynomially. Hence $F$ and $H$ are $C^\infty$ and satisfy the growth assumptions of
Lemma~\ref{lemma:gauss_cov}. From the first-derivative formulas above and
$|f_a|,|f_b|,|f_c|\le\operatorname{Lip}(f)$, we get
\[
|\partial_{x_i}F|\le\operatorname{Lip}(f)\big(2|x_i|+|y_i|\big),\qquad
|\partial_{y_i}F|\le\operatorname{Lip}(f)\big(|x_i|+2|y_i|\big).
\]
Substituting $\bm{z}_I=\bm{Z}_I$ gives
$|\partial_{x_i}F(\bm{Z}_I)|\le\operatorname{Lip}(f)(2|X_i|+|Y_i|)$. By weak stationarity,
$\|X_i\|_{L^2}=\|X_0\|_{L^2}$ and $\|Y_i\|_{L^2}=\|Y_0\|_{L^2}$, which are finite since $Z$ is
Gaussian. Let $c=2\|X_0\|_{L^2}+2\|Y_0\|_{L^2}$. The triangle inequality in $L^2$ yields
\[
\|\partial_{x_i}F(\bm{Z}_I)\|_{L^2}
\le\operatorname{Lip}(f)\big(2\|X_i\|_{L^2}+\|Y_i\|_{L^2}\big)
\le c\operatorname{Lip}(f).
\]
In the same way $\|\partial_{y_i}F(\bm{Z}_I)\|_{L^2}\le c\operatorname{Lip}(f)$ for $i\in I$,
and $\|\partial_{x_j}H(\bm{Z}_J)\|_{L^2},\|\partial_{y_j}H(\bm{Z}_J)\|_{L^2}
\le c\operatorname{Lip}(g)$ for $j\in J$.
Since $\bm{T}_s=q(\bm{Z}_s)$, we have $f(\bm{T}_I)=F(\bm{Z}_I)$ and
$g(\bm{T}_J)=H(\bm{Z}_J)$. Hence
\begin{equation}\label{eq:gauss_main}
\begin{aligned}
\big|\cov\big(f(\bm{T}_I),g(\bm{T}_J)\big)\big|
&=\big|\cov\big(F(\bm{Z}_I),H(\bm{Z}_J)\big)\big|\\
&\le c^2\operatorname{Lip}(f)\operatorname{Lip}(g)\sum_{i\in I}\sum_{j\in J}\sum_{a,b=1}^{2}
\big|\cov\big(Z_i^{(a)},Z_j^{(b)}\big)\big|\\
&=c^2\operatorname{Lip}(f)\operatorname{Lip}(g)\sum_{i\in I}\sum_{j\in J}\sum_{a,b=1}^{2}
\big|\gamma^{ab}(j-i)\big|\\
&\le c^2\operatorname{Lip}(f)\operatorname{Lip}(g)\,(|I|\wedge|J|)\,v_r.
\end{aligned}
\end{equation}
The first inequality applies Lemma~\ref{lemma:gauss_cov} with $\bm{U}=\bm{Z}_I$ and
$\bm{V}=\bm{Z}_J$, whose coordinates are the $Z_i^{(a)}$ and $Z_j^{(b)}$, and
bounds each $L^2$ norm of a first derivative of $F$ by $c\operatorname{Lip}(f)$ and of $H$ by
$c\operatorname{Lip}(g)$. For the last inequality, for each fixed $i$ the lags $j-i$ are distinct
and satisfy $\|j-i\|_\infty\ge \operatorname{dist}(I,J)=r$, since distances are measured in $\|\cdot\|_\infty$,
giving the bound $|I|v_r$ for the sum. Fixing $j$ instead gives $|J|v_r$. This proves
\eqref{eq:gauss_goal} with $K=c^2$ for smooth $f,g$ with bounded derivatives, and hence, by the
reduction above, for all bounded Lipschitz $f,g$.

Thus $\{\bm{T}_s,s\in\mathbb{Z}^d\}$ is $(BL,c^2v)$-dependent.
\end{proof}

\begin{proof}[Proof of Proposition \ref{prop:gaussian}]
By Lemmas~\ref{lemma:cox} and \ref{lemma:abs_sum} and Remark~\ref{remark:cox_bl}, Conditions
\ref{ass:N1} and \ref{ass:N2} are equivalent, and Condition~\ref{ass:N3} implies
Condition~\ref{ass:N2}. These results include the convention allowing finitely many initial
infinite dependence coefficients. It suffices to assume Condition~\ref{ass:N1}.

By Assumptions~\ref{ass:zd_stationary} and \ref{ass:zd_squareint}, $Z$ is weakly stationary. The
standardized field $\{[\tilde X_s,\tilde Y_s]',s\in\mathbb{Z}^d\}$ is an affine transformation of
$Z$, so it is again a weakly stationary Gaussian field, and it is $(BL,\theta)$ by
Lemma~\ref{lemma:bl}. Applying Lemma~\ref{lemma:gauss_square} to it shows that $\tilde T$ is
$(BL,\theta)$. Moreover, $\tilde X$ and $\tilde Y$ are weakly stationary and, by
Lemma~\ref{lemma:bl}, $(BL,\theta)$, so they each have summable covariance by
Lemma~\ref{lemma:abs_sum}. This is exactly Assumption~\ref{ass:bl}.
\end{proof}

\begin{proof}[Proof of Proposition \ref{prop:mdep}]
Since $Z$ is $m$-dependent, so is $\tilde T$, which is therefore $(BL,\theta)$ by Remark
\ref{BL-general}.
Moreover, $\tilde X$ and $\tilde Y$ are weakly stationary by Assumptions \ref{ass:zd_stationary}
and \ref{ass:zd_squareint} and $(BL,\theta)$ by Remark \ref{BL-general}, so they each have summable covariance by Remark \ref{remark:finite_sus}.
Gaussian fields have finite moments of all orders, so $m$-dependent Gaussian fields satisfy
Assumption \ref{ass:zd_squareint}.
\end{proof}

{The proof of Lemma~\ref{lemma:indep} conditions on one of two independent fields and relies on
the following standard fact.}

\begin{lemma}\label{lemma:conditional_exp}

Let $\bm{X} \in \mathbb{R}^a$ and $\bm{Y} \in \mathbb{R}^b$ be independent random vectors, and let
$h: \mathbb{R}^{a+b} \rightarrow \mathbb{R}$ be a bounded Borel function.
Then $g: \mathbb{R}^a \rightarrow \mathbb{R}$, $g(\bm{x}) = E[h(\bm{x}, \bm{Y})]$, is bounded and
Borel, and
$E[h(\bm{X}, \bm{Y}) \mid \bm{X}] = g(\bm{X})$ almost surely.
\end{lemma}

\begin{proof}[Proof of Lemma \ref{lemma:conditional_exp}]
{
Let $\mu$ and $\nu$ be the distributions of $\bm{X}$ and $\bm{Y}$.
By independence, the distribution of $(\bm{X}, \bm{Y})$ is $\mu \times \nu$.
Since $h$ is bounded, $g(\bm{x}) = \int_{\mathbb{R}^b} h(\bm{x}, \bm{y}) \, \nu(d\bm{y})$ is
bounded, and it is Borel by Fubini's Theorem.
Here and below, $\mathcal{B}(\mathbb{R}^d)$ denotes the Borel $\sigma$-field on $\mathbb{R}^d$,
and $\sigma(\cdots)$ denotes the $\sigma$-field generated by $(\cdots)$.
Thus, $g(\bm{X})$ is $\sigma(\bm{X})$-measurable.
Let $A \in \sigma(\bm{X})$.
Then $A = \{\omega: \bm{X}(\omega) \in B\}$ for some $B \in \mathcal{B}(\mathbb{R}^a)$.
Since $\mathbbm{1}_A = \mathbbm{1}_B(\bm{X})$ and $(\bm{X}, \bm{Y})$ has distribution $\mu \times
\nu$, Fubini's Theorem gives
\begin{align*}
E[h(\bm{X}, \bm{Y}) \mathbbm{1}_A]
&= E[h(\bm{X}, \bm{Y}) \mathbbm{1}_B(\bm{X})] \\
&= \int_{\mathbb{R}^a} \int_{\mathbb{R}^b} h(\bm{x}, \bm{y}) \mathbbm{1}_B(\bm{x}) \,
\nu(d\bm{y}) \, \mu(d\bm{x}) \\
&= \int_{\mathbb{R}^a} g(\bm{x}) \mathbbm{1}_B(\bm{x}) \, \mu(d\bm{x}) \\
&= E[g(\bm{X}) \mathbbm{1}_A].
\end{align*}
Thus, $E[h(\bm{X}, \bm{Y}) \mid \bm{X}] = g(\bm{X})$ almost surely.
}
\end{proof}

\begin{proof}[Proof of Lemma \ref{lemma:indep}]
We first prove the case for PA random fields.
Let $\pi_1, \pi_2: \mathbb{R}^2 \rightarrow \mathbb{R}$ be the projections $\pi_1(x, y) = x$ and
$\pi_2(x, y) = y$.
Both are Borel and coordinate-wise nondecreasing.
Hence, if $\{[X_s,Y_s]'\}$ is PA, then $\{X_s\} = \{\pi_1([X_s,Y_s]')\}$ and
$\{Y_s\} = \{\pi_2([X_s,Y_s]')\}$ are both PA by Lemma \ref{lemma:pa}.
This direction does not use independence.

Conversely, suppose $\{X_s\}$ and $\{Y_s\}$ are both PA.
Let $I, J \subset \mathbb{Z}^d$ be finite disjoint sets, and let $f \in M(2|I|)$ and
$g \in M(2|J|)$.
Reordering the arguments of $f$ so that those for $X_I$ come before those for $Y_I$ keeps $f$ in
$M(2|I|)$.
So we may write $f$ as $f(x_I, y_I)$ with $x_I, y_I \in \mathbb{R}^{|I|}$, where $x_I$ is the
argument for $X_I$ and $y_I$ is the argument for $Y_I$, and similarly write $g$ as $g(x_J, y_J)$.
Define $\psi: \mathbb{R}^{|I|} \times \mathbb{R}^{|J|} \rightarrow \mathbb{R}$,
$\tilde f: \mathbb{R}^{|I|} \rightarrow \mathbb{R}$, and $\tilde g: \mathbb{R}^{|J|} \rightarrow
\mathbb{R}$ by
$$
\psi(x_I, x_J) = \cov(f(x_I, Y_I), g(x_J, Y_J)), \quad
\tilde f(x_I) = E[f(x_I, Y_I)], \quad
\tilde g(x_J) = E[g(x_J, Y_J)].
$$
For fixed $x_I$ and $x_J$, the functions $y_I \mapsto f(x_I, y_I)$ and $y_J \mapsto g(x_J, y_J)$
belong to $M(|I|)$ and $M(|J|)$.
Since $\{Y_s\}$ is PA, $\psi(x_I, x_J) \geq 0$ for all $x_I, x_J$.
Moreover, $\tilde f$ is bounded because $f$ is bounded, Borel by Fubini's Theorem, and
coordinate-wise nondecreasing because $f$ is.
Thus, $\tilde f \in M(|I|)$, and likewise $\tilde g \in M(|J|)$.
Since $\{X_s\}$ is PA, $\cov(\tilde f(X_I), \tilde g(X_J)) \geq 0$.

Since $\{X_s\}$ and $\{Y_s\}$ are independent, $(X_I, X_J)$ is independent of $(Y_I, Y_J)$.
Applying Lemma \ref{lemma:conditional_exp} to the bounded functions $fg$, $f$, and $g$ gives,
almost surely,
\begin{gather*}
E[f(X_I, Y_I) \mid X_I, X_J] = \tilde f(X_I),\qquad
E[g(X_J, Y_J) \mid X_I, X_J] = \tilde g(X_J),\\
\cov(f(X_I, Y_I), g(X_J, Y_J) \mid X_I, X_J) = \psi(X_I, X_J).
\end{gather*}
By the law of total covariance,
$$
\cov(f(X_I, Y_I), g(X_J, Y_J)) = E[\psi(X_I, X_J)] + \cov(\tilde f(X_I), \tilde g(X_J)) \geq 0.
$$
Thus, $\{[X_s,Y_s]' \in \mathbb{R}^2, s \in \mathbb{Z}^d\}$ is PA.

The proof for NA random fields is identical with all inequalities reversed.

For $(BL,\theta)$ random fields, the projections $\pi_1$ and $\pi_2$ are $1$-Lipschitz, so the
forward direction follows from Lemma \ref{lemma:bl}.
Conversely, suppose $\{X_s\}$ and $\{Y_s\}$ are $(BL,\theta)$ with sequences $\theta^X$ and
$\theta^Y$.
Let $I, J \subset \mathbb{Z}^d$ be finite disjoint sets with $\operatorname{dist}(I, J) = r$, and let
$f \in BL(2|I|)$ and $g \in BL(2|J|)$.
Permuting the arguments preserves $BL(\cdot)$ and the Lipschitz constant with respect to the
$\ell_1$ norm, so we write $f(x_I, y_I)$ and $g(x_J, y_J)$ as above.
For fixed $x_I$, the function $y_I \mapsto f(x_I, y_I)$ is bounded and Lipschitz with constant at
most $\operatorname{Lip}(f)$, and similarly for $g$.
Hence $|\psi(x_I, x_J)| \leq \operatorname{Lip}(f) \operatorname{Lip}(g) (|I| \wedge |J|)
\theta^Y_r$ for all $x_I, x_J$.
Also, $|\tilde f(x_I) - \tilde f(x_I^*)| \leq E|f(x_I, Y_I) - f(x_I^*, Y_I)| \leq
\operatorname{Lip}(f) \|x_I - x_I^*\|_1$, so $\tilde f \in BL(|I|)$ with
$\operatorname{Lip}(\tilde f) \leq \operatorname{Lip}(f)$, and likewise for $\tilde g$.
Hence $|\cov(\tilde f(X_I), \tilde g(X_J))| \leq \operatorname{Lip}(f) \operatorname{Lip}(g)
(|I| \wedge |J|) \theta^X_r$.
The same decomposition by the law of total covariance then gives
$$
|\cov(f(X_I, Y_I), g(X_J, Y_J))| \leq \operatorname{Lip}(f) \operatorname{Lip}(g)
(|I| \wedge |J|) (\theta^X_r + \theta^Y_r).
$$
The sequence $\theta^X + \theta^Y$ is nonincreasing and tends to $0$.
Thus, $\{[X_s,Y_s]'\}$ is $(BL,\theta)$ with $\theta = \theta^X + \theta^Y$.

\end{proof}

\begin{proof}[Proof of Lemma \ref{lemma:product}]
Taking the pointwise maximum of the sequences of $X$ and $Y$, we may assume that both fields are
$(BL,\theta)$ with the same sequence $\theta$.
Let $I, J \subset \mathbb{Z}^d$ be finite disjoint sets with $\operatorname{dist}(I,J)=r$, and let
$f \in BL(|I|)$ and $g \in BL(|J|)$. Write
\[
Q_I = (X_u Y_u)_{u \in I}, \quad Q_J = (X_v Y_v)_{v \in J}.
\]
We decompose
\[
\cov(f(Q_I), g(Q_J))
=
A + B,
\]
where
\[
A = \cov\big( E[f(Q_I)\mid X], E[g(Q_J)\mid X] \big),
\]
\[
B = E\big[ \cov(f(Q_I), g(Q_J) \mid X) \big].
\]
Now we bound A and B separately.
We first consider A.
Define
\[
h_I(x_I) = E[ f(x_I \odot Y_I) ],
\quad
h_J(x_J) = E[ g(x_J \odot Y_J) ],
\]
where $\odot$ denotes coordinate-wise multiplication.
By independence of $X$ and $Y$, we have
\[h_I(x_I) = E\left[ f(x_I \odot Y_I) |X=x\right], \quad
h_J(x_J) = E\left[ g(x_J \odot Y_J) |X=x\right].\]
Since $f$ is measurable and bounded, by \citet[Example 4.1.6]{durrett2019probability}
we have \[h_I(X_I) = E\left[ f(X_I \odot Y_I) |X\right],\] and similarly \[h_J(X_J) =
E\left[ g(X_J \odot Y_J) |X\right].\] Then
\[
A = \cov(h_I(X_I), h_J(X_J)).
\]
For $x_I, x_I' \in \mathbb{R}^{|I|}$,
\begin{align*}
|h_I(x_I) - h_I(x_I')|
&\le
E\big| f(x_I \odot Y_I) - f(x_I' \odot Y_I) \big|\\
&\leq E\big|\operatorname{Lip}(f) \|x_I \odot Y_I - x_I' \odot Y_I \|_1\big|\\
&=\operatorname{Lip}(f)\, E\| (x_I - x_I') \odot Y_I \|_1\\
&\leq \operatorname{Lip}(f) \sum_{i \in I} E|Y_i|\, |x_i - x_i'|\\
&\le
m_Y\, \operatorname{Lip}(f)\, \|x_I - x_I'\|_1,
\end{align*}
where $m_Y = \sup_i E|Y_i| \le \sqrt{M_Y}$ by Cauchy--Schwarz.
Thus
\[
\operatorname{Lip}(h_I) \le m_Y \operatorname{Lip}(f),
\qquad
\operatorname{Lip}(h_J) \le m_Y \operatorname{Lip}(g).
\]
By $(BL,\theta)$-dependence of $X$,
\begin{equation}\label{eq:gaussian-bl-A}
|A|
\le
m_Y^2\, \operatorname{Lip}(f)\operatorname{Lip}(g) (|I|\wedge |J|)\theta_r
\le
M_Y\, \operatorname{Lip}(f)\operatorname{Lip}(g) (|I|\wedge |J|)\theta_r.
\end{equation}
Now we consider B.
Conditioning on $X=x$, define
\[
F_x(y_I) = f(x_I \odot y_I),
\qquad
G_x(y_J) = g(x_J \odot y_J).
\]
Without loss of generality we assume $Y$ is PA or NA.
Since $F_x$ and $G_x$ are Lipschitz functions, by \citet[Chapter 1, Theorem 5.3]{bulinski2007limit}, we have
\[
|\cov(F_x(Y_I), G_x(Y_J)\mid X=x)|
\le
\sum_{i \in I} \sum_{j \in J}
\operatorname{Lip}_i(F_x) \operatorname{Lip}_j(G_x)\, |\cov(Y_i,Y_j)|.
\]
For each $i\in I$, varying only the $i$-th coordinate, we have
\begin{align*}
|F_x(y_I) -F_x(y_I')| &= |f(x_I \odot y_I)-f(x_I \odot y_I')|\\
&\leq \operatorname{Lip}(f)\|x_I \odot y_I-x_I \odot y_I'\|_1\\
&= \operatorname{Lip}(f)|x_i y_i-x_i y_i'|\\
&\leq \operatorname{Lip}(f)|x_i||y_i- y_i'|.
\end{align*}
Thus,  $\operatorname{Lip}_i(F_x) \le |x_i| \operatorname{Lip}(f)$ and
$\operatorname{Lip}_j(G_x) \le |x_j| \operatorname{Lip}(g)$. Thus
\[
|\cov(F_x(Y_I), G_x(Y_J)\mid X=x)|
\le
\operatorname{Lip}(f)\operatorname{Lip}(g)
\sum_{i \in I} \sum_{j \in J}
|x_i||x_j|\, |\cov(Y_i,Y_j)|.
\]

Define $h: \mathbb{R}^{2(|I|+|J|)}\to \mathbb{R}$, $h(x_I, x_J, y_I,
y_J)=[F_x(y_I)-E\{F_x(Y_I)\}]G_x(y_J)$.
Then we have
\[
|E(h(X_I, X_J, Y_I,Y_J)|X=x)|
\le
\operatorname{Lip}(f)\operatorname{Lip}(g)
\sum_{i \in I} \sum_{j \in J}
|x_i||x_j|\, |\cov(Y_i,Y_j)|.
\]
Note that $h$ is measurable since $F_x$ and $G_x$ are continuous, and that $E(|h(X_I, X_J, Y_I,
Y_J)|)<\infty $ since $F_x$ and $G_x$ are bounded.
Thus, by \citet[Example 4.1.6]{durrett2019probability},
\begin{align*}
|\cov(f(Q_I), g(Q_J)\mid X)| &=|E(h(X_I, X_J, Y_I,Y_J)\mid X)| \\
&\le
\operatorname{Lip}(f)\operatorname{Lip}(g)
\sum_{i \in I} \sum_{j \in J}
|X_i||X_j|\, |\cov(Y_i,Y_j)|  \quad \text{a.s.}
\end{align*}
Taking expectations on both sides and applying Cauchy--Schwarz,
\[E|X_i X_j| \le (EX_i^2)^{1/2}(EX_j^2)^{1/2}
\le M_X, \]
so
\[|B|\le M_X\, \operatorname{Lip}(f)\operatorname{Lip}(g)
\sum_{i \in I} \sum_{j \in J} |\cov(Y_i,Y_j)|.
\]
By Lemma \ref{lemma:bl-cov}, which applies since $Y$ is $(BL,\theta)$ and square
integrable with $\sup_{s}E[Y_s^2]\le M_Y$,
\begin{equation}\label{eq:gaussian-bl-B}
|B| \le M_X\, \operatorname{Lip}(f)\operatorname{Lip}(g)
(|I|\wedge |J|)\theta_r.
\end{equation}
Combining \eqref{eq:gaussian-bl-A} and \eqref{eq:gaussian-bl-B}, we have
\[\big|\cov(f(Q_I), g(Q_J))\big|\le (M_X + M_Y)\operatorname{Lip}(f)\operatorname{Lip}(g)
(|I|\wedge |J|)\theta_r.
\]
Thus $\{Q_s\}$ is $(BL,\tilde{\theta})$-dependent with
\[
\tilde{\theta}_r = (M_X + M_Y)\theta_r.
\]
\end{proof}

\begin{proof}[Proof of Lemma \ref{lemma:square_dying}]
Let $c_M(v)=\max\{-M,\min(M,v)\}$ be the clipping map from the proof of
Lemma~\ref{lemma:bl-cov}, and $h_M(v)=c_M(v)^2=\min(v^2,M^2)$. For $a,b\in\mathbb{R}$,
$|h_M(a)-h_M(b)|=|c_M(a)+c_M(b)|\,|c_M(a)-c_M(b)|\le2M|a-b|$, so $h_M\in BL(1)$ with
$\operatorname{Lip}(h_M)\le2M$. Let $k\in\mathbb{Z}^d$ with $\|k\|_\infty=r$ and $\theta_r<\infty$.
Applying $(BL,\theta)$-dependence with $I=\{0\}$ and $J=\{k\}$ gives
\[
\big|\cov\big(h_M(V_0),h_M(V_k)\big)\big|\le4M^2\theta_r.
\]
Adding and subtracting $\cov(h_M(V_0),V_k^2)$ gives
\[
\begin{aligned}
\cov(V_0^2,V_k^2)-\cov\big(h_M(V_0),h_M(V_k)\big)
={}&\cov\big(V_0^2-h_M(V_0),V_k^2\big)\\
&+\cov\big(h_M(V_0),V_k^2-h_M(V_k)\big).
\end{aligned}
\]
For square-integrable $U$ and $W$, the Cauchy--Schwarz inequality gives
\[
|\cov(U,W)|\le\{\operatorname{Var}(U)\operatorname{Var}(W)\}^{1/2}.
\]
Since
$\operatorname{Var}(U)=E[U^2]-(EU)^2\le E[U^2]$, and similarly for $W$, it follows that
\[
|\cov(U,W)|\le\|U\|_{L^2}\|W\|_{L^2},\qquad\text{where }\|U\|_{L^2}=(EU^2)^{1/2}.
\]
Since $0\le V_0^2-h_M(V_0)\le V_0^2\,\mathbbm{1}\{|V_0|>M\}$, we have
\[
\|V_0^2-h_M(V_0)\|_{L^2}\le\tau(M),\qquad\text{where }\tau(M)=\big(E[V_0^4\mathbbm{1}\{|V_0|>M\}]\big)^{1/2}.
\]
By
strict stationarity, $V_k$ has the same distribution as $V_0$, so the same bound holds for
$\|V_k^2-h_M(V_k)\|_{L^2}$, and $\|V_k^2\|_{L^2}=(E[V_0^4])^{1/2}$. Also
$\|h_M(V_0)\|_{L^2}\le\|V_0^2\|_{L^2}=(E[V_0^4])^{1/2}$ since $0\le h_M(V_0)\le V_0^2$. Hence
\[
\big|\cov(V_0^2,V_k^2)\big|\le4M^2\theta_r+2\,(E[V_0^4])^{1/2}\,\tau(M).
\]
Since $\theta_r\to0$, for every $M>0$,
\[
\limsup_{\|k\|_\infty\to\infty}\big|\cov(V_0^2,V_k^2)\big|\le2\,(E[V_0^4])^{1/2}\,\tau(M).
\]
The left-hand side does not depend on $M$. Since $E[V_0^4]<\infty$, dominated convergence gives
$\tau(M)\to0$ as $M\to\infty$, so
$\limsup_{\|k\|_\infty\to\infty}|\cov(V_0^2,V_k^2)|=0$, and thus $\cov(V_0^2,V_k^2)\to0$ as
$\|k\|_\infty\to\infty$. Finally, strict stationarity of $V$ implies strict stationarity of
$V^2$, which has finite second moments since $E[V_0^4]<\infty$. Hence $V^2$ is weakly stationary,
so $\cov(V_s^2,V_{s+k}^2)=\cov(V_0^2,V_k^2)$ for all $s$, and $V^2$ has vanishing covariance.
\end{proof}

\begin{proof}[Proof of Proposition \ref{prop:indep}]
First, we show that Condition~\ref{ass:I1} implies Assumption~\ref{ass:finite}. By
Assumption~\ref{ass:zd_stationary}, $X$ and $Y$ are strictly stationary, and by
Assumption~\ref{ass:zd_squareint} they have finite fourth moments. The maps
$v\mapsto(v-\mu_x)/\sigma_x$ and $v\mapsto(v-\mu_y)/\sigma_y$ are Lipschitz and nondecreasing, so
by Lemma~\ref{lemma:bl}, $\tilde X$ and $\tilde Y$ are $(BL,\theta)$, and by
Lemma~\ref{lemma:pa}, either $\tilde X$ or $\tilde Y$ is PA or NA. Since $\tilde X$ and $\tilde Y$
are weakly stationary, they each have summable covariance by
Remark~\ref{remark:finite_sus}. They are independent with $E[\tilde X_s^2]=E[\tilde Y_s^2]=1$ for
all $s$, so the moment condition of Lemma~\ref{lemma:product} holds, and by that lemma
$\tilde X\tilde Y$ is $(BL,\theta)$. Finally, $\tilde X$ and $\tilde Y$ are strictly stationary,
$(BL,\theta)$ and have finite fourth moments, so by Lemma~\ref{lemma:square_dying}, $\tilde X^2$
and $\tilde Y^2$ each have vanishing covariance. This proves Assumption~\ref{ass:finite}.

Second, Condition~\ref{ass:I2} implies Condition~\ref{ass:I1} by \citet[Chapter 1, Theorem
5.17]{bulinski2007limit}, since PA and NA fields with finite second moments and summable covariance
are $(BL,\theta)$ (Remark~\ref{remark:cox_bl}). Finally,
Condition~\ref{ass:I3} implies Condition~\ref{ass:I2} by Lemma~\ref{lemma:cox}.
\end{proof}

\subsection{Asymptotic theory for time series under mixing conditions}\label{append:mixing}
For $d=1$, limit theorems for time series have historically been established under strong mixing
conditions.
The main text uses $(BL,\theta)$-dependence instead, because mixing coefficients are
hard to bound \citep{doukhan1999new} and few mixing random fields are known for $d>1$
\citep[Chapter 1, Section 5]{bulinski2007limit}.
In this section, we give versions of Theorems~\ref{theorem:zd} and \ref{theorem:zd_rho0} in which the
$(BL,\theta)$-dependence in Assumptions~\ref{ass:bl} and \ref{ass:finite} is replaced by a strong
mixing condition.
We use the notation of Section~\ref{sec:pre} with $d=1$ and $U_n=\{1,\dots,n\}$, so that
$|U_n|=n$.

\begin{definition}[Strong mixing coefficient]
For a strictly stationary vector-valued time series $\{\bm{V}_t\in\mathbb{R}^p:t\in \mathbb{Z}\}$, the
strong mixing coefficient is
$$\alpha(k)=\sup\left\{ |P(A \cap B)-P(A) P(B)|: A \in \mathcal{F}_{-\infty}^{-k}, B \in
\mathcal{F}_0^{\infty}\right\}, \quad 0 \leq k \in \mathbb{Z},$$
where $\mathcal{F}_{i}^{j}=\sigma(\bm{V}_t: i\leq t\leq j)$ for $-\infty\leq i\leq j\leq \infty$.
\end{definition}
Note that $0\leq\alpha(k)\leq1/4$ for all $k$, since $|P(A \cap B)-P(A)P(B)| =
|\cov(\mathbbm{1}_A, \mathbbm{1}_B)|\leq \{\operatorname{Var}(\mathbbm{1}_A)
\operatorname{Var}(\mathbbm{1}_B)\}^{1/2}\leq 1/4$.

{
\begin{definition}[Summable strong mixing]\label{def:mixing}
A strictly stationary vector-valued time series $\{\bm{V}_t\in\mathbb{R}^p:t\in \mathbb{Z}\}$ with
strong mixing coefficient $\alpha(k)$ is called summable strong mixing if for some $r>2$,
$E\|\bm{V}_0\|_2^r<\infty$ and $\sum_{k=0}^{\infty} \alpha(k)^{1-2 /r}<\infty$.
\end{definition}}

\begin{remark}\label{remark:m_dependent_mixing}
If $\{\bm{V}_t:t\in \mathbb{Z}\}$ is strictly stationary and $m$-dependent, then
$\mathcal{F}_{-\infty}^{-k}$ and $\mathcal{F}_0^{\infty}$ are independent for $k>m$, so $\alpha(k)=0$
for $k>m$.
Hence it is summable strong mixing as soon as $E\|\bm{V}_0\|_2^r<\infty$ for some $r>2$.
\end{remark}

The following lemmas relate summable strong mixing to the conditions used in the main text.
{
\begin{lemma}\label{lemma:cov_mixing}
If a real-valued time series $\{V_t:t\in \mathbb{Z}\}$ is summable strong mixing, then it has summable covariance.
\end{lemma}

\begin{proof}[Proof of Lemma \ref{lemma:cov_mixing}]
Let $r>2$ be as in Definition~\ref{def:mixing}.
Since $E|V_0|^r<\infty$, $V_0$ has a finite variance, so $\{V_t\}$ is weakly stationary.
For $k\geq1$, $V_{-k}$ is $\mathcal{F}_{-\infty}^{-k}$-measurable and $V_0$ is
$\mathcal{F}_0^{\infty}$-measurable, so by the covariance inequality \cite[Theorem
14.13]{hansen2022econometrics} and stationarity,
$|\cov(V_0, V_k)| = |\cov(V_{-k}, V_0)| \leq 8(E|V_0|^r)^{2/r}\alpha(k)^{1-2/r}$.
Since $\cov(V_0,V_{-k})=\cov(V_0,V_k)$, we have
$$\sum_{k\in\mathbb{Z}}|\cov(V_0, V_k)| \leq \operatorname{Var}(V_0) +
16(E|V_0|^r)^{2/r}\sum_{k=1}^{\infty}\alpha(k)^{1-2/r}<\infty.$$
\end{proof}}

{Strong mixing is preserved under measurable transformations, which lets us pass from $Z$ to
$\tilde T$, $\tilde X$ and $\tilde Y$.}

\begin{lemma}
\label{lemma:mixing}
For a vector-valued time series $\{\bm{V}_t \in \mathbb{R}^p:t\in \mathbb{Z}\}$, $p\in \mathbb{N}$,
let $\alpha(k)$ be its strong mixing coefficient. Let $f: \mathbb{R}^p \rightarrow \mathbb{R}^q$ be a measurable function,
 $q\in \mathbb{N}$, and let $\tilde{\alpha}(k)$ be the strong mixing coefficient of $\left\{f(\bm{V}_t)\in \mathbb{R}^q:t\in \mathbb{Z}\right\}$,
 then $\tilde{\alpha}(k) \leq \alpha(k)$ for all $0\leq k \in \mathbb{Z}.$
 If $f$ is bijective, then $\tilde{\alpha}(k) = \alpha(k)$.
\end{lemma}

\begin{proof}[Proof of Lemma \ref{lemma:mixing}]
    First note that for every $t\in \mathbb{Z}$, we have
    $\sigma(f(\bm{V}_t)) \subseteq \sigma(\bm{V}_t)$. This is because
    $\sigma(f(\bm{V}_t)) = \left\{\bm{V}_t^{-1}(f^{-1}(B)):B \in \mathcal{B}(\mathbb{R}^q)\right\}$
    and
    $\sigma(\bm{V}_t) = \left\{\bm{V}_t^{-1}(B):B \in \mathcal{B}(\mathbb{R}^p)\right\}$.
    For every $B \in \mathcal{B}(\mathbb{R}^q)$, since $f$ is measurable,
    we have $f^{-1}(B) \in \mathcal{B}(\mathbb{R}^p)$. Since $\bm{V}_t$ is
    a random vector (i.e., measurable), we have
    $\bm{V}_t^{-1}(f^{-1}(B)) \in \sigma(\bm{V}_t)$. Therefore,
    $\sigma(f(\bm{V}_t)) \subseteq \sigma(\bm{V}_t)$. If $f$ is bijective,
    then $f^{-1}:  \mathbb{R}^q \rightarrow \mathbb{R}^p$ is measurable
    \citep[see e.g.][Corollary 15.2]{kechris2012classical}. Applying the
    previous result with
    $\left\{f(\bm{V}_t)\in \mathbb{R}^q:t\in \mathbb{Z}\right\}$ and
    $f^{-1}$, we have
    $\sigma(f^{-1}(f(\bm{V}_t))) \subseteq \sigma(f(\bm{V}_t))$ and thus
    $\sigma(\bm{V}_t) \subseteq \sigma(f(\bm{V}_t))$. So in this case,
    $\sigma(\bm{V}_t) = \sigma(f(\bm{V}_t))$. Then the result follows by
    the definition of the strong mixing coefficient.
\end{proof}

{The following lemma lets the mixing condition be checked on $X$ and $Y$ separately when
they are independent (Remark~\ref{remark:1d_assumption}).}

\begin{lemma}\label{lemma:mixing_indep}
Let $\alpha_1(k), \alpha_2(k), \alpha_3(k)$ be the strong mixing
coefficients of $\{X_s\in \mathbb{R}, s\in \mathbb{Z}\}$,
$\{Y_s\in \mathbb{R}, s\in \mathbb{Z}\}$ and
$\{[X_s,Y_s]'\in \mathbb{R}^2, s\in \mathbb{Z}\}$, respectively. Suppose
$\{X_s\in \mathbb{R}, s\in \mathbb{Z}\}$ and
$\{Y_s\in \mathbb{R}, s\in \mathbb{Z}\}$ are independent. Then
$\sum_{k=0}^{\infty} \alpha_1(k)^{1-2 /r_1}<\infty$ and
$\sum_{k=0}^{\infty} \alpha_2(k)^{1-2 /r_2}<\infty$ for some
$r_1,r_2>2$ if and only if
$\sum_{k=0}^{\infty} \alpha_3(k)^{1-2 /r_3}<\infty$ for some $r_3>2$.
\end{lemma}
\begin{proof}[Proof of Lemma \ref{lemma:mixing_indep}]
If $\sum_{k=0}^{\infty} \alpha_3(k)^{1-2 /r_3}<\infty$ for some
$r_3>2$, then we can take $r_1=r_2=r_3$ and the result follows since
$\alpha_1\leq \alpha_3$, $\alpha_2\leq \alpha_3$ by
Lemma \ref{lemma:mixing}. For the other direction,  by
\citet[Theorem 5.1]{bradley2005basic}, we have
$\alpha_3(k)\leq \alpha_1(k)+\alpha_2(k).$ Take
$r_3 = \max(r_1,r_2).$ Then
\begin{align*}
\sum_{k=0}^{\infty} \alpha_3(k)^{1-2 /r_3}&\leq\sum_{k=0}^{\infty} (\alpha_1(k)+\alpha_2(k))^{1-2 /r_3}\\
&\leq\sum_{k=0}^{\infty} \alpha_1(k)^{1-2 /r_3}+\alpha_2(k)^{1-2 /r_3}\\
&\leq\sum_{k=0}^{\infty} \alpha_1(k)^{1-2 /r_1}+\alpha_2(k)^{1-2 /r_2}\\
&<\infty.
\end{align*}
\end{proof}

\begin{lemma}
\label{lemma:moment}
For a random vector $\bm{V}=[V_1, \dots, V_p]'\in \mathbb{R}^p$, if $E|V_i|^k < \infty$ for all $i
=1, \dots, p$, then $E\|\bm{V}\|_q^k< \infty$ for any $q\geq1$.
\end{lemma}

\begin{proof}[Proof of Lemma \ref{lemma:moment}]
Let $\tilde{V} = \max_{i}{|V_i|}$. We have
$$\|\bm{V}\|_q^k = \left(\sum_{i=1}^p |V_i|^q\right)^{k/q}\leq
(p\tilde{V}^{q})^{k/q}=p^{k/q}\tilde{V}^{k}\leq p^{k/q}\sum_{i=1}^p|V_i|^k.$$
Therefore, $E\|\bm{V}\|_q^k\leq p^{k/q}\sum_{i=1}^pE|V_i|^k <\infty$.
\end{proof}

We now state the assumptions, which parallel Assumptions~\ref{ass:bl} and \ref{ass:finite}.
{
\begin{assumption}[Weak dependence under mixing, general $\rho$]\label{ass:mixing}
$\tilde T$ is summable strong mixing, and $\tilde X$ and $\tilde Y$ each have summable covariance.
\end{assumption}

\begin{assumptionStar}[Weak dependence under mixing, $\rho=0$ case]\label{ass:mixing_rho0}
$\tilde X\tilde Y$ is summable strong mixing, $\tilde X$ and $\tilde Y$ each have
summable covariance, and $\tilde X^2$ and $\tilde Y^2$ each have vanishing covariance.
\end{assumptionStar}}
As for Assumptions~\ref{ass:bl} and \ref{ass:finite}, Assumption~\ref{ass:mixing} implies
Assumption~\ref{ass:mixing_rho0}: each coordinate of $\tilde T$ is summable strong mixing by
Lemma~\ref{lemma:mixing} applied with the coordinate projections, and hence $\tilde X^2$ and
$\tilde Y^2$ each have summable covariance by Lemma~\ref{lemma:cov_mixing},
which implies vanishing covariance.
Both assumptions keep summable covariance for $\tilde X$ and $\tilde Y$, because the proofs need
the rate of Lemma~\ref{lemma:finite_sus_rate} for their sample means, whereas consistency, and hence
vanishing covariance (Corollary~\ref{cor:bernstein_wlln}), suffices for $\tilde X^2$ and $\tilde Y^2$.

{
\begin{remark}\label{remark:1d_assumption}
Both assumptions hold if $Z$ is summable strong mixing with some $r>2$ such that $E|X_0|^{2r}<\infty$
and $E|Y_0|^{2r}<\infty$.
Indeed, $\tilde T$, $\tilde X$ and $\tilde Y$ are measurable functions of $Z$, so their strong mixing
coefficients are bounded by that of $Z$ (Lemma~\ref{lemma:mixing}).
Moreover, $E|\tilde X_0|^{2r}<\infty$ and $E|\tilde Y_0|^{2r}<\infty$, so $E|\tilde X_0^2|^r$,
$E|\tilde Y_0^2|^r$ and $E|\tilde X_0\tilde Y_0|^r\leq (E|\tilde X_0|^{2r}E|\tilde
Y_0|^{2r})^{1/2}$ are finite, and thus $E\|\tilde{\bm{T}}_0\|_2^r<\infty$ by
Lemma~\ref{lemma:moment}.
Similarly, $E|\tilde X_0|^r$ and $E|\tilde Y_0|^r$ are finite.
Hence $\tilde T$, $\tilde X$ and $\tilde Y$ are summable strong mixing, and $\tilde X$ and $\tilde Y$
each have summable covariance by Lemma~\ref{lemma:cov_mixing}.
If $X$ and $Y$ are independent, the condition on $Z$ can be checked on $X$ and $Y$ separately by
Lemma~\ref{lemma:mixing_indep}, and by Remark~\ref{remark:m_dependent_mixing}, it holds for
$m$-dependent time series with the above moments.
\end{remark}}

\begin{theorem}
\label{theorem:1d}
Let $d=1$ and $U_n=\{1,\dots,n\}$.
Under Assumptions \ref{ass:zd_stationary}, \ref{ass:zd_squareint} and \ref{ass:mixing}, we have
$\sqrt{n}(\hat \rho-\rho) \xrightarrow{d} \mathcal{N}(0, \sigma^2)$ as $n\to\infty$, where
$\sigma^2=\sum_{k\in\mathbb{Z}}h_k {= \sum_{k\in\mathbb{Z}}\cov(L_0,L_k)}$ with $h_k$
{and $L_s$} as in Theorem~\ref{theorem:zd}.
\end{theorem}

\begin{proof}[Proof of Theorem \ref{theorem:1d}]
We follow the proof of Theorem~\ref{theorem:zd} and only change the step that uses the CLT for
$(BL,\theta)$-dependent random fields.
Let $\bm{w}$, $\bm{\mu}=[1,1,\rho]'$ and $\bm{B}_s=\tilde{\bm{T}}_s-\bm{\mu}$ be as in that proof.
The time series $B=\{\bm{B}_s:s\in\mathbb{Z}\}$ is centered and strictly stationary.
Since $\bm{B}_s$ differs from $\tilde{\bm{T}}_s$ by a constant, both generate the same
$\sigma$-fields, so $B$ has the same strong mixing coefficient as $\tilde T$.
In addition, $E\|\bm{B}_0\|_2^r\leq 2^{r-1}(E\|\tilde{\bm{T}}_0\|_2^r+\|\bm{\mu}\|_2^r)<\infty$.
Hence $B$ is summable strong mixing by Assumption~\ref{ass:mixing}, and by the CLT for strong mixing
time series \cite[Theorem 14.15]{hansen2022econometrics},
$$\sqrt{n}(\bm{w}-\bm{\mu})= \frac{1}{\sqrt{n}}\sum_{s=1}^{n}\bm{B}_s \xrightarrow{d} \mathcal{N}(0,
\Omega),$$
where $\Omega = \sum_{k\in \mathbb{Z}}{\cov(\bm{B}_0, \bm{B}_{k})}$.
{Every entry of this series converges absolutely.
Indeed, for $k\neq0$ and components $l,m$, the covariance inequality \citep[Theorem
14.13]{hansen2022econometrics} and stationarity give, as in the proof of
Lemma~\ref{lemma:cov_mixing},
\[
|\cov(B_0^{(l)}, B_k^{(m)})|\leq 8(E|B_0^{(l)}|^r)^{1/r}(E|B_0^{(m)}|^r)^{1/r}\alpha(|k|)^{1-2/r},
\]
where $E|B_0^{(l)}|^r\le E\|\bm{B}_0\|_2^r<\infty$, and $\sum_{k=1}^\infty\alpha(k)^{1-2/r}<\infty$.}
The remaining steps of the proof of Theorem~\ref{theorem:zd} use only
Lemma~\ref{lemma:finite_sus_rate} on $\tilde X$ and $\tilde Y$, which applies since they each have summable covariance by Assumption~\ref{ass:mixing} {and $|U_n|=n\to\infty$,}
and the delta method.
They therefore carry over unchanged: \eqref{eq:tilde_w} holds, $\hat\rho=f(\tilde{\bm{w}})$ with
{$f$ as defined there}, and $\sqrt{n}(\hat\rho-\rho)\xrightarrow{d}\mathcal{N}(0,
\bm{\alpha}(\rho)'\Omega\,\bm{\alpha}(\rho))$, where $\bm{\alpha}(\rho)=\nabla f(\bm{\mu})$ and
$\bm{\alpha}(\rho)'\Omega\,\bm{\alpha}(\rho)=\sum_{k\in\mathbb{Z}}h_k$ with $h_k$ as in
\eqref{eq:variance}.
\end{proof}

Since the proof of Corollary~\ref{cor:zd_gaussian} only evaluates $h_k$, it also holds under the
assumptions of Theorem~\ref{theorem:1d}.
For $d=1$, $h_k=h_{-k}$ by stationarity and $h_0$ {in \eqref{eq:variance_gaussian}}
reduces to $(1-\rho^2)^2$ for Gaussian $Z$, so {$\sigma^2$ in
Corollary~\ref{cor:zd_gaussian}} can be written as
\begin{align*}
    \sigma^2&=
\left(1-\rho^2\right)^2+\rho^2 \sum_{k=1}^{\infty}\left(\rho_{x x, k}^2+\rho_{yy, k}^2+\rho_{xy,
k}^2+\rho_{xy,-k}^2\right)\\
&\quad - 2 \rho \sum_{k=1}^{\infty}\left(\rho_{xx, k}+\rho_{yy, k}\right)\left(\rho_{xy,
k}+\rho_{xy,-k}\right)
+ 2 \sum_{k=1}^{\infty}\left(\rho_{x x, k} \rho_{yy, k}+\rho_{x y, k} \rho_{x y,-k}\right).
\end{align*}

{Similarly, Theorem~\ref{theorem:zd_rho0} holds under the weaker
Assumption~\ref{ass:mixing_rho0}.}

\begin{theorem}
\label{theorem:1d_rho0}
Let $d=1$ and $U_n=\{1,\dots,n\}$.
Under Assumptions \ref{ass:zd_stationary}, \ref{ass:zd_squareint} and \ref{ass:mixing_rho0}, suppose
$\rho=0$, then we have $\sqrt{n}\hat \rho \xrightarrow{d} \mathcal{N}(0, \sigma^2)$ as
$n\to\infty$, where $\sigma^2=\sum_{k\in\mathbb{Z}}h_k$ with $h_k$ as in \eqref{eq:variance_0}.
\end{theorem}

\begin{proof}[Proof of Theorem \ref{theorem:1d_rho0}]
We follow the proof of Theorem~\ref{theorem:zd_rho0} and only change the CLT step.
Since $\rho=0$, $E[\tilde X_s\tilde Y_s]=0$, and $\tilde X\tilde Y$ is strictly stationary and
summable strong mixing by Assumption~\ref{ass:mixing_rho0}.
By the CLT for strong mixing time series \cite[Theorem 14.15]{hansen2022econometrics},
$$\sqrt{n}w_3= \frac{1}{\sqrt{n}}\sum_{s=1}^{n} \tilde{X}_s\tilde{Y}_s \xrightarrow{d}
\mathcal{N}(0, \sigma^2),$$
where $\sigma^2 = \sum_{k\in \mathbb{Z}}E[\tilde{X}_0\tilde{Y}_0\tilde{X}_{k}\tilde{Y}_{k}]$
{converges absolutely by Lemma~\ref{lemma:cov_mixing}, since $\tilde X\tilde Y$ is
centered}.
The remaining steps use only Lemma~\ref{lemma:finite_sus_rate} on $\tilde X$ and
$\tilde Y$, which each have summable covariance, and
Corollary~\ref{cor:bernstein_wlln} on $\tilde X^2$ and $\tilde Y^2$, which each have vanishing covariance, both
by Assumption~\ref{ass:mixing_rho0}, and therefore carry over unchanged.
\end{proof}

Likewise, Corollaries~\ref{cor:zd_rho0_gaussian} and \ref{cor:zd_indep} hold under the assumptions
of Theorem~\ref{theorem:1d_rho0}.
For $d=1$, {$\sigma^2$ in Corollary~\ref{cor:zd_rho0_gaussian}} can be written as $\sigma^2 = 1+ 2
\sum_{k=1}^{\infty}(\rho_{x x, k} \rho_{yy, k}+\rho_{x y, k} \rho_{x y,-k})$.

\subsection{Nonnegativity of the variance estimators}\label{append:nonneg}

Section~\ref{sec:var_special} states that normalizing the sample autocovariances by $|U_n|$
instead of $|A(k)|$ makes the estimators \rfgeneral{}, \rfzero{} and \rfindep{} nonnegative.
Lemma~\ref{lemma:nonneg} gives the identity behind this, and the paragraph after its proof applies
it to the three estimators.

\begin{lemma}\label{lemma:nonneg}
Let $Q=\{l\in\mathbb{Z}^d:\ 0\le l_j\le m_j,\ j=1,\dots,d\}$ be the box with side lengths
$m_1+1,\dots,m_d+1$ and corner at the origin, and let $B_b=\{b+l:\ l\in Q\}$ be $Q$ shifted by
$b\in\mathbb Z^d$. Let $v_s\in\mathbb{R}$ for $s\in U_n$, set $v_s=0$ for $s\notin U_n$, and let
$\gamma_v(k)=|U_n|^{-1}\sum_{s\in A(k)}v_sv_{s+k}$. Then
$$
\sum_{|k_j|\le m_j}w_{\bm m}(k)\,\gamma_v(k)
=\frac{1}{|U_n|\prod_{j=1}^d(m_j+1)}\sum_{b\in\mathbb{Z}^d}\Big(\sum_{t\in B_b}v_t\Big)^2\ \ge0.
$$
Moreover, let $x_s,y_s\in\mathbb R$ for $s\in U_n$, set $x_s=y_s=0$ for $s\notin U_n$, and let
$$
\gamma_x(k)=\frac{1}{|U_n|}\sum_{s\in A(k)}x_sx_{s+k},\qquad
\gamma_y(k)=\frac{1}{|U_n|}\sum_{s\in A(k)}y_sy_{s+k}.
$$
Then $\sum_{|k_j|\le m_j}w_{\bm m}(k)\,\gamma_x(k)\gamma_y(k)\ge0$.
\end{lemma}

\begin{proof}[Proof of Lemma~\ref{lemma:nonneg}]
Note that $\gamma_v(k)=|U_n|^{-1}\sum_{s\in\mathbb Z^d}v_sv_{s+k}$, since
$v_sv_{s+k}=0$ unless $s\in A(k)$, and similarly for $\gamma_x$ and $\gamma_y$.

For the first claim, we first write the Bartlett weight as a count of pairs in the box $Q$.
In coordinate $j$, the number of pairs $(l_j,l_j')$ in $\{0,\dots,m_j\}$ with $l_j'-l_j=k_j$ is
$m_j+1-|k_j|$ if $|k_j|\le m_j$ and $0$ otherwise. The coordinates are independent, so
$$
\sum_{l\in Q}\sum_{l'\in Q}\mathbbm{1}\{l'-l=k\}
=\prod_{j=1}^d(m_j+1-|k_j|)=\prod_{j=1}^d(m_j+1)\cdot w_{\bm m}(k)
$$
if $|k_j|\le m_j$ for all $j$, and the count is $0$ otherwise.
Let $C=|U_n|\prod_{j=1}^d(m_j+1)$. Then
\begin{align*}
\sum_{|k_j|\le m_j}w_{\bm m}(k)\,\gamma_v(k)
&=\frac1C\sum_{|k_j|\le m_j}\ \sum_{l\in Q}\sum_{l'\in Q}\mathbbm{1}\{l'-l=k\}
\sum_{s\in\mathbb Z^d}v_sv_{s+k}\\
&=\frac1C\sum_{k\in\mathbb Z^d}\ \sum_{l\in Q}\sum_{l'\in Q}\mathbbm{1}\{l'-l=k\}
\sum_{s\in\mathbb Z^d}v_sv_{s+k}\\
&=\frac1C\sum_{l\in Q}\sum_{l'\in Q}\ \sum_{k\in\mathbb Z^d}\mathbbm{1}\{l'-l=k\}
\sum_{s\in\mathbb Z^d}v_sv_{s+k}\\
&=\frac1C\sum_{l\in Q}\sum_{l'\in Q}\ \sum_{s\in\mathbb Z^d}v_sv_{s+l'-l}\\
&=\frac1C\sum_{l\in Q}\sum_{l'\in Q}\ \sum_{b\in\mathbb Z^d}v_{b+l}v_{b+l'}\\
&=\frac1C\sum_{b\in\mathbb Z^d}\Big(\sum_{l\in Q}v_{b+l}\Big)\Big(\sum_{l'\in Q}v_{b+l'}\Big)\\
&=\frac1C\sum_{b\in\mathbb Z^d}\Big(\sum_{t\in B_b}v_t\Big)^2\ \ge0.
\end{align*}

For the second claim, write the product of the two autocovariances as a double sum and substitute
$s'=s+h$:
$$
\gamma_x(k)\gamma_y(k)
=\frac{1}{|U_n|^2}\sum_{s\in\mathbb Z^d}\sum_{s'\in\mathbb Z^d}x_sx_{s+k}\,y_{s'}y_{s'+k}
=\frac{1}{|U_n|^2}\sum_{h\in\mathbb Z^d}\sum_{s\in\mathbb Z^d}\big(x_sy_{s+h}\big)\big(x_{s+k}y_{s+k+h}\big).
$$
For each shift $h$, let $z^{(h)}_s=x_sy_{s+h}$. This field vanishes outside $U_n$, and the inner sum
is $\sum_{s\in\mathbb Z^d}z^{(h)}_sz^{(h)}_{s+k}=|U_n|\,\gamma_{z^{(h)}}(k)$. Hence
$\gamma_x(k)\gamma_y(k)=|U_n|^{-1}\sum_{h\in\mathbb Z^d}\gamma_{z^{(h)}}(k)$, and
$$
\sum_{|k_j|\le m_j}w_{\bm m}(k)\,\gamma_x(k)\gamma_y(k)
=\frac{1}{|U_n|}\sum_{h\in\mathbb Z^d}\ \sum_{|k_j|\le m_j}w_{\bm m}(k)\,\gamma_{z^{(h)}}(k)\ \ge0,
$$
since each inner sum is nonnegative by the first claim applied to $v=z^{(h)}$.
\end{proof}

{Fix a realization of the data, so that $\hat L_s$, $\hat X_s$ and $\hat Y_s$ are real
numbers, set to $0$ outside $U_n$ as in Section~\ref{sec:var_impl}.
With the normalization $|U_n|$, \rfgeneral{} and \rfzero{} become
$$
\sum_{|k_j|\le m_j}w_{\bm m}(k)\,\gamma_v(k)
$$
with $v_s=\hat L_s$ and $v_s=\hat X_s\hat Y_s$, respectively, and \rfindep{} becomes
$$
\sum_{|k_j|\le m_j}w_{\bm m}(k)\,\gamma_x(k)\gamma_y(k)
$$
with $x_s=\hat X_s$ and $y_s=\hat Y_s$.
These are nonnegative by the first and second claims of Lemma~\ref{lemma:nonneg}.
Since the realization is arbitrary, the three estimators are nonnegative for every sample.}

\subsection{Additional simulation results}\label{append:simu}
{
Figures~\ref{fig:coverage_rho_chisq} and \ref{fig:coverage_rho_gaussian} show the coverage rates
of the confidence intervals of Section~\ref{sec:simu} for all values of $\rho$ in the baseline
design, for fields with FWHM $2$.
Figure~\ref{fig:type1_symmetric} shows the rates of rejecting $H_0:\rho=0$ under the shift
design of Section~\ref{sec:simu_design}.
}

\setcounter{figure}{0}
\renewcommand{\thefigure}{S\arabic{figure}}
\renewcommand{\theHfigure}{S\arabic{figure}}

\begin{figure}[htbp]
    \centering
    \includegraphics[width=0.9\textwidth]{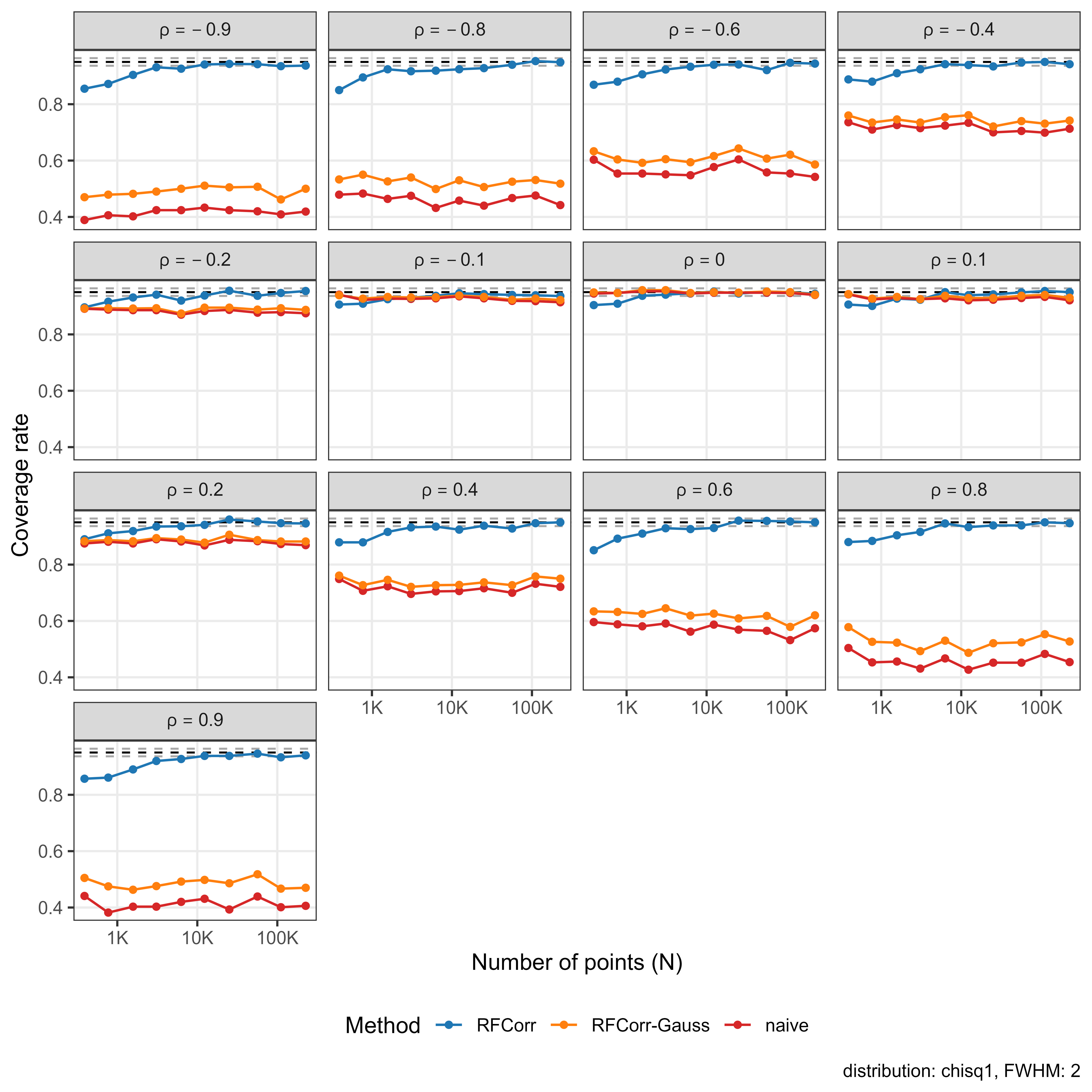}
    \caption{Coverage rate of $95\%$ confidence intervals against the number of points
    $N$ (log scale) under the baseline design, for $\chi^2_1$ fields with FWHM $2$ in $d=2$.
    Panels give the true correlation $\rho$.
    Dashed lines are as in Figure~\ref{fig:coverage}.}
    \label{fig:coverage_rho_chisq}
\end{figure}

\begin{figure}[htbp]
    \centering
    \includegraphics[width=0.9\textwidth]{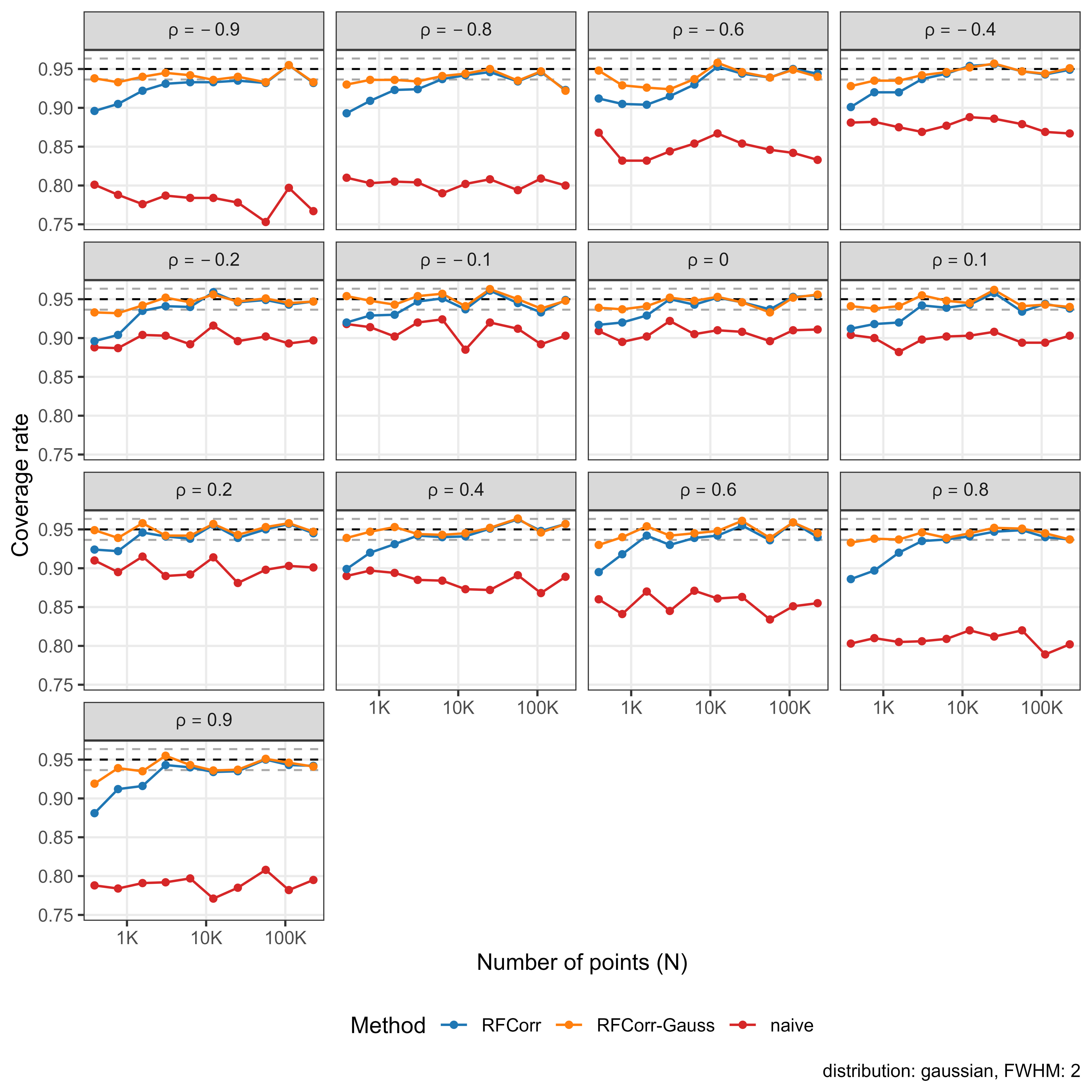}
    \caption{Coverage rate of $95\%$ confidence intervals against the number of points
    $N$ (log scale) under the baseline design, for Gaussian fields with FWHM $2$ in $d=2$.
    Panels give the true correlation $\rho$.
    Dashed lines are as in Figure~\ref{fig:coverage}.}
    \label{fig:coverage_rho_gaussian}
\end{figure}

\begin{figure}[htbp]
    \centering
    \includegraphics[width=0.9\textwidth]{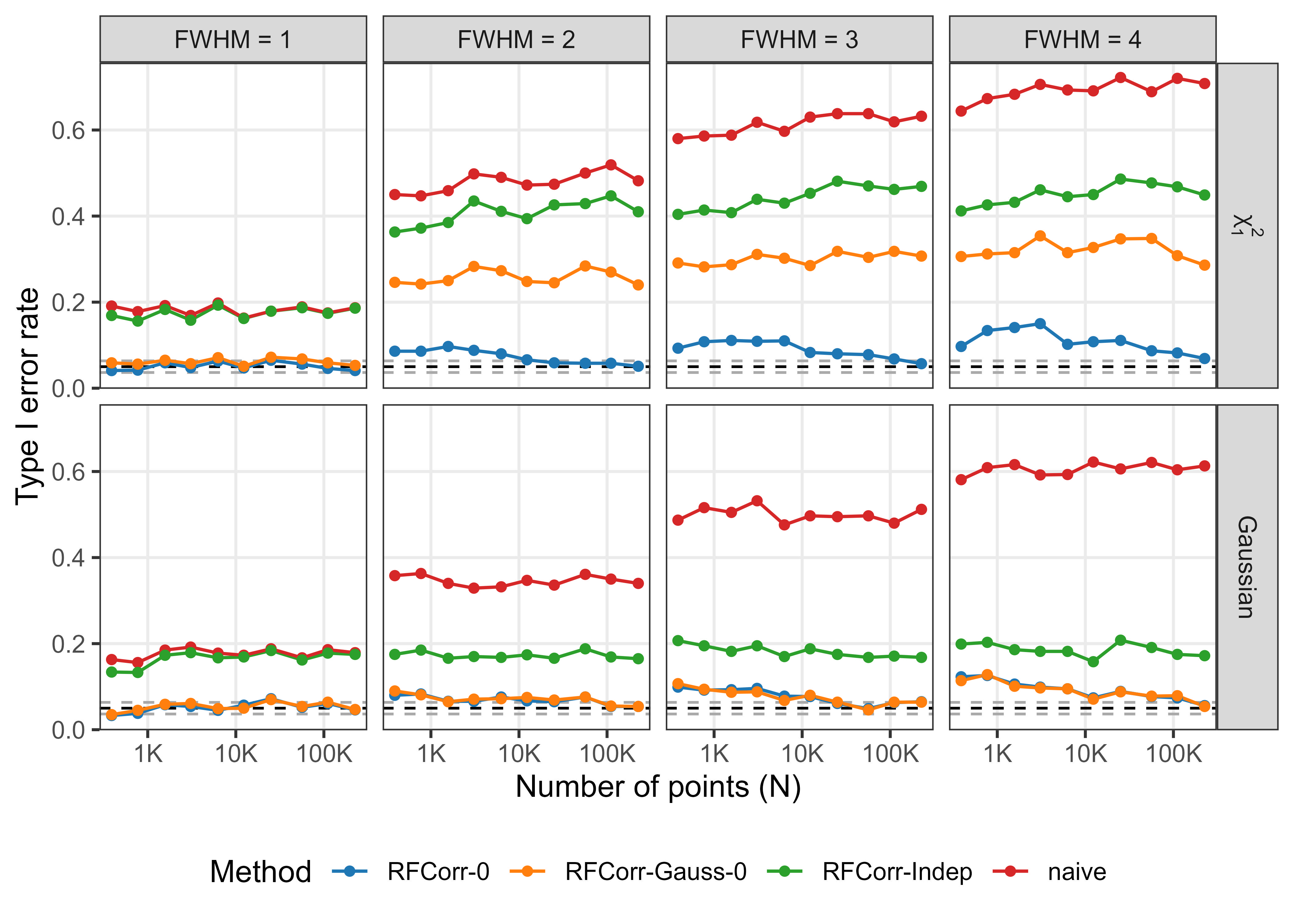}
    \caption{Rate of rejecting $H_0:\rho=0$ at level $0.05$ against the number of
    points $N$ (log scale) under the shift design of Section~\ref{sec:simu_design}, where $X$
    and $Y$ are uncorrelated but dependent, for $\chi^2_1$ fields (top) and Gaussian fields
    (bottom) in $d=2$.
    $H_0$ holds, so the rates of \rfzero{}, \rfzerogaussian{} and the naive method are Type~I error
    rates, whereas $H_0^{\text{indep}}$, tested by \rfindep{}, is false.
    Columns and dashed lines are as in Figure~\ref{fig:type1}.}
    \label{fig:type1_symmetric}
\end{figure}

\end{document}